\documentclass{amsart}
\usepackage{graphicx,indentfirst,amsfonts, amscd, amssymb, mathrsfs}

\title[A GIT Construction of Moduli Spaces of Stable Maps]{A Geometric Invariant Theory Construction of Moduli Spaces of Stable Maps}
\author[E. Baldwin]{Elizabeth Baldwin}
\author[D. Swinarski]{David Swinarski}

\newtheorem{theorem}{Theorem}[section]
\newtheorem{definition}[theorem]{Definition}

\newtheorem{claim}[theorem]{Claim}
\newtheorem{lemma}[theorem]{Lemma}
\newtheorem{amplification}[theorem]{Amplification}
\newtheorem{corollary}[theorem]{Corollary}
\newtheorem{proposition}[theorem]{Proposition}

\providecommand{\C}{\mathbb{C}}
\providecommand{\R}{\mathbb{R}}
\providecommand{\Z}{\mathbb{Z}}
\newcommand{\N}{\mathbb{N}}

\newcommand{\Spec}{\mbox{\rm Spec }}
\newcommand{\Proj}{\mbox{\rm Proj }}

\newcommand{\Pro}{\mathbf{P}}
\newcommand{\pwpr}{\Pro(W) \! \times \Pro^{r}}
\newcommand{\pw}{\Pro(W)}
\newcommand{\pr}{\Pro^r}
\newcommand{\Osheaf}{\mathcal{O}}
\newcommand{\opwopr}[2]{\mathcal{O}_{\Pro(W)}(#1)\otimes\mathcal{O}_{\Pro^r}(#2)}
\newcommand{\opw}{\Osheaf_{\Pro(W)}}
\newcommand{\opr}{\Osheaf_{\Pro^r}}
\newcommand{\Os}[1]{\mathcal{O}_{#1}}

\newcommand{\Hilb}{\mbox{\it{Hilb}}}
\newcommand{\dblq}{/\!/}
\newcommand{\Mg}{\overline{\mathcal{M}}_g}
\newcommand{\Mgn}{\overline{\mathcal{M}}_{g,n}}
\newcommand{\phat}{\hat{\rho}}
\newcommand{\mhat}{\hat{m}}
\newcommand{\Mhat}{\hat{M}}
\newcommand{\bhat}{\hat{\beta}}

\newcommand{\red}{\mbox{\tiny red}}

\newcommand{\Hom}{\mbox{\rm Hom}}
\newcommand{\Pic}{\mbox{\rm Pic}}

\newcommand{\cp}{a}
\newcommand{\acanon}{\omega^{\otimes a}}

\newcommand{\Maps}{\overline{\mathcal{M}}_{g,n}(\mathbf{P}^r,d)}
\newcommand{\curves}{\overline{\mathcal{M}}_{g,n}}
\newcommand{\pts}{x_1,\ldots ,x_n}
\newcommand{\sections}{\sigma_1,\ldots ,\sigma_n}

\newcommand{\alfa}{P(m,\mhat)}

\newcommand{\Omm}{\Omega^{m,\mhat}}
\newcommand{\barO}{\overline{\Omega}^{m,\mhat}}

\newcommand{\m}{m,\mhat,m'}

\begin{document}

\begin{abstract}   
We construct the moduli spaces of stable maps, $\Maps$, via geometric invariant theory (GIT).  This construction is only valid over $\Spec \C$, but a special case is a GIT presentation of the moduli space of stable curves of genus $g$ with $n$ marked points, $\curves$; this is valid over $\Spec \Z$.  In another paper by the first author, a small part of the argument is replaced, making the result valid in far greater generality.  Our method follows that used in the case $n=0$ by Gieseker in \cite{G}, to construct $\Mg$, though our proof that the semistable set is nonempty is entirely different.
\end{abstract}

\maketitle

\section{Introduction}
This paper gives a geometric invariant theory (GIT) construction of the Kontsevich-Manin moduli spaces of stable maps $\Maps$, for any values of $(g,n,d)$ such that smooth stable maps exist.  From this we derive a GIT construction of all such moduli spaces of stable maps $\overline{\mathcal{M}}_{g,n}(X,\beta)$, where $X$ is a projective variety and $\beta$ is a discrete invariant, understood as the homology class of the stable maps.  Although the first part of the construction closely follows Gieseker's construction in \cite{G} of the moduli spaces of stable curves $\Mg$, our proof that there exist GIT-semistable $n$-pointed maps uses an entirely different approach.

Some results of this paper are valid over $\Spec \Z$.  The GIT construction of $\Maps$ is in fact only presented in this paper over $\C$, though it can be extended to work much more generally (see \cite{e_paper_2}).  However, a special case of what we prove here is a GIT construction of the moduli spaces of $n$-pointed curves, $\Mgn$, which works over $\Spec \Z$.  A GIT construction of $\Mgn$ does not seem to have been published previously for $n>0$.

When constructing moduli spaces via GIT, one usually writes down a parameter space of the desired objects together with some extra structure, and then takes a quotient.  In our case, following the construction of \cite{FP}, this extra structure involves an embedding of the domain curve in projective space.  Given a $n$-pointed stable map $f:(C,\pts)\rightarrow\Pro^r$, we define the natural ample line bundle 
\[
\mathcal{L}:= \omega_C(x_1+\cdots +x_n)\otimes f^*\Osheaf_{\Pro^r}(c)
\]
on $C$, where $c$ is a sufficiently large positive integer, as shall be discussed in \ref{bound_for_c}.  Choose $a$ sufficiently large so that $\mathcal{L}^{a}$ is very ample.  We fix a vector space of dimension $h^0(C,\mathcal{L}^{a})$ and denote it by $W$.  A choice of isomorphism $W\cong H^0(C,\mathcal{L}^{a})$ induces an embedding $(C,\pts)\subset\Pro(W)$ and the graph of the map $f$ is a curve $(C,\pts)\subset\pwpr$.

For our parameter space with extra structure, then, we start with the Hilbert scheme $\Hilb(\pwpr)\times\prod_n(\pwpr)$, where the final factors represent the marked points.  There is a projective subscheme, $I$, the incidence subscheme where the $n$ points lie on the curve.  This is in fact the Hilbert scheme of $n$-pointed curves in $\pwpr$. We identify a locally closed subscheme $J\subset I$, corresponding to stable maps which have been embedded as described above.  This subscheme is identified by Fulton and Pandharipande; they remark (\cite{FP} Remark 2.4) that $\Maps$ is a quotient of $J$ by the action of $SL(W)$, and should be presentable via GIT, though they follow a different method.

The main theorems of this paper are stated at the beginning of Section \ref{state_thm}.  The GIT quotient $\bar{J} \dblq_L SL(W)$ is isomorphic to $\Maps$ over $\C$, for a narrow but nonempty range of linearisations $L$; if we set $r=d=0$ we obtain $\Mgn$ over $\Z$.  We prove these results for $n=0$ by generalising Gieseker's technique, and then use induction on $n$.

As alternative constructions of these spaces exist, it is natural to ask why one would go to the (considerable) trouble of constructing them via GIT, especially since the construction of this paper depends on the construction of $\Maps$ as a coarse moduli space given in \cite{FP}.  However, this paper paves the way for a construction independent of \cite{FP}, over a much more general base, laid out in \cite{e_paper_2}.  The potential stability theorem laid out here (Theorem \ref{pot_stab_without_lines}) is more generally applicable; it enables one  to construct moduli spaces of stable curves and stable maps with weighted marked points \cite{david_paper_2}.  The original motivation behind this construction was to use it as a tool for studying $\curves$, by constructing that as a GIT quotient of a subscheme of $\Maps$; see \cite{second_quotient}.  Also, once one has a space constructed via GIT, one may vary the defining linearisation to obtain birational transformations of the quotient.  Such methods may be relevant to study maps arising from the minimal model program for $\curves$ and $\Maps$, cf.\ \cite{hassett-hyeon, hyeon-lee}.

The layout of this paper is as follows: Section \ref{review} is a brief review of background material on the theory of moduli, geometric invariant theory, and stable curves and maps.  Much of the material in this section is standard.  However, we need to extend some of the theory of variation of GIT.  Thaddeus \cite{Th} and Dolgachev and Hu \cite{DH} have a beautiful picture of the way in which GIT quotients vary with linearisation.  Unfortunately, these results are only proved for projective varieties, sometimes with the extra condition of normality.  In addition the results of \cite{DH} are only given over $\C$.  We wish to make use of small parts of this theory in the setting of projective schemes over a field $k$, 
and so we make the elementary extensions necessary in Section \ref{variation_of_git}.

Let us summarise the material we shall need from the theory of variation of GIT\@.  We work in the real vector space of `virtual linearisations' generated by $G$-linearised line bundles.  One may extend the definitions of stability and semistability to virtual linearisations.  We take the cone within this space spanned by ample linearisations. Now, suppose a convex region within this ample cone has the property that no virtual linearisation in it defines a strictly semistable point.  Variation of GIT tells us that all virtual linearisations in the convex region define the same semistable set.

In Section \ref{maps_sect} we review some basic facts about stable curves and stable maps.

Our construction begins in Section \ref{geninfo}, where we define the scheme $J$ described above, and prove that there exists a family $\mathcal{C}\rightarrow J$ with the local universal property for the moduli problem of stable maps.  Here we also lay out in detail the strategy for the rest of the paper.  

Our aim is to show that for some range of virtual linearisations, GIT semistability implies GIT stability and $\bar{J}^{ss}=J$.  However, it will be sufficient for us to show that the semistable set $\bar{J}^{ss}$ is non-empty and contained in $J$.  For, by definition, elements of $J$ have finite stabiliser groups, and so all GIT semistable points will be GIT stable points if $\bar{J}^{ss}\subset J$.  An argument involving the construction of $\Maps$ from \cite{FP} allows us to conclude that if $\bar{J}^{ss}$ is a non-empty subset of $J$ then the quotient $\bar{J}\dblq SL(W)$ must be the entirety of $\Maps$.

As this argument uses the construction of \cite{FP}, which is only given over $\Spec \C$, we can only claim to have constructed $\Maps$ over $\Spec \C$.  However, this is only a short cut which we use for brevity in this paper.  An alternative argument is presented by the first author in \cite{e_paper_2}, which allows us to conclude from $\emptyset \neq \bar{J}^{ss}\subset J$ that $\bar{J}\dblq SL(W)$ is $\Maps$ over a more general base.

Within this paper, in the special case where $r=d=0$, we obtain $\curves$.  Gieseker's construction of $\Mg$ in \cite{G} works over $\Spec \Z$ (although there it is only stated to work over any algebraically closed field).  We may use induction to show that the same is true for our GIT presentation  of $\curves$.

In Section \ref{git_set-up} we describe the range of virtual linearisations and general GIT set-up to be used.  The longest part of the paper follows.  In Section \ref{ss_implies_potstab} we gradually refine our choice of virtual linearisation so that GIT semistability of a $n$-pointed map implies that it is 'potentially stable'.  The definition of precisely what is meant by this, and the corresponding theorem, can be found in Section \ref{reduced_sect}.  With this description of possible semistable curves we are able, in Section \ref{prop_and_proj}, to show that GIT semistable curves in $\bar{J}$ are indeed in $J$, at least for a carefully defined range of virtual linearisations.  All that is left is to prove non-emptiness of  the semistable set.

A further important fact may be deduced at this stage.  We have a range of virtual linearisations, which is a convex set in the vector space described above.  For this range, semistability is equivalent to stability.  It follows that the semistable set is the same for the whole of our range.  Thus non-emptiness need only be proved for one such virtual linearisation.

In Section \ref{state_thm} we complete the construction, by proving this non-emptiness.  This is done by induction on the number $n$ of marked points.  Section \ref{state_thm} is therefore divided into two: the base case and the inductive step.  In Section \ref{base_subsection}, we follow the methods of Gieseker and show that smooth maps are GIT semistable when $n=0$.  This gives us the required non-emptiness.

The inductive step follows a more novel approach, and is laid out in Section \ref{ind_step}.  Given a moduli stable map of genus $g$ with $n$ marked points, we attach an elliptic curve at the location of one of the markings,  to obtain a new stable map of genus $(g+1)$, with $(n-1)$ marked points.  Induction tells us that this has a GIT semistable model, so we have verification of the numerical criterion for GIT semistability for this map.  This implies GIT semistability of the original stable map for a virtual linearisation within the specified range.  We use the constancy of the GIT quotient for the whole of the range to deduce the result.

 As we talk here about spaces of maps from curves of differing genera and numbers of marked points, it is necessary to extend the notation $J$ to $J_{g,n,d}$ to specify which space we refer to.  The crucial result can then be summarised as
\[
\bar{J}_{g+1,n-1,d}^{ss}=J_{g+1,n-1,d}\Longrightarrow \bar{J}_{g,n,d}^{ss}=J_{g,n,d}.
\]

In the special case of curves, this induction constructs the moduli space $\overline{\mathcal{M}}_{0,n}$ for every $n\geq 3$; the base case for the induction in this case is $\overline{\mathcal{M}}_{n,0}$.

In \cite{Swin}, Swinarski gave a GIT construction of $\overline{\mathcal{M}}_{g,0}(X,\beta)$, the moduli spaces of stable maps without marked points.  Baldwin extended this in \cite{e_thesis} to marked points.  This paper brings together the results from those two theses.    Finally, we note that Parker has recently given a very different GIT construction of $\overline{\mathcal{M}}_{0,n}(\Pro^r,d)$ as a quotient of $\overline{\mathcal{M}}_{0,n}(\Pro^r\times\Pro^1,(d,1))$ in \cite{parker}.

It is a great pleasure for both authors to thank Frances Kirwan, who supervised Swinarski as an M.Sc. student and Baldwin as a D.Phil. student, for posing this problem and for all her support as we worked on it.  Elizabeth Baldwin would further like to thank Nigel Hitchin and Peter Newstead for their very useful comments on \cite{e_thesis}, and Carel Faber for encouragement and support.  She is also grateful to the Engineering and Physical Sciences Research Council for funding her doctorate and the Mittag-Leffler Institute for a very stimulating stay there.  David Swinarski would like to thank the Marshall Aid Commemoration Commission, which funded his work at Oxford through a Marshall Scholarship.  He is also grateful to Brian Conrad, Ian Morrison, Rahul Pandharipande, and Angelo Vistoli for patiently answering questions.  

\section{Background Material}\label{review}
There is a certain amount of background material which we must review.  Almost all of this section is standard, although in Section \ref{var_git} we must extend some results on variation of geometric invariant theory, to work for arbitrary schemes over a base of any characteristic.

\subsection{Moduli and quotients}

We shall take the definitions of coarse and fine moduli spaces to be standard.  However, our construction will rely on families with one particular very nice property, so that an orbit space quotient of their base is a coarse moduli space.  We recall the required definitions now.

\begin{definition}[\cite{New} page 37]
Given a moduli problem, a family $\mathcal{X}\rightarrow S$ is said to have the \emph{local universal property}\index{local universal property} if, for any other family $\mathcal{X}'\rightarrow S'$ and any $s\in S'$, there exists a neighbourhood $U$ of $s$ in $S'$ and a morphism $\phi:U\rightarrow S$ such that $\phi^*\mathcal{X}\sim \mathcal{X}'|_U$. 
\end{definition}

Suppose that we also have a group action on the base space $S$, such that orbits correspond to equivalence classes for the moduli problem.  Some sort of quotient seems a good candidate as a moduli space, but unfortunately in most cases the naive quotient will not exist as a scheme.  What we require instead is a categorical quotient:

\begin{definition}[\cite{GIT} Definition 0.5]
\label{categoricalquotient}
Let $\sigma:G\times X \rightarrow X$ be the action of an algebraic group $G$ on a scheme $X$.  A \emph{categorical quotient}\index{categorical quotient} of $X$ by $G$ is a pair $(Y,\phi)$ where $Y$ is a scheme and $\phi:X\rightarrow Y$ a morphism, such that:
\begin{enumerate}
\item[(i)]
$\phi$ is constant on the orbits of the action, i.e.\ $\phi\circ\sigma=\phi\circ p_2$ as maps from $G\times X$ to $Y$;
\item[(ii)]
given any other pair $(Z,\psi)$ for which (i) holds, there is an unique morphism $\chi:Y\rightarrow Z$ such that $\psi=\chi\circ\phi$.
\end{enumerate}

A categorical quotient $(Y,\phi)$ is an \emph{orbit space}\index{orbit space} if in addition the geometric fibres of $\phi$ are precisely the orbits of the geometric points of $X$.
\end{definition}

By definition, a categorical quotient $(Y,\phi)$ is unique up to isomorphism, and $\phi$ is a surjective morphism.  Now we see that these definitions are enough to provide coarse moduli spaces, as formalised in the following proposition:

\begin{proposition}[\cite{New} Proposition 2.13]
\label{propertiesoffamily}\label{newstead}
Suppose that the family $\mathcal{X}\rightarrow S$ has the local universal property for some moduli problem, and that the algebraic group $G$ acts on $S$, with the property that $\mathcal{X}_s \sim \mathcal{X}_t$ if and only if $G\cdot s = G\cdot t$.  Then:
\begin{enumerate}
\item[(i)]
any coarse moduli space is a categorical quotient of $S$ by $G$;
\item[(ii)]
a categorical quotient of $S$ by $G$ is a coarse moduli space if and only if it is an orbit space.\hfill $\Box$ 
\end{enumerate}
\end{proposition}


\subsection{Geometric invariant theory}\label{git}	
Geometric invariant theory (GIT) is a method to construct categorical quotients.  Details of the theory may be found in \cite{GIT}, and the results are extended over more general base in \cite{sesh_completeness} and \cite{sesh}.  More gentle introductions may be found in \cite{New} and \cite{Mukai}.  We here state the key concepts.
A quotient in geometric invariant theory depends not only on an algebraic group action on a projective scheme, but also on a linearisation of that action:
\begin{definition}
If a reductive algebraic group $G$ acts on a projective scheme $X$ and $p:L\rightarrow X$ is a line bundle, then a \emph{linearisation}\index{linearisation} of the action of $G$ with respect to $L$ is an action $\sigma:G\times L\rightarrow L$ such that:
\begin{enumerate}
\item[(i)]
for all $y\in L, \ g\in G$, $p(\sigma(g,y))=\sigma(g,p(y))$;
\item[(ii)]
for all $x\in X$, $g\in G$, the map $L_x\rightarrow L_{gx}$ defined by $y\mapsto \sigma(g,y)$ is linear.
\end{enumerate}
This defines an \emph{L-linear action} of $G$ on $X$.
\end{definition}

Line bundles together with linearisations of the action of $G$ form a group, which we denote by $\Pic^G(X)$.  An $L$-linear action of $G$ on $X$ induces an action of $G$ on the space of sections of $L^{r}$, where $r$ is any positive integer.
If $L$ is ample, then the quotient scheme we obtain is 
\begin{equation}
X\dblq_L G:=\Proj \bigoplus_{n=0}^{\infty} H^0(X,L^{\otimes n})^G.\label{X//G}\index{$X\dblq_L G$}
\end{equation}

This is a categorical quotient, but not necessarily of the whole of $X$.  The rational map $X\dashrightarrow X\dblq_L G$ is only defined at those $x\in X$ where there exists a section $s\in H^0(X,L^{\otimes n})^G$ such that $s(x)\neq 0$.  We must identify this open subscheme, and also discover to what extent the categorical quotient is an orbit space.  Accordingly we make extra definitions.  In the following, we may work with schemes defined over any universally Japanese ring, and in particular over any field or over $\Z$.
\begin{definition}[cf.\ \cite{GIT} Definition 1.7, and \cite{sesh} Proposition 7 and Remark 9]
Let $G$ be a reductive algebraic group, with an $L$-linear action on the projective scheme $X$.
\begin{enumerate}
\item[(i)]
A geometric point $x\in X$ is \emph{semi-stable}\index{GIT semistable} (with respect to $L$ and $\sigma$) if there exists $s\in H^0(X,L^{\otimes n})^G$ for some $n\geq 0$, such that $s(x)\neq 0$ and the subset $X_s$ is affine.  The open subset of $X$ whose geometric points are the semistable points is denoted $X^{ss}_{\sigma}(L)$\index{$X^{ss}(L)$}.
\item[(ii)]
A geometric point $x\in X$ is \emph{stable}\index{GIT stable} (with respect to $L$ and $\sigma$) if there exists $s\in H^0(X,L^{\otimes n})^G$ for some $n\geq 0$, such that $s(x)\neq 0$ and the subset $X_s$ is affine, the action of $G$ on $X_s$ is closed, and the stabiliser $G_x$ of $x$ is $0$-dimensional.  The open subset of $X$ whose geometric points are the stable points is denoted $X^s_{\sigma}(L)$\index{$X^s(L)$}.
\item[(iii)]
A point $x\in X$ which is semi-stable but not stable is called \emph{strictly semistable}\index{GIT strictly semistable}.
\end{enumerate}
\end{definition}
In particular, we shall use:
\begin{corollary}[\cite{GIT} page 10]\label{fin_stab_implies_stab}
If $G_x$ is finite for all $x\in X^{ss}(L)$, then $X^{ss}(L)=X^s(L)$.\hfill$\Box$
\end{corollary}
Now the main theorem of GIT is as follows:
\begin{theorem}[\cite{GIT} Theorem 1.10, \cite{sesh} Theorem 4 and Remark 9]
Let $X$ be a projective scheme, and $G$ a reductive algebraic group with an $L$-linear action on $X$.
\begin{enumerate}
\item[(i)]
A categorical quotient $(X\dblq_L G,\phi)$ of $X^{ss}(L)$ by $G$ exists.
\item[(ii)]
There is an open subset $Y^s$ of $X\dblq_L G$ such that $\phi^{-1}(Y^s)=X^s(L)$ and $(Y^s,\phi)$ is an orbit space of $X^s(L)$.
\item[(iii)]
If $L$ is ample, then $X\dblq_L G$ is  equal to $\Proj \bigoplus_{n=0}^{\infty} H^0(X,L^{\otimes n})^G$.
\hfill$\Box$
\end{enumerate}
\end{theorem}
This may all be summarised in the diagram:
\[
\begin{array}{ccccc}
X^s(L) 		&	\stackrel{open}{\subset}&	X^{ss}(L)	&	\stackrel{open}{\subset}&	X \\
\downarrow	&				&	\downarrow	&				&\\
X^s(L)/G	&	\stackrel{open}{\subset}&	X\dblq_L G.	&				&
\end{array}
\]

\label{num_crit_convention}
Stability and semi-stability are difficult to prove directly; fortunately the analysis is made much easier by utilising one parameter subgroups (1-PSs) of $G$, i.e.\ homomorphisms $\lambda:\mathbb{G}_m\rightarrow G$.  This is the so-called Hilbert-Mumford numerical criterion.  It is not used by Seshadri in \cite{sesh}; although these techniques probably do work for schemes over $\Z$, we shall only need them to apply GIT over a fixed base field, $k$.

In the following we use the conventions of Gieseker in \cite{G}, which are equivalent to but different from those of \cite{GIT}.  Note that throughout this paper, we shall use Grothendieck's convention that if $V$ is a vector space then $\Pro(V)$ is the collection of equivalence classes (under scalar action) of the nonzero elements of the dual space $V^\vee$.

Let $\lambda:\mathbb{G}_m\rightarrow G$ be a 1-PS of $G$\@.  Set $x_{\infty}:=\lim_{t\rightarrow 0}\lambda(t^{-1})\cdot x$.  The group $\lambda(\mathbb{G}_m)$ acts on the fibre $L_{x_{\infty}}$ via some character $t\mapsto t^{R}$.  Then set
\[
\mu^L(x,\lambda):=R.\index{$\mu^L(x,\lambda)$}
\]
From this perspective, one may see clearly that the map $L\mapsto \mu^L(x,\lambda)$ is a group homomorphism $\Pic^G(X)\rightarrow \Z$.

For ample line bundles we have an alternative view. Suppose $L$ is very ample, and consider $X$ as embedded in $\Pro(H^0(X,L))=:\Pro$.  We have an induced action of $\lambda(\mathbb{G}_m)$ on $H^0(X,L)$.  Pick a basis $\{e_0,\ldots ,e_{N}\}$ of $H^0(X,L)$ such that for some $r_0\leq\cdots\leq r_{N}\in\mathbb{Z}$
\[
\lambda(t)e_i=t^{r_i}e_i \mbox{ for all } t\in \mathbb{G}_m.
\]
If $\{e_0^\vee,\ldots ,e_N^\vee\}$ is the dual basis for $H^0(X,L)^\vee$, then the action of $\lambda(t)$ on $H^0(X,L)^\vee$ is given by the weights $-r_0,\ldots,-r_{N}$.

A point $x\in X$ is represented by some non-zero $\hat{x}=\sum_{i=0}^N x_i e_i^\vee\in H^0(X,L)^\vee$.  Let 
\[
R':=\min\{r_i|x_i\neq 0\}=-\max\{-r_i|x_i\neq 0\}. 
\]
Then $-R'$ is the maximum of the weights for $\hat{x}$, and so $x_{\infty}$ is represented by $\hat{x}_{\infty}:=\sum_{r_i=R'}x_i e_i^\vee$.  The fibre $L_{x_{\infty}}$ is spanned by $\{e_i(x_{\infty})|0\leq i \leq N\}$, but the nonzero part of this set is $\{e_i(x_{\infty})|r_i=R'\}$, where by definition $\lambda(\mathbb{G}_m)$ acts via the character $t\mapsto t^{R'}$.  Thus:
\[
\mu^L(x,\lambda)=R'=\min\{r_i|x_i\neq 0\}.
\]
We shall refer to the set $\{r_i:x_i\neq 0\}$ as the $\lambda$-weights of $x$.
The crucial property is that, for ample linearisations, semistability may be characterised in terms of these minimal weights:
\begin{theorem}[\cite{GIT} Theorem 2.1]
Let $k$ be a field.  Let $G$ be a reductive algebraic group scheme over $k$, with a $L$-linear action on the projective scheme $X$ (defined over $k$), where $L$ is ample.  Then:
\begin{eqnarray*}
x\in X^{ss}(L)	& \Longleftrightarrow	& \mu^L(x,\lambda)\leq 0 \mbox{ for all 1-PS $\lambda\neq 0$} \\
x\in X^{s}(L)	& \Longleftrightarrow	& \mu^L(x,\lambda)< 0 \mbox{ for all 1-PS $\lambda\neq 0$} .
\end{eqnarray*}\hfill$\Box$
\end{theorem}

\subsection{Variation of GIT}\label{var_git}\label{variation_of_git}	
The semistable set depends on the choice of linearisation of the group action.  The nature of this relationship is explored in the papers of Thaddeus \cite{Th} and Dolgachev and Hu \cite{DH} on the variation of GIT\@.  Unfortunately for us, these papers deal only with GIT quotients of projective varieties, sometimes requiring the extra condition of normality.  We wish to present $\Maps$ as a GIT quotient $\bar{J}\dblq_L SL(W)$, where the scheme $\bar{J}$ will be defined in Section \ref{defns}; as we already know that $\Maps$ is in general neither reduced nor irreducible, we cannot expect $\bar{J}$ to have either of these properties.  

It seems likely that much of the theory of variation of GIT  extends to general projective schemes.  We shall here extend the small part that we shall need to use; it is easier to prove non-emptiness of the semistable set $\bar{J}^{ss}$ if we have a certain amount of freedom in the precise choice of linearisation.  We do not need the full picture of `walls and chambers' as defined by Dolgachev and Hu in \cite{DH} and Thaddeus in \cite{Th}; developing this theory for general projective schemes would take more work, so we shall not do so here.
We shall follow the methods of \cite{DH}, though we depart from their precise conventions.

We shall assume that $X$ is a projective scheme over a field $k$.  This is more convenient than working over $\Z$ and shall be sufficient for our final results.  We further specify for all of the following that the character group $\Hom (G,\mathbb{G}_m)$ is trivial; then there is at most one $G$-linearisation for any line bundle (\cite{GIT} Proposition 1.4).  In particular this holds for $G=SL(W)$. 

We shall write $\Pic^G(X)_{\R}$ for the vector space $\Pic^G(X)\otimes\R$.  We shall refer to general elements of $\Pic^G(X)_{\R}$ as `virtual linearisations'\index{linearisation!virtual} of the group action, and denote them with a lower case $l$ to distinguish them from true linearisations, $L$.  

We shall review the construction of the crucial function $M^\bullet(x):\Pic^G(X)_{\mathbb{R}}\rightarrow \mathbb{R}$.  As the map $\mu^{\bullet}(x,\lambda):\Pic^G(X)\rightarrow\Z$ is a group homomorphism, it may be naturally extended to 
\[
\mu^{\bullet}(x,\lambda):\Pic^G(X)_{\R}\rightarrow \Z\otimes\R=\R.
\]
The numerical criterion applies for ample linearisations, so we shall be most interested in the convex cone of these:
\begin{definition}\label{ample_cone}
The {\rm ample cone}\index{ample cone} $\mathbf{A}^G(X)_{\R}$\index{$\mathbf{A}^G(X)_{\R}$} is the convex cone in $\Pic^G(X)_{\R}$ spanned by ample line bundles possessing a $G$-linearisation.  
\end{definition}

Let $T$ be a maximal torus of $G$, and let $W=N_G(T)/T$ be its Weyl group. Let $\mathcal{X}_*(G)$ be the set of non-trivial one-parameter subgroups of $G$.  Note that $\mathcal{X}_*(G)=\bigcup_{g\in G} \mathcal{X}_*(gTg^{-1})$.  If $\dim T=n$ then we can identify $\mathcal{X}_*(T)\otimes\mathbb{R}$ with $\mathbb{R}^n$.  Let $\|\cdot\|$ be a $W$-invariant Euclidean norm on $\mathbb{R}^n$.  Then for any $\lambda$ in $\mathcal{X}_*(G)$, define $\|\lambda \|:=\|g\lambda g^{-1}\|$ where $g\lambda g^{-1} \in \mathcal{X}_*(T)$.  For any 1-PS $\lambda\neq 0$, any $x\in X$ and virtual linearisation $l\in\Pic^G(X)_{\R}$, we may set
\[
\bar{\mu}^l(x,\lambda)	:= \frac{\mu^l(x,\lambda)}{\|\lambda\|}.\index{$\mu^L(x,\lambda)$!$\bar{\mu}^l(x,\lambda)$}
\]
Now our crucial function may be defined:
\begin{definition}
The function $M^{\bullet}(x):\Pic^G(X)_{\R}\rightarrow \R\cup\{\infty\}$ is defined:
\[
M^l(x)	:= \sup_{\lambda\in\mathcal{X}_*(G)} \bar{\mu}^l(x,\lambda).\index{$M^{\bullet}(x)$}
\]
\end{definition}
It is a result of Mumford that if $L$ is an ample line bundle then $M^{L}(x)$ is finite (\cite{GIT}, Proposition 2.17; recall that \cite{GIT} treats GIT over an arbitrary base field $k$). 
We observe that $M^{\bullet}(x)$ is a positively homogeneous lower convex function on $\Pic^G(X)_{\R}$.  Namely, if $\alpha\in\R_{\geq 0}$ then $M^{\alpha l}(x)=\alpha M^l(x)$; and if $l,l'\in\Pic^G(X)_{\R}$, then:
\begin{eqnarray*}
M^{l\otimes l'}(x)
	&=&\sup_{\lambda\in\mathcal{X}_*(X)}\left(\bar{\mu}^l(x,\lambda) + \bar{\mu}^{l'}(x,\lambda)\right)\\
	&\leq& \sup_{\lambda\in\mathcal{X}_*(X)}\bar{\mu}^l(x,\lambda) + \sup_{\lambda\in\mathcal{X}_*(X)}\bar{\mu}^{l'}(x,\lambda) \\
	&=&M^l(x) + M^{l'}(x).
\end{eqnarray*}
Thus in particular $M^{\bullet}(x)$ is finite-valued on the whole of $\mathbf{A}^G(X)_{\R}$.

The numerical criterion may be re-stated in terms of $M^{L}(x)$.  If $L$ is an ample linearisation of the group action, then:
\begin{eqnarray*}
x\in X^{ss}(L) 	& \Longleftrightarrow 	& M^{L}(x)\leq 0 \\
x\in X^{s}(L)	& \Longleftrightarrow	& M^{L}(x) < 0.
\end{eqnarray*}
We use $M^{\bullet}(x)$ to extend naturally the definitions of stability and semistability for virtual linearisations $l\in \mathbf{A}^G(X)_{\mathbb{R}}$.  
\begin{definition}\index{GIT semistable!with respect to virtual linearisations}\index{GIT stable!with respect to virtual linearisations}
Let $l\in \mathbf{A}^G(X)_{\R}$.  Then:
\begin{eqnarray*}
X^{ss}(l)&:=&\{x\in X : M^l(x)\leq 0\}\index{$X^{ss}(L)$!$X^{ss}(l)$} \\
X^{s}(l)&:=&\{x\in X : M^l(x)<0\}.\index{$X^s(L)$!$X^s(l)$}
\end{eqnarray*}
\end{definition}
Observe (using lower convexity for (iii)):
\begin{lemma}\label{now_obvious}Suppose $l\in\mathbf{A}^G(X)_{\R}$.
\begin{enumerate}
\item[(i)] A point $x\in X$ is semistable with respect to $l$ if and only if it is $\lambda$-semistable with respect to $l$, for all 1-PSs $\lambda$.  
\item[(ii)]A point $x\in X$ is strictly semistable with respect to $l$ precisely when $M^l(x)=0$. 
\item[(iii)]If $x$ is semistable with respect to $l_1,\ldots,l_k\in\Pic^G(X)_{\R}$ then it is semistable with respect to all virtual linearisations in the convex hull of $l_1,\ldots,l_k$.\hfill$\Box$
\end{enumerate}
\end{lemma}
Now we may prove the result that we shall need:
\begin{proposition}[cf.\ \cite{DH} Theorem 3.3.2.]\label{needed_variation} 
Suppose $\mathbf{H}\subset \mathbf{A}^G(X)_{\R}$ is a convex region satisfying $X^{ss}(l)=X^s(l)$ for all $l\in\mathbf{H}$.  It follows that $X^s(l)=X^{ss}(l)=X^{ss}(l')=X^s(l')$ for all $l,l'\in\mathbf{H}$.
\end{proposition}

\begin{proof} Let $x\in X$ be arbitrary.  It follows from the assumptions that the function $M^{\bullet}(x)$ is non-zero in $\mathbf{H}$.  Let $l,l'\in\mathbf{H}$ and let $V$ be the vector subspace of $\Pic^G(X)_{\R}$ spanned by $l$ and $l'$.  This has a basis consisting either of $l$ and $l'$, or just of $l$; use this basis to define a norm and hence a topology on $V$.  Now, since $M^{\bullet}(x)$ is positively homogeneous lower convex, the restriction
\[
M^{\bullet}(x):V\cap\mathbf{A}^G(X)_{\R}\rightarrow \R
\]
is a continuous function.  Let $L$ be the line between $l$ and $l'$.  Then $L\subset V\cap\mathbf{H}\subset V\cap\mathbf{A}^G(X)_{\R}$, so $M^{\bullet}(x)$ is non-zero and continuous on $L$.  Thus it does not change sign; either $x\in X^{s}(l'')$ for every $l''\in L$ or $x\notin X^{ss}(l'')$  for every $l''\in L$.  This holds for all $x\in X$ and so in particular $X^s(l)=X^s(l')$.
\end{proof}

{\it Remark.}  In \cite{Th} and \cite{DH} the authors use the fact that algebraically equivalent line bundles give rise to the same semistable sets, and so work in the N\'eron-Severi group:
\[
\mbox{\rm NS}^G(X):=\frac{\Pic^G(X)}{\Pic^G_0(X)}.
\]
The advantage is that in many cases, this is known to be a finitely generated abelian group, and so $NS^G(X)\otimes\R$ is a finite dimensional vector space.  However, this finite generation does not appear to have been proved in sufficient generality for our purposes.  We could show (\cite{e_thesis} Proposition 1.3.4) that the group homomorphism $\mu^{\bullet}(x,\lambda):Pic^G(X)\rightarrow\Z$ descends to $\mbox{\rm NS}^G(X)$, but as we would be left with a possibly infinite dimensional vector space in any case we have not troubled with the extra definitions and results.

\subsection{Stable curves and stable maps}\label{maps_sect}
We now turn our attention to the specific objects that we shall study.  These are the moduli spaces of stable curves and of stable maps.

The moduli space $\overline{\mathcal{M}}_{g,n}$ of Deligne--Mumford stable pointed curves is by now very well-known.  We shall not rehearse all the definitions here and instead simply cite Knudsen's work  \cite{knudsen3}; lots of background and context is given in \cite{vakil_Mgn}.  The only terminology we shall use which may not be completely standard is the following: a \emph{prestable curve} is a connected reduced projective curve whose singularities (if there are any) are nodes. 

The moduli spaces of stable maps, $\overline{\mathcal{M}}_{g,n}(X,\beta)$ parametrise isomorphism classes of certain maps from pointed nodal curves to $X$ (this will be made precise below).  They were introduced as a tool for calculating Gromov-Witten invariants, which are used in enumerative geometry and quantum cohomology.

Fix a projective scheme $X$.  The discrete invariant $\beta$ may intuitively be understood as the class of the push-forward $f_*[C]\in H_*(X;\Z)$.  In this paper we shall only complete the construction of moduli of stable maps over $\C$, and so there is no harm in taking this as the definition of $\beta$.   For the case of more general schemes $X$ see \cite{BM} Definition 2.1.  

We may define our moduli problem:  
\begin{definition}
\begin{enumerate}
\item[(i)]
A \emph{stable map}\index{maps!stable} of genus $g$, degree $d$ and homology class $\beta$ is a 
map $f:(C,\pts) \rightarrow X$ where $C$ is an $n$-pointed prestable curve of genus $g$, the homology class $f_*[C]=\beta$, and the following stability conditions are satisfied: if $C'$ is a nonsingular rational component of $C$ and $C'$ is mapped to a point by $f$, then $C'$ must have at least three {\emph special points} (either marked points or nodes); if $C'$ is a component of arithmetic genus $1$ and $C'$ is mapped to a point by $f$, then $C'$ must contain at least one special point.  
(Note that since we require the domain curves $C$ to be connected, the stability condition on genus one components is automatically satisfied except in $\overline{\mathcal{M}}_{1,0}(X,0)$, which is empty).
\item[(ii)]
A \emph{family of stable maps}\index{maps!family of stable} 
\[ 
\begin{array}{cc} 
	\mathcal{X}		 & \stackrel{f}{\rightarrow} X \\
	{\varphi}\downarrow\uparrow \sigma_i	&  \\
	S					& 
\end{array}
\] 
is a family $(\mathcal{X}\stackrel{\varphi}{\rightarrow} S,\sections)$ of pointed prestable curves together with a morphism $f: \mathcal{X} \rightarrow X$ such that $f_*[\mathcal{X}]=\beta$, and satisfying $f|_{\mathcal{X}_s}:(\mathcal{X}_s,\sigma_1(s),\ldots ,\sigma_n(s))\rightarrow X$ is a stable map for each $s\in S$.
\item[(iv)]
Two families $(\mathcal{X}\stackrel{\varphi}{\rightarrow} S,\sections,f)$ and $(\mathcal{X}'\stackrel{\varphi'}{\rightarrow} S,\sigma_1'\ldots,\sigma_n',f')$ of stable maps are \emph{equivalent} if there is an isomorphism $\tau:\mathcal{X}\cong \mathcal{X}'$ over $S$, compatible with sections, such that $f'\circ\tau=f$.
\end{enumerate}
\end{definition}

Note that $\overline{\mathcal{M}}_{g,n}(X,0)$ is simply $\curves\times X$.  In this sense the Kontsevich spaces generalise the moduli spaces of stable curves.  However, although there is an open subscheme $\mathcal{M}_{g,n}(X,\beta)$ corresponding to maps from smooth curves, in general it is not dense in the moduli space of stable maps---in general, $\overline{\mathcal{M}}_{g,n}(X,\beta)$ is reducible and has components corresponding entirely to nodal maps.  

In addition, it is very important to note that the domain of a stable map is not necessarily a stable curve!  It may have rational components with fewer than three special points (though such components cannot be collapsed by $f$).   The dualising sheaf may not be ample, even after twisting by the marked points.  We use a sheaf that provides an extra twist to all components which are not collapsed:
\[
\mathcal{L}:=\omega_C(x_1+\cdots +x_n)\otimes f^*\Osheaf_{\Pro^r}(c),
\]
where $c$ is a positive integer, whose magnitude we will discuss below:

In the special case where $X=\Pro^r$, we may fix an isomorphism $H_2(X;\Z)\cong\Z$ and denote $\beta$ by an integer $d\geq 0$.  For smooth stable maps to exist we require $2g-2+n+3d>0$; we shall only consider these cases.  

{\it Remark on the magnitude of $c$}.\label{bound_for_c}\label{note_on_c}
We require $\mathcal{L}$ to be ample on a nodal map if and only if the map is stable.  This is certainly true if $c\geq 3$, as then $\mathcal{L}$ is positive on all rational components which are not collapsed by the map or have at least three special points.  However, unless we are in the case $g=n=0$, all rational components have at least one special point, and so $c\geq 2$ will suffice for us.  If $g=n=0$ we in addition ensure that $\mathcal{L}$ is positive on irreducible curves (which now have no special points); this holds when $cd\geq 3$, so it is only in the case $(g,n,d)=(0,0,1)$ that we require $c\geq 3$.

We shall construct $\Maps$ by GIT.  A corollary is a GIT construction of $\overline{\mathcal{M}}_{g,n}(X,\beta)$.  An existing construction (not by GIT) is crucial to our proof:
\begin{theorem}[\cite{FP} Theorem 1]
Let $X$ be a projective algebraic scheme over $\mathbb{C}$, and let $\beta\in H_2(X;\Z)^+$.  There exists a projective, coarse moduli space $\overline{\mathcal{M}}_{g,n}(X,\beta)$.\hfill$\Box$
\end{theorem}
In \cite{Swin}, Swinarski gave a GIT construction of $\overline{\mathcal{M}}_{g,0}(X,\beta)$, the moduli spaces of stable maps without marked points.  Baldwin extended this in \cite{e_thesis} to marked points; this seemingly innocent extension turns out to be very difficult in GIT, because finding a linearisation with the required properties becomes much more subtle.  This paper brings together the results from those two theses.

We gather together a few more facts that have been proven about these spaces.  The most progress has been made in the case $g=0$.  If $X$ is a nonsingular convex projective variety, e.g.\ $\Pro^r$, then $\overline{\mathcal{M}}_{0,n}(X,\beta)$ is an orbifold projective variety; when non-empty it has 
the ``expected dimension''.  However, moduli spaces for stable maps of higher genera have fewer such nice properties.  Kim and Pandharipande have shown in \cite{KP} that, if $X$ is a homogeneous space $G/P$, where $P$ is a parabolic subgroup of a connected complex semisimple algebraic group $G$, then $\overline{\mathcal{M}}_{g,n}(X,\beta)$ is connected.  Little more can be said even when $X=\Pro^r$; the spaces $\overline{\mathcal{M}}_{g,n}(\Pro^{r},d)$ are in general reducible, non-reduced, and singular. 
Further, Vakil has shown in \cite{murphy} that 
every singularity of finite type over $\Z$ appears in one of the moduli spaces $\Maps$.

\section{Constructing $\Maps$: Core Definitions and Strategy} \label{geninfo}\label{construction}

We use the standard isomorphism $H_{2}(\Pro^{r}) \cong \Z$ throughout.  For any $(g,n,d)$ such that stable maps exist, we write
$\Maps$ for the moduli space of stable maps of degree $d$ from $n$-pointed genus $g$ curves into $\Pro^{r}$, as defined in Section \ref{maps_sect}.  We wish to construct this moduli space via geometric invariant theory.

The structure of the main theorem of this paper is given in this section; we shall summarise it briefly here.  We shall define a subscheme $J$ of a Hilbert scheme, such that $J$ is the base for a locally universal family of stable maps.  A group $G$ acts on $J$ such that orbits of the action correspond to isomorphism classes in the family, and hence an orbit space of $J$ by $G$, if it exists, will be precisely $\Maps$.  

The group action extends to the projective scheme $\bar{J}$, which is the closure of $J$ in the relevant Hilbert scheme.  Given any linearisation $L$ of this action, we may form a GIT quotient $\bar{J}\dblq_L G$.  Such a quotient is a categorical quotient of the semistable set $\bar{J}^{ss}(L)$, and is in addition an orbit space if all semistable points are stable.  Thus if we can show that there exists a linearisation $L$ of the action of $G$ on $\bar{J}$ such that
\[
\bar{J}^{ss}(L)=\bar{J}^{s}(L)=J,
\]
then we will have proved that $\bar{J}\dblq_L G \cong \Maps$. 

\subsection{The schemes $I$ and $J$}\label{defns}
We start by defining the scheme desired, $J=J_{g,n,d}$.  Note that all the quantities and spaces defined in the following depend on $(g,n,d)$, but that we shall only decorate them with subscripts when it is necessary to make the distinction.

{\bf Notation.}  Given a morphism $f: (C,x_1,\ldots x_n) \rightarrow \Pro^{r}$, where $C$ is nodal, write
\[
\mathcal{L}:= \omega_{C}(x_1+\cdots+x_n) \otimes f^{*}(\mathcal{O}_{\Pro^{r}}(c)).\index{$\mathcal{L}$}
\]
where $c$ is a positive integer satisfying $c\geq 2$ unless $(g,n,d)=(0,0,1)$, in which case we require $c\geq 3$, as discussed in Section \ref{maps_sect}.  Then $\mathcal{L}$ is ample on $C$ if and only if $C$ is a stable map. 
If $\cp\geq 3$\index{$a$} then $\mathcal{L}^{\cp}$ is very ample and
$h^{1}(C,\mathcal{L}^{\cp})=0$.  However, larger values of $\cp$ will be required for us to complete our GIT construction; it is shown in \cite{Schubert} that cusps are GIT stable for $\cp=3$.  We shall assume for now that $\cp\geq 5$, although it will become apparent that further refinements are needed in some cases.  Define
\[
e:=\deg(\mathcal{L}^{\cp}) = \cp (2g-2+n+cd),\index{$e$}
\]
so $h^{0}(C,\mathcal{L}^{\cp}) = e-g+1$.  We will work a lot with projective space of dimension $e-g$, so it is convenient to define
\[
N:=e-g.\index{$N$}
\]

Note a corollary of our assumptions is that $e\geq \cp g$\label{e>ag}.  If $g\geq 2$ then this follows from the inequality $2g-2+n+cd\geq 2g-2\geq g$.  If $g\leq 1$, we may see $e=\cp(2g-2+n+cd)\geq \cp \geq \cp g$.  If $g=0$ it will be more useful to estimate $e\geq \cp$.  In any case, it follows that $N\geq 4$ since $a\geq 5$.

Let $W$\index{$W$} be a vector space over $k$ of dimension $N+1$.  Then an isomorphism $W \cong H^{0}(C, \mathcal{L}^{\cp})$ induces an embedding $C \hookrightarrow \Pro(W)$ (recall that our convention is that $\Pro(V)$ is the set of equivalence classes of nonzero linear forms on $V$).  Now the graph of $f$ is an $n$-pointed nodal curve $(C,\pts)\subset\pwpr$, of bidegree $(e,d)$\@.  Its Hilbert polynomial is:
\[
P(m,\mhat):=em+d\mhat-g+1.\index{$P(m,\mhat)$}
\]

Let $\Hilb(\pwpr)$\index{$\Hilb(\pwpr)$} be the Hilbert scheme of curves in
$\pwpr$ with Hilbert polynomial $P(m,\mhat)$.  We append $n$ extra factors of $\pwpr$ to give the locations of the marked points.  Thus, given a stable map $f: (C,\pts) \rightarrow \Pro^{r}$ and a choice of isomorphism $W \cong H^{0}(C, \mathcal{L}^{\cp})$, we obtain an associated point in $\Hilb(\pwpr)\times(\pwpr)^{\times n}$.

We write $\mathcal{C}
\stackrel{\varphi}{\rightarrow} \Hilb(\pwpr)$\index{$\mathcal{C}$} for the universal family.  This may be extended:
\begin{equation}\label{uni_family}
(\varphi,1_{\prod_{i=1}^n\pwpr}):\mathcal{C}\times\prod_{i=0}^n(\pwpr)\rightarrow 
	\Hilb(\pwpr)\times\prod_{i=0}^n(\pwpr).
\end{equation}

\begin{definition}[\cite{FP} page 58]\label{i}\index{$I$}
The scheme 
\[
I\subset \Hilb(\pwpr)\times\prod_{i=1}^n(\pwpr) 
\]
is the closed incidence subscheme consisting of $(h,\pts)\subset\Hilb(\pwpr)\times(\pwpr)^{\times n}$ such that $x_1,\ldots,x_n$ lie on $\mathcal{C}_h$.
\end{definition}
We restrict the family (\ref{uni_family}) over $I$. 
This restriction is the universal family of $n$-pointed curves in $\pwpr$, possessing $n$ sections $\sections$, giving the marked points.  \index{$(\tilde{\varphi}:\tilde{\mathcal{C}}\rightarrow I,\sections)$}  Next we consider the subscheme of $I$ corresponding to $\cp$-canonically embedded stable maps:
\begin{definition}[\cite{FP} page 58]\label{j}\index{$J$}
The scheme $J\subset I$ is the locally closed subscheme consisting of those $(h,\pts)\in I$ such that:
\begin{enumerate}
\item[(i)] $(\mathcal{C}_{h},\pts)$ is \emph{prestable}, i.e.\ $\mathcal{C}_h$ is projective, connected, reduced and nodal, and $\pts$ are non-singular, distinct points on $\mathcal{C}_h$;
\item[(ii)] the projection map $\mathcal{C}_{h} \rightarrow \Pro(W)$ is a non-degenerate embedding;
\item[(iii)] $ (\mathcal{O}_{\Pro(W)}(1) \otimes \mathcal{O}_{\Pro^{r}}(1))|_{\mathcal{C}_{h}} $  and $(\omega_{\mathcal{C}_{h}}^{\cp}(\cp x_1+\cdots +\cp x_n) \otimes \mathcal{O}_{\Pro^{r}}(c \cp+1))|_{\mathcal{C}_{h}} $ are isomorphic.
\end{enumerate}
We denote by $\bar{J}$ the closure of $J$ in $I$.
\end{definition}
That $J$ is indeed a locally closed subscheme is verified in \cite{FP} page 58.  As $I$ is a projective scheme, it follows that this is also the case for $\bar{J}$.  If $(C,\pts)\subset\pwpr$ is an $n$-pointed curve, we shall say that it is represented in $J$ if $C=\mathcal{C}_h\subset\pwpr$ and $(h,\pts)\in J$.  

We restrict our universal family over $J$, denoting it $(\tilde{\varphi}^J:\tilde{\mathcal{C}}^J\rightarrow J,\sections,p_r)$\index{$(\tilde{\varphi}^J:\tilde{\mathcal{C}}^J\rightarrow J,\sections,p_r)$}, where $p_r:\tilde{\varphi}^J\rightarrow \pr$ is projection from the universal family (a subscheme of $\pwpr\times J$) to $\pr$.
The reader may easily check that this is a family of stable maps.

We extend the natural action of $SL(W)$ on $\Pro(W)$ to $\pwpr$ by defining it to be trivial on the second factor.  This
induces an action on $\Hilb(\pwpr)$, which we extend `diagonally' on $\Hilb(\pwpr)\times (\pwpr)^{\times n}$.  Note that $I$ and $J$ are invariant under this action.  Our main object of study is the $SL(W)$ action on $\bar{J}$, but we shall approach this by studying the $SL(W)$ action on $I$.

As $W$ is a vector space over $k$, we have so far only defined $I$ and $J$ as schemes over the field $k$.  By replacing each use of $\Pro(W)$ with $\Pro^N$, we may define all our schemes over $\Z$.  However, we shall continue to use $W$ for notational convenience.


\subsection{Strategy for the construction of $\Maps$}\label{strategy}

Our aim is to apply Proposition \ref{newstead} to the family $(\tilde{\varphi}^{J}:\tilde{\mathcal{C}}^{J}\rightarrow J,\sections,p_r)$.  Therefore, we check in this section that it has the local universal property and that orbits correspond to isomorphism classes.  We will first need the following lemma:


\begin{lemma} \label{locfree} 
Suppose $(\pi:\mathcal{X} \rightarrow S,\sections,f)$ is a family of stable $n$-pointed maps to $\Pro^{r}$.  Then $\pi_{*}(\omega_{\mathcal{X}/S}^{\cp}(\cp  a_1\sigma_1(S)+\cdots + \cp  a_n\sigma_n(S)) \otimes f^{*}(\mathcal{O}_{\Pro^{r}}(c\cp)))$ is a locally free $\Osheaf_S$-module, where $\sigma_i(S)$ denotes the divisor defined by the image of $\sigma_i$ in $\mathcal{X}$.
\end{lemma}

\begin{proof}The morphism $\pi:\mathcal{X}\rightarrow S$ is proper and flat, so the $\Os{\mathcal{X}}$-module $\omega_{\mathcal{X}_{s}}^{ \cp}(\sigma_1(s)+\cdots +\sigma_n(s)) \otimes f^{*}_{s}(\mathcal{O}_{\Pro^{r}}(c \cp))$ is flat over $S$ for all $s$ in $S$.  Recall that we have chosen $\cp$ so that 
\[
H^{1} \left(\mathcal{X}_{s},
\omega_{\mathcal{X}_{s}}^{\cp}(\cp \sigma_1(s)+\cdots +\cp \sigma_n(s)) \otimes f^{*}_{s}(\mathcal{O}_{\Pro^{r}}(c \cp))\right)=0
\]
for all $s \in S$.  Since $\mathcal{X}_{s}$ is a curve, all higher cohomology groups are also zero.  The hypotheses of \cite{EGA} Corollary III.7.9.9 are met, and we conclude that 
\[
R^{0}\pi_{*}(\omega_{\mathcal{X}/S}^{\cp}(\cp \sigma_1(S)+\cdots +\cp \sigma_n(S)) \otimes
f^{*}(\mathcal{O}_{\Pro^{r}}(c \cp)))
\]
is locally free.  
\end{proof}


\begin{proposition}
\label{locuni} 
\begin{enumerate}
\item[(i)]
$(\tilde{\varphi}^J:\tilde{\mathcal{C}}^J\rightarrow J,\sections,p_r)$ has the
local universal property for the moduli problem $\Maps$.
\item[(ii)]$(\mathcal{C}_{h},x_1,\ldots ,x_n,p_r|_{\mathcal{C}_h})\cong(\mathcal{C}_{h'},x'_1,\ldots ,x'_n,p_r|_{\mathcal{C}_{h'}})$ if and only if $(h,x_1,\ldots ,x_n)\in J$ and $(h',x'_1 ,\ldots x'_n)\in J$ lie in the same orbit under the action of $SL(W)$.
\end{enumerate}
\end{proposition}

\begin{proof}
(i) Suppose 
\[ 
\begin{array}{ccc} 
	\mathcal{X} & 		 \stackrel{f}{\longrightarrow} &\Pro^{r} \\
	\pi\downarrow\uparrow \sigma_i	& & \\
	S					&& 
\end{array}
\]
is an $n$-pointed family of stable maps to $\Pro^{r}$.   For any $s_{0} \in S$ we seek an open neighbourhood $V \ni s_{0}$ and a morphism $V \stackrel{\psi}{\rightarrow} J$ such that $\psi^{*}(\tilde{\mathcal{C}}^J) \sim_{\mbox{\tiny fam}} \mathcal{X}|_{V}$.  

Pick a basis for 
\[
H^{0}\left(\mathcal{X}_{s_{0}},\omega_{\mathcal{X}_{s_0}}^{\cp}(\cp  \sigma_1(s_0)+\cdots +\cp  \sigma_n(s_0)) \otimes f^{*}_{s_{0}}(\mathcal{O}_{\Pro^{r}}(c \cp))\right),
\]
and a basis for $W$.  We showed that $\pi_{*}(\omega_{\mathcal{X}/S}^{\cp}(\cp  \sigma_1(S)+\cdots \cp \sigma_n(S)) \otimes f^{*}(\mathcal{O}_{\Pro^{r}}(c \cp)))$ is locally free; it will be free on a sufficiently small neighbourhood $V \ni s_{0}$.  Then by \cite{Hart} III.12.11(b) and Lemma \ref{locfree} there is an induced basis of  
\[
H^{0} \! \left(\mathcal{X}_{s}, \omega_{\mathcal{X}_{s}}^{\cp}(\cp \sigma_1(s)+\cdots +\cp  \sigma_n(s)) \otimes f^{*}_{s}(\mathcal{O}_{\Pro^{r}}(c \cp))\right) 
\]
for each $s \in V$.  With our basis for $W$, then, this defines a map $\mathcal{X}_{s} \stackrel{\iota_{s}}{\rightarrow} \Pro(W)$ for each $s \in V$.  These fit together as a morphism $\iota:\mathcal{X}|_{V} \rightarrow \Pro(W)$.  

We have given $\mathcal{X}|_{V} \rightarrow V$ the structure of a family of $n$-pointed curves in $\pwpr$ parametrised by $V$:
\[ 
\begin{array}{ccc} 
	\mathcal{X}|_{V}& 			\stackrel{(\iota,f)}{\longrightarrow}&	\pwpr. \\
{\pi_{\mathcal{X}_V}}\downarrow\uparrow \sigma_i|_V 		& & \\
	V 							& & 
\end{array}
\]
By the universal properties of $I$, there is a unique
morphism $\psi: V \rightarrow I$ such
that $\mathcal{X}|_{V} \sim_{\mbox{\tiny fam}}
\psi^{*}(\tilde{\mathcal{C}})$.  Finally, observe that $\psi(V) \subset J$.

(ii)  If $(h,\pts)$ and $(h',x'_1,\ldots ,x'_n)$ are in the same orbit of the action of $SL(W)$ on $J$, it is immediate that as stable maps, $(\mathcal{C}_{h},x_1,\ldots ,x_n,p_r|_{\mathcal{C}_h})\cong(\mathcal{C}_{h'},x'_1,\ldots ,x'_n,p_r|_{\mathcal{C}|_{h'}})$.  Conversely, suppose we are given an isomorphism  $\tau:(\mathcal{C}_{h},x_1,\ldots ,x_n)\cong(\mathcal{C}_{h'},x'_1,\ldots ,x'_n)$; we wish to extend $\tau$ to be an automorphism of projective space $\Pro(W)$.  This is possible because the curves $\mathcal{C}_h$ and $\mathcal{C}_{h'}$ embed non-degenerately in $\Pro(W)$ via $\mathcal{L}_{\mathcal{C}_h}^{\cp}$ and $\mathcal{L}_{\mathcal{C}_{h'}}^{\cp}$ respectively.  As $\mathcal{L}_{\mathcal{C}_{h}}^{\cp}$ and $\mathcal{L}_{\mathcal{C}_{h'}}^{\cp}$ are canonically defined on isomorphic stable maps, we have an isomorphism $\mathcal{L}_{\mathcal{C}_h}^{\cp}\cong \tau^*(\mathcal{L}_{\mathcal{C}_{h'}}^{\cp})$, which induces an isomorphism on the spaces of global sections, from which we obtain an automorphism of $\Pro(W)$.
\end{proof}


\begin{corollary} \label{catquodone} If there exists a linearisation $L$ of the action of $SL(W)$ on $\bar{J}$ such that $\bar{J}^{ss}(L)=\bar{J}^{s}(L)=J$, then $\bar{J} \dblq_L SL(W)\cong\Maps$.
\end{corollary}
\begin{proof}
If $\bar{J}^{ss}(L)=J$, then $\bar{J} \dblq_L SL(W)$ is a categorical quotient of $J$, and if $\bar{J}^{ss}(L)=\bar{J}^{s}(L)$, the quotient is an orbit space.  The result follows from Proposition \ref{locuni} and Proposition \ref{newstead}.
\end{proof}
In the following GIT construction, we will first seek a linearisation $L$ such that $\bar{J}^{ss}(L)\subset J$.  This has many useful implications, which we shall explore now.  In particular, it gives us the first half of the desired equality: $\bar{J}^{ss}(L)=\bar{J}^{s}(L)$.

\begin{proposition}\label{Jss=Js}
Suppose there exists a linearisation $L$ such that $\bar{J}^{ss}(L)\subseteq J$.  Then $\bar{J}^{ss}(L)=\bar{J}^s(L)$.  
\end{proposition}
\begin{proof}
Every point of $J$ corresponds to a moduli stable map, and so has finite stabiliser.  The result follows from Corollary \ref{fin_stab_implies_stab}.
\end{proof}

This leaves us with the second required equality, that $\bar{J}^{ss}=J$.  In this paper, we shall prove it using the existing construction of $\Maps$ over $\C$, given by Fulton and Pandharipande in \cite{FP}.  This is not necessary (see \cite{e_paper_2}), but for brevity we take this shortcut for now.  
 
\begin{corollary} \label{catquo} There exists a map $j:J \rightarrow \Maps$, which is an orbit space for the $SL(W)$ action, and in particular a categorical quotient.  The morphism $j$ is universally closed.
\end{corollary}  
\begin{proof} 
$\Maps$ is a coarse moduli space (\cite{FP} Theorem 1).  By Proposition \ref{locuni} $J$ carries a local universal family, and  $SL(W)$ acts on $J$ such that orbits of the group action correspond to equivalence classes of stable maps.  The existence of the orbit space morphism $j:J\rightarrow\Maps$ follows by Proposition \ref{newstead}.   Universal closure of $j$ is a consequence of \cite{FP} Proposition 6.  
\end{proof}
Now, together with showing that $\bar{J}^{ss}(L)\subset J$, it suffices to show that the semistable set is non-empty:

\begin{theorem} \label{wishful_thinking}
Suppose that, for some linearisation $L$ of the $SL(W)$-action, $\emptyset\neq\bar{J}^{ss}(L)\subset J$.  Then, over $\C$:
\begin{enumerate}
\item[(i)] $\bar{J}\dblq_L SL(W) \cong \Maps$;
\item[(ii)] $\bar{J}^{ss}(L)=\bar{J}^s(L)=J$.
\end{enumerate}
\end{theorem}
\begin{proof}(i) Write $\bar{J}^{ss}$ for $\bar{J}^{ss}(L)$; this is an open subset of $\bar{J}$, and so $J - \bar{J}^{ss}$ is closed in $J$.  From Corollary \ref{catquo}, we have a closed morphism $j:J\rightarrow \Maps$.  Therefore $j(J - \bar{J}^{ss})$ is closed in $\Maps$.  However, $\bar{J}^{ss}$ is $SL(W)$-invariant, and $j$ is an orbit space morphism, so $j(\bar{J}^{ss})$ and $j(J - \bar{J}^{ss})$ are disjoint subsets of $\Maps$; an orbit space morphism is surjective, so these subsets make up the whole of $\Maps$.  It follows that $j(\bar{J}^{ss})=\Maps- j(J-\bar{J}^{ss})$, and is thus an open subset.
 
Now since $j$ is an orbit space and $\bar{J}^{ss}$ is $SL(W)$-invariant, it follows that the inverse image $j^{-1}j(\bar{J}^{ss})=\bar{J}^{ss}$.  The property of being a categorical quotient is local on the base, and we showed that $j(\bar{J}^{ss})$ is open in $\Maps$ so the restriction  $j|_{\bar{J}^{ss}} : \bar{J}^{ss}=j^{-1}j(\bar{J}^{ss}) \rightarrow j(\bar{J}^{ss})$ is a categorical quotient of $\bar{J}^{ss}$ for the $SL(W)$ action.  However, so is $\bar{J}\dblq SL(W)$.  A categorical quotient is unique up to isomorphism, so
\[
\bar{J} \dblq SL(W) \cong j(\bar{J}^{ss}).
\]
The scheme $\bar{J}$ projective so by construction $\bar{J} \dblq SL(W)$ is projective.  Thus $j(\bar{J}^{ss})$ is also projective, hence closed as a subset of $\Maps$.

We have shown $j(\bar{J}^{ss})$ is open and closed as a subset of $\Maps$.  It must be a union of connected components of $\Maps$.  However, Pandharipande and Kim have shown (main theorem, \cite{KP}) that $\Maps$ is connected over $\C$.  We assumed that $\bar{J}^{ss}\neq \emptyset$.  Hence
\[
\bar{J} \dblq SL(W) \cong \Maps.
\]

(ii) We showed that $\bar{J}^{ss}=\bar{J}^s$ in Proposition \ref{Jss=Js}.  We saw in the proof of (i) that $\Maps=j(\bar{J}^{ss})\sqcup j(J-\bar{J}^{ss})$ and that $j(\bar{J}^{ss})=\Maps$.  It follows that $j(J-\bar{J}^{ss})=\emptyset$, whence $\bar{J}^{ss}=J$.
\end{proof}

Moreover, proving that $\bar{J}^{ss}(L)=\bar{J}^s(L)=J$ is adequate to provide a GIT construction of the moduli spaces $\overline{\mathcal{M}}_{g,n}(X,\beta)$, where $X$ is a general projective variety.  The proof of this corollary follows the lines of \cite{FP} Lemma 8, and is given in detail as \cite{e_thesis} Corollary 3.2.8.

\begin{corollary}[\cite{e_thesis} Corollary 3.2.8, cf.\ \cite{FP} Lemma 8]\label{general_x}
Let $X$ be a projective variety defined over $\C$, with a fixed embedding to projective space $X\stackrel{\iota}{\hookrightarrow}\Pro^r$, and let $\beta\in H_2(X;\Z)^+$.  Let $g$ and $n$ be non-negative integers.  If $\beta=0$ then suppose in addition that $2g-2+n\geq 1$.  Let $\iota_*(\beta)=d\in H_2(\Pro^r;\Z)^+$.  Let $J$ be the scheme from \ref{j} corresponding to the moduli space $\Maps$.  Suppose, in the language of Theorem \ref{wishful_thinking}, that there exists a linearisation $L$ such that $\bar{J}^{ss}(L)=\bar{J}^s(L)=J$.

Then there exists a closed subscheme $J_{X,\beta}$\index{$J_{X,\beta}$} of $J$, such that 
\[
\bar{J}_{X,\beta} \dblq_{L|_{\bar{J}_{X,\beta}}} SL(W) \cong \overline{\mathcal{M}}_{g,n}(X,\beta),
\]
where $\bar{J}_{X,\beta}$ is the closure of $J_{X,\beta}$ in $\bar{J}$.\hfill$\Box$
\end{corollary}


In Section \ref{ss_implies_potstab} we shall prove that $\bar{J}^{ss}(L)\subset J$, for a suitable range of linearisations $L$.  Non-emptiness of $\bar{J}^{ss}(L)$ is dealt with in Section \ref{state_thm}.  This uses induction on the number $n$ of marked points.  For $n=0$ we show that $\bar{J}^{ss}$ is non-empty by showing that smooth maps are stable, following Gieseker.  We then apply Theorem \ref{wishful_thinking} to see that all moduli stable maps have a GIT semistable model.  However, the inductive step follows a different route, and in fact we are able to prove that $\bar{J}^{ss}=J$ directly, and then apply Corollary \ref{catquodone}, which unlike Theorem \ref{wishful_thinking} does not depend on the construction of \cite{FP}.

The generality of our proof thus is limited by the base case.  The GIT quotient $\bar{J} \dblq SL(W)$ may be defined over quite general base schemes.  Indeed, the Hilbert scheme $\Hilb(\pwpr)$ is projective over $\Spec \Z$, so our $\bar{J}$ is too, and thus we obtain a projective quotient $\bar{J}\dblq SL(W)$ over $\Spec\Z$.  However, Theorem \ref{wishful_thinking} depends on the cited results of Fulton, Kim, and Pandharipande that a projective scheme $\Maps$ exists, is a coarse moduli space for this moduli problem, and is connected; the relevant papers \cite{FP} and \cite{KP} present their results only over $\C$; we can only claim to have constructed $\Maps$ over $\Spec\C$.

In the special case of $\overline{\mathcal{M}}_{g,n}=\overline{\mathcal{M}}_{g,n}(\Pro^0,0)$ our base case is Gieseker's construction of $\Mg$, which does indeed work over $\Spec \Z$; our construction of $\curves$ thus works in this generality.

A modification of our argument is given in \cite{e_paper_2}, making our construction work independently from that of \cite{FP}, and over the more general base $\Spec \Z[p_1^{-1}\cdots p_j^{-1}]$ (where $p_1,\ldots,p_j$ are all prime numbers less than or equal to the degree $d$.
However, it is not clear whether our quotient, or an open set thereof, coarsely represents a desirable functor over $\Z$.  If the characteristic is less than the degree, a stable map may itself be an inseparable morphism.  The proof that the moduli functor is separated then fails (\cite{BM} Lemma 4.2).  Such maps can be left out if one wishes to obtain a Deligne-Mumford stack (cf. \cite{murphy} p. 2); this however will not be proper, and so our projective quotient cannot be its coarse moduli space.

Gieseker's construction of $\Mg$ begins analogously to ours.  However, the end of his argument is different from the proof of Theorem \ref{wishful_thinking} above, in two ways.  Gieseker shows directly that all smooth curves are GIT semistable, and then uses a deformation argument and semistable replacement to show that all Deligne-Mumford stable curves have $SL(W)$-semistable Hilbert points.  However, the proof that smooth curves are GIT semistable does not provide a good enough inequality to prove stability of $n$-pointed curves or maps; hence our inductive argument.  Further, not all stable maps can be smoothed (the smooth locus is not in general dense in $\Maps$), so that aspect of the argument has also needed modification. 
	

\section{The Range of Linearisations to be Used} \label{beg} \label{git_set-up}
Now we shall define linearisations of the action of $SL(W)$ on $\bar{J}$ and $I$, which were defined in \ref{j}.  Most of our analysis of GIT semistability is valid for general curves in $I$.  Accordingly it makes sense to prove results in this greater generality, so we shall define linearisations on $I$ and restrict them to $\bar{J}$.  If $L$ is a linearisation of the group action on $I$, then $\bar{J}^{ss}(L|_{\bar{J}})=\bar{J}\cap I^{ss}(L)$ by \cite{GIT} Theorem 1.19.

\subsection{The linearisations $L_{m,\mhat,m'}$ and their hull $\mathbf{H}_M(I)$}\label{lin_and_hull}

\label{mapHilbpts}  Let $(h,x_1,\ldots ,x_n) \in I \subset \Hilb(\pwpr)\times(\pwpr)^{\times n}$.  It is natural to start by defining line bundles on the separate factors, namely $\Hilb(\pwpr)$ and $n$ copies of $\pwpr$, and then take the tensor product of the pull-backs of these.  

{\bf Notation.} This notation will be used throughout the following:
\begin{eqnarray*}
Z_{m,\mhat}	&:=&H^0(\pwpr,\opwopr{m}{\mhat})\index{$Z_{m,\mhat}$} \\
P(m,\mhat)	&:=&em+d\mhat-g+1;\index{$P(m,\mhat)$}
\end{eqnarray*}
$P(m,\mhat)$ is the Hilbert polynomial of a genus $g$ curve $C\subset\pwpr$, of bidegree $(e,d)$.\\

We first define our line bundles on $Hilb(\pwpr)$.  Let $h\in\Hilb(\pwpr)$.  If $m$ and $\mhat$\index{m}\index{\mhat} are sufficiently large, then $h^{1}(\mathcal{C}_{h},
\mathcal{O}_{\Pro(W)}(m)\otimes\mathcal{O}_{\Pro^{r}}(\mhat)|_{\mathcal{C}_{h}})=0$ and the restriction map
\[ 
\label{phat} \phat^C_{m,\mhat}: H^{0}(\pwpr,\mathcal{O}_{\Pro(W)}(m)\otimes\mathcal{O}_{\Pro^{r}}(\mhat)) 
	\rightarrow H^{0}(\mathcal{C}_{h},\mathcal{O}_{\Pro(W)}(m)\otimes\mathcal{O}_{\Pro^{r}}(\mhat) |_{\mathcal{C}_{h}})\index{$\phat^C_{m,\mhat}$}
\] 
is surjective (cf.\ Grothendieck's `Uniform $m$ Lemma', \cite{HM} 1.11).  The Hilbert polynomial $P(m,\mhat)$ will be equal to $h^{0}(\mathcal{C}_{h},\mathcal{O}_{\Pro(W)}(m)\otimes\mathcal{O}_{\Pro^{r}}(\mhat) |_{\mathcal{C}_{h}})$, so that $\bigwedge^{P(m,\mhat)} \phat^C_{m,\mhat}$ gives a point of
\[
\Pro\Big(\bigwedge^{P(m,\mhat)} H^{0}(\pwpr,\mathcal{O}_{\Pro(W)}(m)\otimes\mathcal{O}_{\Pro^{r}}(\mhat))\Big)
 	=\Pro\Big(\bigwedge^{P(m,\mhat)}Z_{m,\mhat}\Big). 
\]

By the `Uniform $m$ Lemma' again, for sufficiently large $m,\mhat$, say $m,\mhat \geq m_{3},$ the Hilbert embedding\index{Hilbert embedding}
\begin{eqnarray}
\hat{H}_{m,\mhat}:\Hilb(\pwpr) 	& \rightarrow & \displaystyle{ \Pro\Big(\bigwedge^{P(m,\mhat)} Z_{m,\mhat}\Big) } \notag \\
\hat{H}_{m,\mhat}:h 		& \mapsto & \Big[\bigwedge^{P(m,\mhat)}\phat^C_{m,\mhat}\Big] \label{hilbert_point}
\end{eqnarray}
is a closed immersion (see Proposition \ref{7facts} below).  
\begin{definition}\label{L_{m,mhat}}
Let the set-up be as above, and let $m,\mhat\geq m_3$.
\index{$L_{m,\mhat}$}The line bundle $L_{m,\mhat}$ on $\Hilb(\pwpr)$ is defined to be the pull-back of the hyperplane line bundle $\Osheaf_{\Pro(\bigwedge^{P(m,\mhat)} Z_{m,\mhat})}(1)$ via the Hilbert embedding 
\[
\hat{H}_{m,\mhat}:\Hilb(\pwpr)\hookrightarrow \Pro\Big(\bigwedge^{P(m,\mhat)} Z_{m,\mhat}\Big).
\]
\end{definition}
Recall that $\overline{\mathcal{M}}_{g,n}(\Pro^0,0)=\curves$.  Whenever we write `assume $m,\mhat\geq m_3$', one should bear in mind that $\mhat$ may be set to zero in the case $r=d=0$.\label{mhat=0}

We identify $L_{m,\mhat}$ with its pull-back to $\Hilb(\pwpr)\times(\pwpr)^{\times n}$.  Now, for $i=1,\ldots,n$, let
\[
p_i:\Hilb(\pwpr)\times(\pwpr)^{\times n}\rightarrow \pwpr
\]
be projection to the $i$th such factor.  Then, for choices $m'_1,\mhat'_1,\ldots,m'_n,\mhat'_n\in\Z$, we may define $n$ line bundles on the product $\Hilb(\pwpr)\times(\pwpr)^{\times n}$:
\[
p_i^*\left(\opwopr{m'_i}{\mhat'_i}\right).
\]
The integers $\mhat'_1,\ldots,\mhat'_n$ will in fact turn out to be irrelevant to our following work.  We shall assume that they are all positive, but suppress them in notation to make things more readable.
\begin{definition}\index{$L_{m,\mhat,m'_1,\ldots,m'_n}$}\index{$L_{m,\mhat,m'}$}\label{L_many_m'_defn}
If $m,\mhat\geq m_3$ and $m'_1,\mhat'_1,\ldots,m'_n,\mhat'_n\geq 1$ then we define the very ample line bundle on $I$:
\begin{equation}\label{L_many_m'}
L_{m,\mhat,m'_1,\ldots,m'_n}:=
	\Big(L_{m,\mhat}\otimes\bigotimes_{i=1}^np_i^*\left(\opwopr{m'_i}{\mhat'_i}\right)\Big)\Big|_{I}.
\end{equation}
If $m'_1=\cdots=m'_n=m'$ then we write this as $L_{m,\mhat,m'}$.  
\end{definition}
These line bundles each possess an unique $SL(W)$-action linearising the action on $I$, which will be described in Section \ref{group_action}.   Our aim is to show that 
\[
\bar{J}\dblq_{L_{m,\mhat,m'}} SL(W) \cong\Maps, 
\]
for a suitable range of choices of $m,\mhat,m'$.  However, in order to prove that $\bar{J}^{ss}(L_{m,\mhat,m'})$ has the desired properties, we shall make use of the theory of variation of GIT (summarised in Section \ref{var_git}).  It is therefore necessary to prove results, not just for certain $L_{m,\mhat,m'}$ but for all virtual linearisations lying within the convex hull of this range in $\Pic^{SL(W)}(I)_{\R}$ or $\Pic^{SL(W)}(\bar{J})_{\R}$.  To make this precise, let $M\subset \N^{3}$ be a set such that, for every $(m,\mhat,m')\in M$, we have $m,\mhat\geq m_3$ and $m'\geq 1$.
\begin{definition}\index{$\mathbf{H}_M(I)$}\index{$\mathbf{H}_M(I)$!$\mathbf{H}_M(\bar{J})$}
We define $\mathbf{H}_M(I)$ to be the convex hull in $\Pic^{SL(W)}(I)_{\R}$ of 
\[
\{L_{m,\mhat,m'} : (m,\mhat,m')\in M\}.
\]
We define $\mathbf{H}_M(\bar{J})\subseteq \Pic^{SL(W)}(\bar{J})_{\R}$ by taking the convex hull of the restrictions of the line bundles to $\bar{J}$.
\end{definition}

As each $L_{m,\mhat,m'}$ possesses an unique lift of the action of $SL(W)$, there is an induced group action on each $l\in\mathbf{H}_M(I)$ or $\mathbf{H}_M(\bar{J})$.


\subsection{The action of $SL(W)$ for linearisations $L_{m,\mhat,m'}$}\label{group_action}	

\index{$L_{m,\mhat,m'}$!action of $SL(W)$ on}We shall describe the $SL(W)$ action on $L_{m,\mhat,m'}$.
Recall that we defined
\[
L_{m,\mhat,m'}:=\Big(L_{m,\mhat}\otimes\bigotimes_{i=1}^np_i^*\left(\opwopr{m'_i}{\mhat'_i}\right)\Big)\Big|_{I}.
\]

The linearisation of the action of $SL(W)$ on the last factors is easy to describe.  The action of $SL(W)$ on $\Os{\Pro(W)}(1)$ is induced from the natural action on $W$ (recall that $H^0(\Pro(W),\Os{\Pro(W)}(1))\cong W$).  The trivial action on $\Pro^r$ lifts to a trivial action on $\Os{\Pro^r}(1)$.  Thus we have an induced action on each $\opwopr{m'_i}{\mhat'_i}$, which may be pulled back along $p_i$ and restricted to the invariant subscheme $I$.

To describe the $SL(W)$ action on $L_{m,\mhat}$ it is easier to talk about the linear action on the projective space $\Pro \big(\bigwedge^{P(m,\mhat)}Z_{m,\mhat} \big)$.  Indeed, recall our conventions for the numerical criterion, and our definition of the function $\mu^L(x,\lambda)$ as given in Section \ref{num_crit_convention}; what we shall wish to calculate are the weights of the $SL(W)$ action on the vector space $\bigwedge^{P(m,\mhat)}Z_{m,\mhat}$. These will enable us to verify stability for a point in $\Pro \big(\bigwedge^{P(m,\mhat)}Z_{m,\mhat}\big)$, where $SL(W)$ acts with the dual action.

Fix a basis $w_{0},\ldots,w_{N}$ for $W=H^{0}(\Pro(W), \mathcal{O}_{\Pro(W)}(1))$ and a basis $f_{0},\ldots,f_{r}$ for $H^0(\Pro^r,\Os{\Pro^r}(1))$.  The group $SL(W)$ acts canonically on $H^0(\Pro(W),\Os{\Pro(W)}(1))$; the action on $H^0(\Pro^r,\Os{\Pro^r}(1))$ is trivial.

We describe the $SL(W)$ action on a basis for $\bigwedge^{P(m,\mhat)}Z_{m,\mhat}$.  Let $\hat{B}_{m,\mhat}$ be a basis for $Z_{m,\mhat}\cong H^{0}(\Pro(W), \mathcal{O}_{\Pro(W)}(m))\otimes  H^{0}(\Pro^{r}, \mathcal{O}_{\Pro^{r}}(\mhat))$ consisting of monomials of bidegree $(m,\mhat)$, where the degree $m$ part is a monomial in $w_0,\ldots ,w_N$ and the degree $\mhat$ part is a monomial in $f_0,\ldots, f_r$.
Then if $\hat{M}_{i} \in \hat{B}_{m,\mhat}$ is given by $w_{0}^{\gamma_{0}} \cdots w_{N}^{\gamma_{N}} f_{0}^{\Gamma_{0}} \cdots f_{r}^{\Gamma_{r}}$, we define $g.\hat{M}_{i} := (g.w_{0})^{\gamma_{0}} \cdots (g.w_{N})^{\gamma_{N}}  f_{0}^{\Gamma_{0}} \cdots f_{r}^{\Gamma_{r}}$.

A basis for $\bigwedge^{P(m,\mhat)}Z_{m,\mhat}$ is given by
\begin{eqnarray}
\bigwedge^{P(m,\mhat)}\hat{B}_{m,\mhat}
	:=  \{ \hat{M}_{i_{1}} \wedge \cdots \wedge \hat{M}_{i_{P(m,\mhat)}} 
	&|& 1 \leq i_{1} < \cdots < i_{P(m,\mhat)} \leq \dim Z_{m,\mhat}, \notag\\
	&& \hat{M}_{i_j}\in\hat{B}_{m,\mhat}\}. \label{wedge_basis}
\end{eqnarray}
The $SL(W)$ action on this basis is given by 
\[
 g.(\hat{M}_{i_{1}} \wedge \cdots \wedge \hat{M}_{i_{P(m,\mhat)}}) 
 	= (g.\hat{M}_{i_{1}}) \wedge \cdots \wedge (g. \hat{M}_{i_{P(m,\mhat)}}).
\]

{\bf Terminology.} We have defined virtual linearisations of the $SL(W)$ action on the scheme $I$.  We may abuse notation, and say that $(C,\pts)\subset\pwpr$ is semistable with respect to a virtual linearisation $l$ to mean that $(h,\pts)\in I^{ss}(l)$ where $(C,\pts)=(\mathcal{C}_h,\pts)$.


\subsection{The numerical criterion for $L_{m,\mhat,m'}$}				

Let $\lambda'$ be a 1-PS of $SL(W)$.  We wish to state the Hilbert-Mumford numerical criterion for our situation.  In fact, if we are careful in our analysis then we need only prove results about the semistability of points with respect to linearisations of the form $L_{m,\mhat,m'}$, as the following key lemma shows.  
\begin{lemma}\label{ignore_hull}
Fix $(h,\pts)\in I$.  Let $M$ be a range of values for $(m,\mhat,m')$.  Suppose that there exists a one-parameter subgroup $\lambda'$ of $SL(W)$, such that
\[
\mu^{L_{m,\mhat,m'}}((h,\pts),\lambda')>0 
\]
for all $(m,\mhat,m')\in M$.  Then $x$ is unstable with respect to $l$ for all $l\in\mathbf{H}_M(I)$.
\end{lemma}
\begin{proof}
Let $l\in\mathbf{H}_M(I)$.  Then $l$ is a finite combination:
\[
l=L_{m_1,\mhat_1,m'_1}^{\alpha_1}\otimes\cdots\otimes L_{m_k,\mhat_k,m'_k}^{\alpha_k},
\]
where $\alpha_1,\cdots ,\alpha_k\in\R_{\geq 0}$ satisfy $\sum \alpha_i =1$, and where $(m_i,\mhat_i,m'_i)\in M$ for $i=1,\ldots ,k$.  We know that $\mu^{L_{m_i,\mhat_i,m'_i}}((h,\pts),\lambda')>0$ for $i=1,\ldots ,k$.  The map $l'\mapsto \mu^{l'}(x,\lambda')$ is a group homomorphism $\Pic^G(I)_{\R}\rightarrow\R$, where $\R$ has its additive structure, so it follows that $\mu^{l}((h,\pts),\lambda')>0$.  Hence $M^l(h,\pts)>0$ and so $(h,\pts)$ is unstable with respect to $l$.
\end{proof}
Note the necessity of the condition that the destabilising 1-PS $\lambda'$ be the same for all $L_{m,\mhat,m'}$ such that $(m,\mhat,m')\in M$.  

Recall the definition of $L_{m,\mhat,m'}$ given in line (\ref{L_many_m'}).  From this and the functorial nature of $\mu^{\bullet}((C,\pts),\lambda)$, we see
\begin{equation}\label{num_crit_different_m'}
\mu^{L_{m,\mhat,m'}}((C,\pts),\lambda)
	= \mu^{L_{m,\mhat}}(C,\lambda)+\sum_{i=1}^n  \mu^{\Osheaf_{\Pro(W)}(1)}(x_i,\lambda)m'_i.
\end{equation}
Let us start, then, with $w_{0},\ldots,w_{N}$, a basis of $W=H^{0}(\Pro(W), \mathcal{O}_{\Pro(W)}(1))$ diagonalising the action of $\lambda'$.  There exist integers $r_{0},\ldots,r_{N}$ such that $\lambda'(t)w_{i} = t^{r_{i}}w_{i}$ for all $t \in k^*$ and $0 \leq i \leq N$.  

In the first place, following the conventions set up in Section \ref{num_crit_convention}, 
\[
\mu^{\Osheaf_{\Pro(W)}(1)}(x_i,\lambda')= \min \{r_j | w_j(x_i)\neq 0 \}.
\]
We calculate $\mu^{L_{m,\mhat}}(C,\lambda)$.
Referring to the notation of the previous subsection, if $\Mhat:=w_{0}^{\gamma_{0}} \cdots w_{N}^{\gamma_{N}} f_{0}^{\Gamma_{0}} \cdots f_{r}^{\Gamma_{r}}$,
then 
\[
\lambda'(t)\Mhat=t^{\sum \gamma_{p}r_{p}}\Mhat.
\]
Accordingly, we define the \mbox{ $\lambda'$- weight} of the monomial $\Mhat$\index{$w_{\lambda'}(\Mhat)$} to be 
\[
w_{\lambda'}(\Mhat) = \sum_{p=0}^{N} \gamma_{p}r_{p}.
\]
Let $\bigwedge^{P(m,\mhat)}\hat{B}_{m,\mhat}$ be the basis for $\bigwedge^{P(m,\mhat)}Z_{m,\mhat}$ given in line (\ref{wedge_basis}).
  Then the $\lambda'$ action on this basis is given by 
\[
\lambda'(t)(\hat{M}_{i_{1}} \wedge \cdots \wedge \hat{M}_{i_{P(m,\mhat)}}) = t^{\theta}(\hat{M}_{i_{1}} \wedge \cdots \wedge \hat{M}_{i_{P(m,\mhat)}}),
\] 
where $\theta := \sum_{j=1}^{P(m,\mhat)} w_{\lambda'}(\hat{M}_{i_{j}})$.  If we write $\hat{H}_{m,\mhat}(h)$ in the basis which is dual to $\bigwedge^{P(m,\mhat)}\hat{B}_{m,\mhat}$:
\begin{eqnarray*}
\lefteqn {\hat{H}_{m,\mhat}(h)} \\
	&&= \left[ \sum_{1 \leq j_{1} < \cdots < j_{P(m,\mhat)}}
		\phat^C_{m,\mhat} (\hat{M}_{j_{1}} \wedge \cdots \wedge \hat{M}_{j_{P(m,\mhat)}})
		\cdot(\hat{M}_{j_{1}} \wedge \cdots \wedge \hat{M}_{j_{P(m,\mhat)}}) ^{\vee} \right],
\end{eqnarray*}
so we may calculate
\[
\mu^{L_{m,\mhat}}(C,\lambda') = \mbox{min} \Big\{ \sum_{j=1}^{P(m,\mhat)} w_{\lambda'}(\hat{M}_{i_{j}})\Big\},
\]
where the minimum is taken over all sequences $1 \leq i_{j} < \cdots < i_{P(m,\mhat)}$ such that $ \phat^C_{m,\mhat} (\hat{M}_{i_{1}} \wedge \cdots \wedge \hat{M}_{i_{P(m,\mhat)}}) \neq 0$.  However, the latter is true precisely when the set $\{\phat^C_{m,\mhat}(\Mhat_{i_1}),\ldots,\phat^C_{m,\mhat}(\Mhat_{i_{P(m,\mhat)}})\}$ is a basis for $H^0(C,\opwopr{m}{\mhat}|_C)$.
 
Putting this together in (\ref{num_crit_different_m'}) we may state the numerical criterion:
$(h,\pts)$ is semistable with respect to $L_{m,\mhat,m'}$ if and only if $\mu^{L_{m,\mhat,m'}}((C,\pts), \lambda') \leq 0$ for all 1-PS $\lambda'$, where 
\[
\mu^{L_{m,\mhat,m'}}((h,\pts), \lambda') = \min \Big\{ \sum_{j=1}^{P(m,\mhat)} w_{\lambda'}(\hat{M}_{i_{j}})+\sum_{l=1}^n r_{k_l}m' \Big\},
\]
and the minimum is taken over all sequences $1 \leq i_{1} < \cdots < i_{P(m,\mhat)}$ such that $\{\phat^C_{m,\mhat}(\Mhat_{i_1}),\ldots,\phat^C_{m,\mhat}(\Mhat_{i_{P(m,\mhat)}})\}$ is a basis for $H^0(C,\opwopr{m}{\mhat}|_C)$, and all basis elements $w_{k_l}$ such that $w_{k_l}(x_l)\neq 0$.

In our applications we will often ``naturally'' write down torus actions on $W$ which highlight the geometric pathologies we wish to exclude from our quotient space.  These will usually be one-parameter subgroups of $GL(W)$ rather than $SL(W)$, as then they may be defined to act trivially on most of the space, which makes it easier to calculate their weights.  We may translate these by using ``$GL(W)$ version'' of the numerical criterion, derived as follows.

Given a 1-PS $\lambda$ of $GL(W)$, we define our ``associated 1-PS $\lambda'$ of $SL(W)$''\index{associated 1-PS of $SL(W)$}.  There is a basis $w_{0},\ldots,w_{N}$ of $H^{0}(\Pro(W), \mathcal{O}_{\Pro(W)}(1))$ so that the action of $\lambda$ is given by $\lambda(t)w_{i} = t^{r_{i}}w_{i}$ where $r_{i} \in \Z$ (the sum $\sum_{p=0}^{N} r_{p}$ is not necessarily zero).  We obtain a 1-PS $\lambda'$ of $SL(W)$ by the rule $\lambda'(t)w_{i} = t^{r_{i}'}w_{i}$ where 
\[
r_{i}' = (N+1)r_{i} - \sum_{i=0}^{N}r_{i}.  
\]
Note that now $\sum_{p=0}^N r'_p=0$.

We use this relationship to convert our numerical criterion for the $\lambda'$-action into one for the $\lambda$ action.   We define the total $\lambda$-weight of a monomial in analogy with that defined for a 1-PS of $SL(W)$.  For some suitable collection of monomials $\hat{M}_{i_1},\ldots,\hat{M}_{i_{P(m,\mhat)}}$, we see,
\begin{eqnarray*}
\lefteqn{\mu^{L_{m,\mhat,m'}}((C,\pts),\lambda')
	= \sum_{j=1}^{P(m,\mhat)} w_{\lambda'}(\hat{M}_{i_{j}}) + \sum_{l=1}^n r'_{k_l}m'} \\
	&&=(N+1) \left(\sum_{j=1}^{P(m,\mhat)} w_{\lambda}(\hat{M}_{i_{j}}) 
		+ \sum_{l=1}^n r_{k_l}m' \right) -(mP(m,\mhat) + nm') \sum_{p=0}^{N} r_{p}.
\end{eqnarray*}
Thus if $\lambda'$ is the 1-PS of $SL(W)$ arising from the 1-PS $\lambda$ of $GL(W)$, the numerical criterion may be expressed as follows: 


\begin{lemma}[cf.\ \cite{G} page 10] \label{numcritlemma}  
Let $(h,\pts)\in I$, let $\lambda$ be a 1-PS of $GL(W)$, and let $\lambda'$ be the associated 1-PS of $SL(W)$.  There exist monomials $\hat{M}_{i_{1}},\ldots,\hat{M}_{i_{P(m,\mhat)}}$ in $\hat{B}_{m,\mhat}$ such that $\{ \phat^{\mathcal{C}_h}_{m,\mhat}(\hat{M}_{i_{1}}) ,\ldots,  \phat^{\mathcal{C}_h}_{m,\mhat}(\hat{M}_{i_{P(m,\mhat)}}) \}$ is a basis of $H^{0}(\mathcal{C}_{h}, \mathcal{O}_{\Pro(W)}(m) \otimes \mathcal{O}_{\Pro^{r}}(\mhat) |_{\mathcal{C}_{h}})$, and there exist basis elements $w_{k_1},\ldots w_{k_n}$ for the $SL(W)$ action such that $w_{k_l}(x_l)\neq 0$, with
\begin{eqnarray}
\mu^{L_{\m}}((C,\pts),\lambda') 
	&=& \Big(\sum_{j=1}^{P(m,\mhat)} w_{\lambda}(\hat{M}_{i_{j}})+ \sum_{l=0}^n w_{\lambda}(w_{k_l})m'\Big)(e-g+1) \notag\\
	 &&- \left(mP(m,\mhat)+nm'\right)\sum_{p=0}^{N} w_{\lambda}(w_p);\label{num_crit}
\end{eqnarray}
moreover this choice of monomials minimises the left hand side of (\ref{num_crit}).
\hfill$\Box$
\end{lemma}
In the course of the construction we progressively place constraints on the set $M$.  In particular, for $(m,\mhat,m')\in M$, we shall be concerned with the values of the ratios $\frac{\mhat}{m}$ and $\frac{m'}{m^2}$.  It may appear surprising at first that $m'$ varies with $m^2$ and not with $m$.  Note, however, that both terms in the right hand side of (\ref{num_crit}) have terms of order $mP(m,\mhat)=em^2+dm\mhat-(g-1)m$ and terms of order $m'$, so in fact it is quite natural that $m'$ is proportional to $m^2$.


\subsection{Fundamental constants and notation}\label{constants}		

We shall now fix some notation for the whole of this paper.
The morphisms $p_{W}: \pwpr \rightarrow \Pro(W)$ and $p_{r}: \pwpr \rightarrow \Pro^{r}$ are projection onto the first and second factors respectively.  Let $C \stackrel{\iota}{\rightarrow} \pwpr$ be inclusion.   We define
\begin{eqnarray*} \label{LW}
L_{W}	&:=&	\iota^{*} p_{W}^{*} \mathcal{O}_{\Pro(W)}(1) \index{$L_W$}\\
L_{r}	&:=&	\iota^{*} p_{r}^{*}\mathcal{O}_{\Pro^{r}}(1). \index{$L_r$} 
\end{eqnarray*}
The following well-known facts are analogous to those given by Gieseker:
\begin{proposition}[cf.\ \cite{G} page 25] \label{7facts}\index{$m_{1}$}\index{$m_{2}$}\index{$m_{3}$}\index{$q_{1}$}\index{$q_{2}$}\index{$q_{3}$}\index{$\mu_{1}$}\index{$\mu_{2}$}
Let $C\subset\pwpr$ have genus $g$ and bidegree $(e,d)$.  There exist positive integers $m_{1}$, $m_{2}$, $m_{3}$, $q_{1}$, $q_{2}$, $q_{3}$, $\mu_{1}$, and $\mu_{2}$ satisfying the following properties.
\begin{enumerate}
\item[(i)] 
For all $m, \mhat > m_{1}$,  $H^{1}(C,L_{W}^{m}) =  H^{1}(C,L_{r}^{\mhat})=  H^{1}(C,L_{W}^{m} \otimes L_{r}^{\mhat})=0$.  Also $H^{1}(\bar{C},\bar{L}_{W}^{m}) =  H^{1}(\bar{C},\bar{L}_{r}^{\mhat})=  H^{1}(\bar{C},\bar{L}_{W}^{m} \otimes \bar{L}_{r}^{\mhat})=0$  and the three restriction maps 
\begin{eqnarray*} 
H^{0}(\Pro(W), \mathcal{O}_{\Pro(W)}(m)) 
	& \rightarrow & H^{0}(p_{W}(C), \mathcal{O}_{p_{W}(C)}(m)) \\
 H^{0}(\Pro^{r}, \mathcal{O}_{\Pro^{r}}(\mhat)) 
 	&  \rightarrow & H^{0}(p_{r}(C), \mathcal{O}_{p_{r}(C)}(\mhat)) \\
 H^{0}(\pwpr, \mathcal{O}_{\Pro(W)}(m)\otimes \mathcal{O}_{\Pro^{r}}(\mhat)) 
 	& \rightarrow & H^{0}(C, L_{W}^{m} \otimes L_{r}^{\mhat})
\end{eqnarray*} 
are surjective. 
\item[(ii)] 
$\mathcal{I}_{C}^{q_{1}} =0$ where $\mathcal{I}_{C}$\index{$\mathcal{I}_C$} is the sheaf of nilpotents in $\mathcal{O}_{C}$.  
\item[(iii)] 
$h^{0}(C,\mathcal{I}_{C}) \leq q_{2}$.
\item[(iv)] 
For every complete subcurve $\tilde{C}$ of $C$, $h^{0}(\tilde{C}, \mathcal{O}_{\tilde{C}}) \leq q_{3}$ and $q_{3} \geq q_{1}$.   
\item[(v)] 
$\mu_{1} > \mu_{2}$ and for every point $P \in C$ and for all integers $x \geq 0$, 
\[ 
\dim \frac{\mathcal{O}_{C,P}}{m^{x}_{C,P}} \leq \mu_{1}x + \mu_{2},
\]
where $\mathcal{O}_{C,P}$ is the local ring of $C$ at $P$ and $m_{C,P}$ is the maximal ideal of $\mathcal{O}_{C,P}$.
\item[(vi)] 
For every subcurve $\tilde{C}$ of $C$, for every point $P \in C$, and for all integers $I$ such that $m_2 \leq i <m$, the cohomology $H^{1}(\tilde{C}, \mathcal{I}_{P}^{m-i}\otimes L_{W \tilde{C}}^{m}\otimes L_{r\tilde{C}}^{\mhat})=0$, where $\mathcal{I}_{P}$ is the ideal subsheaf of $\mathcal{O}_{\tilde{C}}$ defining $P$.
\item[(vii)] 
For all integers $m, \mhat \geq m_{3}$ the map 
\begin{eqnarray*} h \!& \! \mapsto \!&\! \hat{H}_{m,\mhat}(h) \\
 \Hilb(\pwpr) \!& \! \rightarrow \! &\!
\Pro\Big(\bigwedge^{P(m,\mhat)} H^{0}(\pwpr, \mathcal{O}_{\Pro(W)}(m)\otimes \mathcal{O}_{\Pro^{r}}(\mhat))\Big)
\end{eqnarray*}
is a closed immersion.\index{Hilbert embedding}  \hfill$\Box$
\end{enumerate}
\end{proposition}

In addition, we define a constant not used by Gieseker:
\begin{eqnarray*} \index{$\bar{g}$}
\bar{g}:= \mbox{min} \{ 0, g_{\bar{Y}} 
	&| & \mbox{$\bar{Y}$ is the normalisation of a complete subcurve $Y$ contained} \\ 
	&& \mbox{in a connected fibre $\mathcal{C}_{h}$ for some $h \in \Hilb(\pwpr)$} \}.
\end{eqnarray*}
$Y$ need not be a proper subcurve.  The maximum number of irreducible components of $Y$ is $e+d$, as each must have positive degree.  Hence a lower bound for $\bar{g}$ is given by $-(e+d)+1$.  One would expect $\bar{g}$ to be negative for most $(g,n,d)$, but we have stipulated in particular that $\bar{g}\leq 0$ as this will be convenient in our calculations.


\section{GIT Semistable Maps Represented in $\bar{J}$ are Moduli Stable} \label{ss_implies_potstab}
We now embark on our proof that $\bar{J}\dblq SL(W)$ is isomorphic to $\Maps$.  Recall Proposition \ref{Jss=Js}; our first goal is to show that $\bar{J}^{ss}\subset J$. We achieve this in this section.  In Sections \ref{first_props} to \ref{reduced_sect} we work with the locus of semistable points in $I$.  Over the course of results \ref{1.0.2hat} - \ref{pot_stab_without_lines} we find a range $M$ of values for $(m,\mhat,m')$, such that if $(C,x_1,\ldots ,x_n)\subset\pwpr$ is semistable with respect to $l\in\mathbf{H}_M(I)$, then $(C,\pts)$ must be close to being a moduli stable map; this `potential stability' is defined formally in Definition \ref{potentially_stable_map}.

For this section, it is only necessary to work over a field; if we prove that equality $\bar{J}^{ss}(L)=\bar{J}^{s}(L)=J$ holds over any field $k$, then equality over $\Z$ follows (see the proof of Theorem \ref{curves}, given at the end of Section \ref{end}). 

In Section \ref{prop_and_proj} we finally restrict attention to $\bar{J}$, and greatly refine the range $M$.  Now the results of Sections \ref{first_props} - \ref{reduced_sect}, together with an application of the valuative criterion of properness, show us that if $l\in\mathbf{H}_M(\bar{J})$ then $\bar{J}^{ss}(l)\subset J$, as required.

\subsection{First properties of GIT semistable maps}\label{first_props}				
\begin{proposition}[cf.\ \cite{G} 1.0.2]\label{1.0.2hat} Let $M$ consist of integer triples  $(m,\mhat,m')$ such that $m,\mhat> m_{3}$ and $m'\geq 1$ with $m > (q_{1}-1)(e-g+1)$.  Let $l\in\mathbf{H}_M(I)$.  Suppose that $(C,\pts)\subset\pwpr$ is connected and $SL(W)$-semistable with respect to $l$.  Then $p_W(C)$ is a nondegenerate curve in $\Pro(W)$, i.e.\ $p_W(C)$ is not contained in any hyperplane in $\Pro(W)$.
\end{proposition}

\begin{proof}It is enough to prove that the composition
\[ H^{0}(\Pro(W), \mathcal{O}_{\Pro(W)}(1)) \rightarrow H^{0}(p_{W}(C)_{\red}, \mathcal{O}_{p_{W}(C)_{\red}}(1)) \rightarrow H^{0}(C_{\red}, L_{W \red})
\] 
is injective.  So suppose that it has non-trivial kernel $W_0$.
Write $N_{0} = \dim W_{0}$.  Choose a basis $w_{0},\ldots ,w_{N}$ of $H^{0}(\Pro(W), \mathcal{O}_{\Pro(W)}(1)) =: W_{1}$ relative to the filtration $0 \subset W_{0} \subset W_{1}$.  Let $\lambda$ be the 1-PS of $GL(W)$ whose action is given by
\begin{displaymath}
\begin{array}{lcrcccccc}
\lambda(t)w_{i} & = & w_{i}, & t \in \C^{*}, & 0 & \leq & i & \leq & N_{0}-1 \\
\lambda(t)w_{i} & = & tw_{i}, & t \in \C^{*}, & N_{0} & \leq & i & \leq & N,
\end{array} 
\end{displaymath}
and let $\lambda'$ be the associated 1-PS of $SL(W)$.  

We wish to show that $(C,\pts)$ is unstable with respect to any $l\in\mathbf{H}_M(I)$.  By Lemma \ref{ignore_hull} it is sufficient to show that $\mu^{L_{\m}}((C,\pts),\lambda')>0$ for all $(\m)\in M$, so pick $(\m)\in M$.

Let $\hat{B}_{m,\mhat}$ be a basis of $H^{0}(\pwpr, \mathcal{O}_{\Pro(W)}(m) \otimes \mathcal{O}_{\Pro^r}(\mhat))$ consisting of monomials of bidegree $(m,\mhat).$
Then by Lemma \ref{numcritlemma} there exist monomials  $\hat{M}_{i_{1}},\ldots ,\hat{M}_{i_{P(m,\mhat)}}$ in $\hat{B}_{m,\mhat}$ such that $\{ \phat^C_{m,\mhat}(\hat{M}_{i_{1}}),\ldots ,\phat^C_{m,\mhat}(\hat{M}_{i_{P(m,\mhat)}}) \}$ is a basis of $H^{0}(C,L_{W}^{m}\otimes L_{r}^{\mhat})$, and basis elements $w_{k_1},\ldots ,w_{k_n}$ satisfying $w_{k_l}(x_{k_l})\neq 0$, with
\begin{eqnarray}
\lefteqn {\mu^{L_{\m}}((C,\pts),\lambda') 
	= \left(\sum_{j=1}^{P(m,\mhat)}w_{\lambda}(\hat{M}_{i_{j}}) + \sum_{l=1}^{n} w_{\lambda}(w_{k_l})m'\right)(e-g+1) }\notag\\
	&&{} \hspace{2in}- \left(mP(m,\mhat)+\sum m_l'\right)\sum_{i=0}^{N}w_{\lambda}(w_{i}). \hspace{.5in} \label{num_crit_nondeg_embed} 
\end{eqnarray}
For each of the $\Mhat_{i_j}$, write $\Mhat_{i_j}=w_{0}^{\gamma_{0}}\cdots w_{N}^{\gamma_{N}}f_{0}^{\Gamma_{0}} \cdots f_{r}^{\Gamma_{r}}$.  

Recall that, if $\mathcal{I}_{C}$ denotes the ideal sheaf of nilpotent elements of $\mathcal{O}_{C}$, then the integer $q_{1}$ satisfies $\mathcal{I}_{C}^{q_{1}}=0$.  Now suppose that $\sum_{i=0}^{N_{0}-1} \gamma_{i} \geq q_{1}$.  It follows that $\phat^C_{m,\mhat}(\Mhat_{i_j})=0$, and so cannot be in a basis for $H^{0}(C,L_{W}^{m}\otimes L_{r}^{\mhat})$.  Thus $ \sum_{i=0}^{N_{0}-1} \gamma_{i} \leq q_{1}-1$.  The 1-PS $\lambda$ acts with weight $1$ on the factors $w_{N_0},\ldots,w_{N}$, and so
\[
w_{\lambda}(\hat{M}_{i_{j}}) \geq m-q_{1}+1.
\]
Our basis consists of $P(m,\mhat)$ such monomials, $\hat{M}_{i_{1}},\ldots ,\hat{M}_{i_{P(m,\mhat)}}$, so we can estimate their total weight:
\[ \sum_{j=1}^{P(m,\mhat)}w_{\lambda}(\hat{M}_{i_{j}}) \geq P(m,\mhat)(m-q_{1}+1).      
\]

We assumed that the $n$ marked points lie on the curve.  Hence if $w_{k_l}(x_l)\neq 0$ then $w_\lambda(w_{k_l})$ must be equal to $1$, so $\sum_{l=1}^n w_\lambda(w_{k_l})m' = nm'$.  Finally,
\[
\displaystyle \sum_{i=0}^{N} w_{\lambda}(w_{i}) = \dim W_{1} - \dim W_{0} = e-g+1 - \dim W_{0} \leq e-g,
\]
because $\dim W_{0} \geq 1$.    

Combining these estimates in (\ref{num_crit_nondeg_embed}), we obtain:
\begin{eqnarray*}
\lefteqn{ \mu^{L_{\m}}((C,\pts),\lambda')}\\
	&\geq& \left(P(m,\mhat)(m-q_{1}+1)+ nm'\right)(e-g+1) -\left(mP(m,\mhat)+nm'\right)(e-g) \\
	&\geq & P(m,\mhat)(m-(q_1-1)(e-g+1))	.
\end{eqnarray*}
However, $P(m,\mhat)$ is positive and by hypothesis $m> (q_{1}-1)(e-g+1)$; thus $\mu^{L_{\m}}((C,\pts),\lambda')>0$.   Recall that $\lambda'$ did not depend on the choice of $(\m)\in M$.  So $\mu^{L_{\m}}((C,\pts),\lambda')>0$ for all $(\m)\in M$.  Thus by Lemma \ref{ignore_hull}, the curve $(C,\pts)$ is not semistable with respect to any $l\in\mathbf{H}_M(I)$.
\end{proof}
Next we would like to show that no components of the curve collapse under the projection $p_W$.  We must refine the choice of virtual linearisation to obtain this result; the proof is spread over Propositions \ref{genred}--\ref{e'>0}.

Let $C_{i}$ be an irreducible component of $C$.  If the morphism $p_{W}|_{C_{i}}$ does not collapse $C_{i}$ to a point then it is finite. 
We shall find a range $M\subset \N^3$, such that if $(C,\pts)\in I^{ss}(l)$ with $l\in\mathbf{H}_M(I)$, and if the morphism $p_{W}|_{C_{i}}$ does not collapse $C_{i}$ to a point, then $p_W|_{C_{i\,\red}}$ is generically 1-1, and $C_i$ is generically reduced.

 We use the following notation.  Write $C=C'\cup Y$, where $C'$ is the union of all irreducible components of $C$ which collapse under $p_{W}$ and $Y = \overline{C-C'}$ is the union of all those that do not. 

Let $p_{W}(C)_{i}$ be the irreducible components of $p_{W}(C)$, for $i = 1,\ldots ,\ell$.  We use these to label the components of $C'$ and $Y$:
\begin{enumerate}
\item[1.] 
Let $C'_{1,1},\ldots,C'_{1,j_{1}'}$, up to $C'_{\ell,1},\ldots,C'_{\ell,j_{\ell}'}$, be the irreducible components of $C'$, labelled so that $p_{W}(C'_{i,j'}) \in p_{W}(C)_{i}$.  If there is a tie (that is, the projection of a component of $C'$ is a point lying on more than one component of $p_{W}(C))$), index it by the smallest $i$.
\item[2.] Let $Y_{1,1},\ldots,Y_{1,j_{1}}$, up to $Y_{\ell,1},\ldots,Y_{\ell,j_{\ell}}$, be the irreducible components of $Y$, so that $p_{W}(Y_{i,j}) = p_{W}(C)_{i}$.  Without loss of generality, we may assume that these are ordered in such a way that, if we set  
\[
b_{i,j}:= \deg p_{W}|_{Y_{i,j} \, \red}, 
\]
then $b_{i,j} \geq b_{i,j+1}$.
\end{enumerate}

Define:
\[\begin{array}{rclcrcl}
e_{W}	&:=&	\deg_{p_{W}(C)_{\red}} \mathcal{O}_{\Pro(W)}(1)	&&	e_{W\,i}&:=&\deg_{p_{W}(C)_{i \, \red}} \mathcal{O}_{\Pro(W)}(1) .
\end{array}\]
Since $p_{W}(C) \subset \Pro(W)$ we have $e_{W\,i} \geq 1$ for $i=1,\ldots ,\ell$.  
Recall that $L_{W}$ denotes $\iota^{*} p_{W}^{*} \mathcal{O}_{\Pro(W)}(1)$ and $L_{r}$ denotes $\iota^{*} p_{r}^{*}\mathcal{O}_{\Pro^{r}}(1)$.  By definition, the degree of $L_W$ on the components $C'_{i,j'}$ is zero, so we define:
\[\begin{array}{rclcrcl}
e_{i,j}	&:=&	\deg_{Y_{i,j \, \red}}L_{W}  		&&&&\\
d_{i,j}	&:=&\deg_{Y_{i,j \, \red}}L_{r} 	&& d'_{i,j'}&:=&\deg_{C'_{i,j'\,\red}}L_r.\\
\end{array}\]
Finally, let $\xi_{i,j}$ be the generic point of $Y_{i,j}$ and $\xi_i$ be the generic point of $p_W(C)_i$.  Write 
\[\begin{array}{rclcrcl}
k_{i,j} &:=& \mbox{length } \mathcal{O}_{Y_{i,j},\xi_{i,j}} &&	k_{i} 	&:=& \mbox{length } \mathcal{O}_{p_{W}(C)_{i},\xi_{i}}.  
\end{array}\]
Then 
\begin{eqnarray*} e &  = & \sum k_{i,j}e_{i,j} \\
e_{i,j} & = & b_{i,j}e_{W\,i} \\
e_{W} & = & \sum k_{i}e_{W\,i} .
\end{eqnarray*}
\begin{proposition}[\cite{G} 1.0.3] \label{genred} Let $M$ consist of $(\m)$ such that $m,\mhat> m_{3}$ and
\[
m > ( g - \frac{1}{2} + e(q_{1}+1) + q_{3} +\mu_{1}m_2)(e-g+1)
\]
with
\[ 
d\frac{\mhat}{m} + n\frac{m'}{m^2} < \frac{1}{4}e - \frac{5}{4}g + \frac{3}{4}.
\]  
Let $l\in\mathbf{H}_M(I)$.  Suppose that $(C,x_1,\ldots ,x_n)$ is connected and semistable with respect to $l$.  Then, in the notation explained above, the morphism $p_{W}|_{Y \, \red}$ is generically 1-1, that is, $b_{i,j} =  1$ and $j_{i} = 1$ for all $i = 1,\ldots ,\ell$.  Furthermore $Y$ is generically reduced, i.e.\ $k_{i,j}=1$ for all $i=1,\ldots ,\ell$ and $j=1,\ldots,j_i$.
\end{proposition}

{\it Remark.}
Since $e\geq ag$ and $a\geq 5$ it follows that $e-5g+3>0$, so the condition on $d\frac{\mhat}{m}+n\frac{m'}{m^2}$ may be satisfied.
\begin{proof}  Suppose not.  Then we may assume that at least one of the following is true: $j_{1} \geq 2$; or $b_{1,1} \geq 2$; or, for some $1 \leq j \leq j_{1}$, we have $k_{1,j} \geq 2$.  The first condition implies that two irreducible components of $Y$ map to the same irreducible component of $p_{W}(C)$.  The second condition implies that a component of $Y$ is a degree $b_{1,1}\geq 2$ cover of its image.  The third condition implies that the subcurve $Y$ is not generically reduced.  

Let $W_{0}$ be the kernel of the restriction map 
\[
H^{0}(\Pro(W), \mathcal{O}_{\Pro(W)}(1)) \rightarrow H^{0}(p_{W}(C)_{1 \, \red}, \mathcal{O}_{p_{W}(C)_{1 \, \red}}(1)). 
\]
{\it Step 1.}  We claim that $W_{0} \neq 0$.  To see this, suppose $W_{0}=0$.  Let $D_{1}$ be a divisor on $p_{W}(C)_{1 \, \red}$ corresponding to the invertible sheaf $\mathcal{O}_{p_{W}(C)_{1 \, \red}}(1)$ and having support in the smooth locus of $p_{W}(C)_{1 \red}$.  We have an exact sequence:
\[
0 \rightarrow \mathcal{O}_{p_{W}(C)_{1 \, \red}} \rightarrow \mathcal{O}_{p_{W}(C)_{1 \, \red}}(1) \rightarrow \mathcal{O}_{D_{1}} \rightarrow 0.  
\]
Then the long exact sequence in cohomology implies that 
\[ 
h^{0}(p_{W}(C)_{1 \, \red}, \mathcal{O}_{p_{W}(C)_{1 \, \red}}(1) ) \leq  h^{0}(p_{W}(C)_{1 \, \red}, \mathcal{O}_{D_{1}}) + h^{0}(p_{W}(C)_{1 \, \red}, \mathcal{O}_{p_{W}(C)_{1 \, \red}}).
\]
Note that $h^{0}(p_{W}(C)_{1 \, \red}, \mathcal{O}_{p_{W}(C)_{1 \, \red}})=1$, and 
\[h^{0}(p_{W}(C)_{1 \, \red}, \mathcal{O}_{D_{1}}) \leq \deg D_{1} = \deg \mathcal{O}_{p_{W}(C)_{1 \, \red}}(1) = e_{W\,1}.
\]

If $W_{0} = 0$, then 
\begin{eqnarray} e-g+1 = h^{0}(\Pro(W), \mathcal{O}_{\Pro(W)}(1))
	& \leq & h^{0}(p_{W}(C)_{1 \, \red}, \mathcal{O}_{p_{W}(C)_{1 \, \red}}(1)) \nonumber \\
 	& \leq & e_{W\,1} +1  \notag\\
\Rightarrow e-g 
	&\leq & e_{W\,1}=\frac{e_{1,1}}{b_{1,1}}. \label{me}
\end{eqnarray}
We show that this statement leads to a contradiction.  First suppose that $b_{1,1} \geq 2$ or $k_{1,1} \geq 2$.  Now we re-arrange (\ref{me}) to find:
\begin{eqnarray*}
k_{1,1}b_{1,1}(e-g) \ \leq \ k_{1,1}e_{1,1} \ = \ e - \sum_{(i,j)\neq (1,1)}k_{i,j}e_{i,j} 
	& \leq & e  \nonumber \\
\Rightarrow (k_{1,1}b_{1,1} -1)e 
	& \leq & k_{1,1}b_{1,1}g.
\end{eqnarray*}  
But our assumptions imply that $\frac{k_{1,1}b_{1,1} -1}{k_{1,1}b_{1,1}} \geq \frac{1}{2}$ and we obtain $\frac{e}{2} \leq g$, a contradiction.
On the other hand, suppose that $b_{1,1}=k_{1,1}=1$ but $j_{1} \geq 2$.  Then 
\[
e_{W\,1} =k_{1,1}b_{1,1}e_{W\,1} = e - \sum_{(i,j)\neq (1,1)}k_{i,j}b_{i,j}e_{W\,i} \leq e - k_{1,2}b_{1,2}e_{W\,1} \leq e - e_{W\,1}, 
\]
i.e.\ $e_{W\,1}\leq \frac{1}{2}e$.  Combining this with (\ref{me}), we again obtain the contradiction $\frac{e}{2} \leq g$.  

{\it Step 2.}  By Step 1 we have that $W_{0} \neq 0$, and in particular that $e_{W\,1} < (e-g)$, as it is line (\ref{me}) which leads to the contradiction.  Write $N_{0} = \dim W_{0}$.  Choose a basis $w_{0},\ldots ,w_{N}$ of  $H^{0}(\Pro(W), \mathcal{O}_{\Pro(W)}(1))$ relative to the filtration $0 \subset W_{0} \subset W_{1} = H^{0}(\Pro(W), \mathcal{O}_{\Pro(W)}(1))$.  Let $\lambda$ be the 1-PS of $GL(W)$ whose action is given by
\begin{displaymath}
\begin{array}{lcrcccccc}
\lambda(t)w_{i} & = & w_{i}, & t \in \C^{*}, & 0 & \leq & i & \leq & N_{0}-1 \\
\lambda(t)w_{i} & = & tw_{i}, & t \in \C^{*}, & N_{0} & \leq & i & \leq & N
\end{array} 
\end{displaymath}
and let $\lambda'$ be the associated 1-PS of $SL(W)$.  

Pick $(\m)\in M$ and suppose $\mu^{L_{\m}}((C,\pts),\lambda')\leq 0$.
We shall show that this leads to a contradiction.

Construct a filtration of $H^{0}(C, L_W^{m}\otimes L_{r}^{\mhat})$ as follows:
For $0 \leq p \leq m$ let $\Omm_p$ be the subspace of $H^{0}(\pwpr, \opwopr{m}{\mhat})$ spanned by monomials of weight less than or equal to $p$.  Let
\[
\barO_p:=\phat^C_{m,\mhat}(\Omm_p) \subset H^{0}(C,L_{W}^{m}\otimes L_{r}^{\mhat}).
\]
We have a filtration of $H^{0}(C, L_W^{m}\otimes L_{r}^{\mhat})$ in order of increasing weight:
\[ 
0 \subseteq \barO_0 \subseteq  \barO_1 \subseteq \cdots \subseteq \barO_m = H^{0}(C, L_W^{m}\otimes L_{r}^{\mhat}).
\]
Write $\bhat_{p} = \dim\barO_p$.

\label{B} Let  $B$ be the inverse image of $p_{W}(C)_{1}$ under $p_{W}$, i.e.\ $B = \bigcup_{j'=1}^{j_{1}'}C'_{1,j'} \cup \bigcup_{j=1}^{j_{1}}Y_{1,j}$.  Let $\tilde{C}$ be the closure of $C-B$ in $C$.  Since $C$ is connected, there is at least one closed point in $B \cap \tilde{C}$.  Choose one such point $P$.  Let 
\[
\rho^{\tilde{C},C}: H^{0}(C, L_W^{m}\otimes L_{r}^{\mhat}) \rightarrow H^{0}(\tilde{C}, L_{W \, \tilde{C}}^{m}\otimes L_{r \, \tilde{C}}^{\mhat})
\]
be the map induced by restriction.  

The following claim is analogous to one of Gieseker and may be proved using a similar argument:
\begin{claim}[cf.\ \cite{G} page 43] \label{claim} $\tilde{C} = \overline{C-B}$ can be given the structure of a closed subscheme of $C$ such that for all $0 \leq p \leq m-q_{1}$,
\[
\barO_p \cap \ker \{\rho^{\tilde{C},C}: H^{0}(C, L_{W}^{m} \otimes L_{r}^{\mhat}) \rightarrow H^{0}(\tilde{C}, L_{W\tilde{C}}^{m} \otimes L_{r\tilde{C}}^{\mhat}) \} = 0.
\]
\end{claim}

Let $\mathcal{I}_{P}$ be the ideal subsheaf of $\mathcal{O}_{\tilde{C}}$ defining the closed point $P$.  We have an exact sequence
\begin{displaymath} 0 \rightarrow \mathcal{I}_{P}^{m-p} \otimes  L_{W \, \tilde{C}}^{m}\otimes L_{r \, \tilde{C}}^{\mhat} \rightarrow  L_{W \, \tilde{C}}^{m}\otimes L_{r \, \tilde{C}}^{\mhat} \rightarrow \mathcal{O}_{\tilde{C}}/ \mathcal{I}_{P}^{m-p}\otimes  L_{W \, \tilde{C}}^{m}\otimes L_{r \, \tilde{C}}^{\mhat} \rightarrow 0.
\end{displaymath}
In cohomology this induces
\begin{eqnarray}
\lefteqn{ 0 	\rightarrow H^{0}(\tilde{C}, \mathcal{I}_{P}^{m-p} \otimes  L_{W \, \tilde{C}}^{m}\otimes L_{r \, \tilde{C}}^{\mhat}) 
		\rightarrow H^{0}(\tilde{C}, L_{W\,\tilde{C}}^{m}\otimes L_{r\,\tilde{C}}^{\mhat}) } \notag\\ 
		&&  \rightarrow  
	H^{0}(\tilde{C},\Osheaf_{\tilde{C}}/ \mathcal{I}_{P}^{m-p}\otimes L_{W\,\tilde{C}}^{m}\otimes L_{r\,\tilde{C}}^{\mhat})
		\rightarrow  
	H^{1}(\tilde{C},\mathcal{I}_{P}^{m-p}\otimes  L_{W\,\tilde{C}}^{m}\otimes L_{r \, \tilde{C}}^{\mhat}) 
		\rightarrow 0. \label{coh}
\end{eqnarray}
The following facts are analogous to those stated by Gieseker in \cite{G} page 44:
\begin{enumerate}
\item[I.] 
$ \begin{array}[t]{rcl}  
	\lefteqn{ h^{0}(\tilde{C}, L_{W\,\tilde{C}}^{m}\otimes L_{r\,\tilde{C}}^{\mhat}) 
	= \chi (L_{W\,\tilde{C}}^{m}\otimes L_{r\,\tilde{C}}^{\mhat}) 
	= \deg_{\tilde{C}}L_{W \, \tilde{C}}^{m}\otimes L_{r \, \tilde{C}}^{\mhat} + \chi(\mathcal{O}_{\tilde{C}}) }\\
		&&\leq (e-\sum k_{1,j}e_{1,j})m + (d- \sum k_{1,j}d_{1,j}- \sum k_{1,j'}d'_{1,j'})\mhat + q_3. \end{array}$
\item[II.] 
Since $\mathcal{O}_{\tilde{C}}/ \mathcal{I}_{P}^{m-p}\otimes  L_{W \, \tilde{C} }^{m}\otimes L_{r \, \tilde{C}}^{\mhat}$ has support only at the point $P\in\tilde{C}$, we estimate  $h^{0}(\tilde{C}, \mathcal{O}_{\tilde{C}}/ \mathcal{I}_{P}^{m-p}\otimes  L_{W \, \tilde{C} }^{m}\otimes L_{r \, \tilde{C}}^{\mhat}) \geq m-p.$
\item[III.] 
For $0\leq p\leq m_2-1$, Proposition \ref{7facts} says that $h^{0}(\tilde{C}, \mathcal{O}_{\tilde{C}}/ \mathcal{I}_{P}^{m-p}\otimes  L_{W \, \tilde{C}}^{m}\otimes L_{r \, \tilde{C}}^{\mhat}) \leq \mu_{1}(m-p) + \mu_{2}$, so from the long exact sequence in cohomology, $h^{1}(\tilde{C}, \mathcal{I}_{P}^{m-p} \otimes  L_{W \, \tilde{C}}^{m}\otimes L_{r \, \tilde{C}}^{\mhat}) \leq \mu_{1}(m-p) + \mu_{2}$.
\item[IV.]  
For $m_2\leq p<m$, we have $h^{1}(\tilde{C}, \mathcal{I}_{P}^{m-p} \otimes  L_{W \, \tilde{C}}^{m}\otimes L_{r \, \tilde{C}}^{\mhat}) =0$ (cf.\ Proposition \ref{7facts}).
\item[V.]  $\rho^{\tilde{C},C}_{m,\mhat}(\barO_p)\subset  H^{0}(\tilde{C}, \mathcal{I}_{P}^{m-p} \otimes  L_{W \, \tilde{C}}^{m}\otimes L_{r \, \tilde{C}}^{\mhat})\subset H^{0}(\tilde{C}, L_{W \, \tilde{C}}^{m}\otimes L_{r \, \tilde{C}}^{\mhat})$.
\end{enumerate}

If $p>m-q_1$ we may make no useful estimate, but if $0\leq p \leq m-q_{1}$ then by Claim \ref{claim} and fact V, we have 
\[
\bhat_{p} = \dim \barO_p \leq h^{0}(\tilde{C}, \mathcal{I}_{P}^{m-p} \otimes  L_{W \, \tilde{C}}^{m}\otimes L_{r \, \tilde{C}}^{\mhat}).
\]
Now the exact sequence (\ref{coh}) tells us
\begin{eqnarray*} 
\lefteqn{ h^{0}(\tilde{C}, \mathcal{I}_{P}^{m-p} \otimes  L_{W \, \tilde{C}}^{m}\otimes L_{r \, \tilde{C}}^{\mhat}) 
	= h^{0}(\tilde{C}, L_{W \, \tilde{C}}^{m}\otimes L_{r \, \tilde{C}}^{\mhat}) }\\
		&&\hspace{.5in}{}- h^{0}(\tilde{C},\mathcal{O}_{\tilde{C}}/\mathcal{I}_{P}^{m-p}\otimes L_{W\,\tilde{C}}^{m}\otimes L_{r\,\tilde{C}}^{\mhat}) 
		+ h^{1}(\tilde{C}, \mathcal{I}_{P}^{m-p} \otimes  L_{W \, \tilde{C}}^{m}\otimes L_{r \, \tilde{C}}^{\mhat}).
\end{eqnarray*}  
Thus, using the facts above, we have:
\begin{displaymath}
\bhat_{p} \leq \left\{ \begin{array}{llll} 
	\lefteqn{\left(e-\sum k_{1,j}e_{1,j}\right)m + \left(d - \sum k_{1,j}d_{1,j} - \sum k_{1,j'}d'_{1,j'}\right)\mhat} \\
			&& {}+q_{3}+p-m+\mu_{1}(m-p) + \mu_{2} \hspace{.4in}
		& \mbox{if $0 \leq p \leq m_2-1$} \\
\vspace{.1in}
	\lefteqn{\left(e-\sum k_{1,j}e_{1,j}\right)m +\left(d-\sum k_{1,j}d_{1,j}-\sum k_{1,j'}d'_{1,j'}\right)\mhat} \\
			&&{} + q_{3} + p - m 
		& \mbox{if $m_2 \leq p \leq m-q_{1}$} \\
\vspace{.1in}
	\lefteqn{ em + d\mhat -g + 1} 
		&&& \mbox{if $m-q_{1} + 1 \leq p \leq m$}. \end{array} \right.
\end{displaymath}

{\it Step 3.}  We wish to estimate $\displaystyle \sum_{j=1}^{\alfa}w_{\lambda}(\hat{M}_{i_{j}})$.  This is larger than $\displaystyle \sum_{p=1}^{m}p(\bhat_{p}-\bhat_{p-1})$, and we proceed:

\begin{eqnarray}
\lefteqn{ \sum_{p=1}^{m}p(\bhat_{p}-\bhat_{p-1}) = m\bhat_{m} - \displaystyle \sum_{p=0}^{m-1}\bhat_{p} } \notag \\
	&\geq & m(em+d\mhat+1-g) - \sum_{p=0}^{m-q_{1}} \left( (e-\sum k_{1,j}e_{1,j}) m \right. \notag\\
	&&{}+ \left.(d-\sum k_{1,j}d_{1,j}-\sum k_{1,j'}d_{1,j'})\mhat + q_{3} + p-m \right) \nonumber \\
	&&{} - \sum_{p=0}^{m_2-1} \left( \mu_{1}(m-p) + \mu_{2} \right) - \sum_{m-q_{1}+1}^{m-1}(em+d\mhat+1-g) \nonumber \\
&= & \left(\sum k_{1,j}e_{1,j}+ \frac{1}{2} \right ) m^{2} 
			+ \left (\sum k_{1,j}d_{1,j}+\sum k_{1,j'}d_{1,j'} \right )(m-q_1+1)\mhat \nonumber \\
	&&{} +  \left( \frac{3}{2}-g -q_{3}
		-\sum k_{1,j}e_{1,j}(q_{1}+1)-\mu_{1}m_2 \right)m \nonumber \\
	&&{} +  (q_{1}-1)\left(g+q_{3}-\frac{q_{1}}{2}-1\right) - \mu_{2}m_2 + \mu_{1}\frac{m_2(m_2-1)}{2} \nonumber \\
&\geq & \left( \sum k_{1,j}e_{1,j}+ \frac{1}{2} \right ) m^{2} -S_{1}m + c_2, \label{noc}
\end{eqnarray}
where 
\begin{eqnarray*}
S_{1} 	& = & g - \frac{3}{2} + \sum k_{1,j}e_{1,j}(q_{1}+1) + q_{3} +\mu_{1}m_2 \\
c_{2} 	& =& (q_{1}-1)(g+q_{3}-\frac{q_{1}}{2}-1) - \mu_{2}m_2 + \mu_{1}\frac{m_2(m_2-1)}{2}.
\end{eqnarray*}

\noindent The inequality (\ref{noc}) follows because the term $ (\sum k_{1,j}d_{1,j}+\sum k_{1,j'}d_{1,j'})\mhat( m -q_{1}+1)$ is positive, since the hypotheses imply $m > q_{1}$.  Finally, we may estimate $c_{2} \geq 0$ since $q_{3} > q_{1}$ and $\mu_{1} > \mu_{2}$ (see Proposition \ref{7facts}), to obtain
\[  
\displaystyle \sum_{j=1}^{\alfa}w_{\lambda}(\hat{M}_{i_{j}})
	\geq \left( \sum k_{1,j}e_{1,j}+ \frac{1}{2} \right ) m^{2} -S_{1}m.
\]

{\it Step 4.} We estimate the weight coming from the marked points.  We know nothing about which components each marked point lies on, so we can simply state that $\sum_{l=1}^n w_{\lambda}(w_{k_l})m' \geq 0$.  Finally, we estimate the sum of the weights:
\begin{eqnarray*}  \sum_{i=0}^{N} w_{\lambda}(w_{i})  = \dim W_{1} - \dim W_{0} & \leq & h^{0}(p_{W}(C_{1,1})_{\red}, \mathcal{O}_{p_{W}(C_{1,1})_{\red}}(1)) \\
& \leq & \deg \mathcal{O}_{p_{W}(C_{1,1}) \, \red}(1) +1 \leq   e_{W\,1} + 1.
\end{eqnarray*}  

{\it Step 5.} We combine the inequalities in Lemma \ref{numcritlemma} to obtain a contradiction as follows:  
\begin{eqnarray*}
\lefteqn {\left( (\sum k_{1,j}e_{1,j}+ \frac{1}{2})m^{2} -S_{1}m\right)(e-g+1)
	-(mP(m,\mhat)-nm')(e_{W\,1}+1) } \\
	&\hspace{2in}\leq & \mu^{L_{\m}}((C.\pts),\lambda') \leq 0 \hspace{.7in}\\
\lefteqn {\Longrightarrow
(\sum k_{1,j}e_{1,j}+\frac{1}{2})(e-g+1)m^2 - (e+\frac{1}{m})(e_{W\,1}+1)m^2 - S_{1}(e-g+1)m } \\
	 &\hspace{2in}\leq & (e_{W\,1}+1)(d\frac{\mhat}{m}+n\frac{m'}{m^2})m^2.
\end{eqnarray*}
We showed that $e_{W\,1}< e-g$, and by hypothesis $m > (g-\frac{1}{2}+e(q_{1}+1)+q_{3}+\mu_{1}m_2)(e-g+1)\geq (S_{1}+1)(e-g+1)$, so we may estimate:
\begin{eqnarray}
\frac{(\sum k_{1,j}e_{1,j}+\frac{1}{2})(e-g+1)- e(e_{W\,1}+1)-1}{(e_{W\,1}+1)} 
	& \leq & d\frac{\mhat}{m}+n\frac{m'}{m^2}.\label{now_tidy_to_show_contradiction}
\end{eqnarray}
Note that since $b_{1,1} \geq 2$ or $k_{1,1} \geq 2$ or $j_{1} \geq 2$ we have $\sum k_{1,j}b_{1,j} \geq 2$.  Thus, 
\[ 
(e_{W\,1} \sum k_{1,j}b_{1,j}+\frac{1}{2})(e-g+1) - e(e_{W\,1}+1)-1 > 0.
\]  
Furthermore the quantity 
\[ \frac{(e_{W\,1} \sum k_{1,j}b_{1,j}+\frac{1}{2})(e-g+1)- e(e_{W\,1}+1)-1}{(e_{W\,1}+1)}
\] 
is minimised when $e_{W\,1}$ takes its smallest value, that is, when $e_{W\,1}=1$.  So
\begin{eqnarray*}
\frac{(e_{W\,1} \sum k_{1,j}b_{1,j}+\frac{1}{2})(e-g+1) - e(e_{W\,1}+1)-1}{(e_{W\,1}+1)}
	&\geq& \frac{\frac{5}{2}(e-g+1)-2e-1}{2} \\
	&=& \frac{1}{4} e - \frac{5}{4}g +\frac{3}{4}.
\end{eqnarray*}

But by hypothesis $d\frac{\mhat}{m} + n\frac{m'}{m^2} < \frac{1}{4}e - \frac{5}{4}g+\frac{3}{4}$; combining this result with (\ref{now_tidy_to_show_contradiction}) gives a contradiction.  This implies that $\mu^{L_{\m}}((C,\pts),\lambda')>0$, and this holds for all $(\m)\in M$.  It follows by Lemma \ref{ignore_hull} that $(C,\pts)$ is not semistable with respect to $l$ for any $l\in \mathbf{H}_M(I)$.  This completes the proof of Proposition \ref{genred}.
\end{proof}

In the following propositions, we often estimate $\mu^{L_{\m}}((C,\pts))$ using a similar techniques each time.  It is thus efficient to prove two results useful for these computations in advance.

Firstly, if $C$ is a general curve, we have an inclusion $i:C_{\red}\hookrightarrow C$.  The reduced curve $C_{\red}$ has normalisation\index{curve!normalisation of} $\pi':\bar{C}_{\red}\rightarrow C_{\red}$.  Following Gieseker in \cite{G} page 22, we define the normalisation $\pi:\bar{C}\rightarrow C$ by letting $\bar{C}:=\bar{C}_{\red}$ and $\pi:=i\circ\pi'$.  Then, whatever the properties of $C$, the curve $\bar{C}$ is smooth and integral.  With these conventions, we may show:


\begin{claim}[cf.\ \cite{G} p.52]\label{I,II,III}\label{3facts}
Let $C$ be a generically reduced curve over $k$; we do not assume it has genus $g$.  Let $\pi:\bar{C}\rightarrow C$ be the normalisation morphism, and let $\mathcal{I}_C$ be the sheaf of nilpotents.
\begin{enumerate}
\item[I.]
Suppose that $C\subset \pwpr$.  Define $L_{W\,C}$ and $L_{r\,C}$ as in Section \ref{LW}, and let $\bar{L}_{W\,\bar{C}}:=\pi^*L_{W\,C}$ and $\bar{L}_{r\,\bar{C}}:=\pi^*L_{r\,C}$.
Let
\[
\pi_{m,\mhat *}:H^0(C,L_{W\,C}^m\otimes L_{r\,C}^{\mhat})
	\rightarrow H^0(\bar{C},\bar{L}_{W\,\bar{C}}^m\otimes\bar{L}_{r\,\bar{C}}^{\mhat})
\]
be the induced morphism.  Then
\[
\dim \ker \pi_{m,\mhat *} = h^0(C,\mathcal{I}_C).
\]
\item[II.]  Suppose that $C$ is reduced.
Let $D$ be an effective divisor on $C$, and let $M$ be an invertible sheaf on $C$ such that $H^1(C,M)=0$.  Then $h^1(C,M(-D))\leq \deg D$.
\item[III.] Suppose that $C$ is integral and smooth, with genus $g_{C}$.  Let $M$ be an invertible sheaf on $C$ with $\deg M \geq 2g_{C} -1$.  Then $H^1(C,M)=0$.\hfill$\Box$
\end{enumerate}
\end{claim}
Secondly, suppose we have a 1-PS $\lambda$ of $GL(W)$ and the associated 1-PS $\lambda'$ of $SL(W)$.  We shall find a lower bound for $\mu^{L_{\m}}((C,\pts),\lambda')$ by filtering the vector space $H^0(C,L_W^m\otimes L_r^{\mhat})$ according to the weight with which $\lambda$ acts.  As our filtration is constructed in the same way every time, we describe it now.

{\bf Notation.} Suppose we have a 1-PS $\lambda$ of $GL(W)$, acting with weights $r_0\leq\cdots\leq r_N$ with respect to the basis $w_0,\ldots,w_N$ for $W$.  We also assume $r_0 \geq 0$; this will be the case in all our applications.  Let $R$ be a positive integer such that $\sum_{i=0}^N r_i\leq R$.  For $0 \leq p \leq m$ let $\Omm_p$ be the subspace of $H^{0}(\pwpr, \opwopr{m}{\mhat})$ spanned by monomials of weight less than or equal to $p$.  Let
\[
\barO_p:=\phat^C_{m,\mhat}(\barO_p) \subset H^{0}(C,L_{W}^{m}\otimes L_{r}^{\mhat}).
\]
We have a filtration of $H^{0}(C, L_W^{m}\otimes L_{r}^{\mhat})$ in order of increasing weight:
\begin{equation}
0 \subseteq \barO_0 \subseteq  \barO_1 \subseteq \cdots \subseteq \barO_m = H^{0}(C, L_W^{m}\otimes L_{r}^{\mhat}).\label{template_filtration}
\end{equation}
Write $\bhat_{p} = \dim\barO_p$.  With these conventions, we show how to estimate the minimal weight, $\mu^{L_{\m}}((C,\pts),\lambda')$:

\begin{lemma}\label{main_calculation}
In the set-up described above, suppose that
\[
\bhat_p \leq (e-\alpha)m +(d-\beta)\mhat + \gamma p + \epsilon_p
\]
where $\alpha,\beta,\gamma,\epsilon_p$ are constants.  Set 
\[
\epsilon:=\frac{1}{m}\sum_{p=0}^{r_Nm-1}\epsilon_p.  
\]
Suppose 
\[
\sum_{j=1}^n w_{\lambda}(w_{i_j})m'=\delta m', 
\]
where $w_{i_1},\ldots,w_{i_n}$ are the basis elements of minimal weight satisfying $w_{i_j}(x_j)\neq 0$. Then
\begin{eqnarray}
\lefteqn {\mu^{L_{\m}}((C,\pts),\lambda')
	\geq \left( (r_N\alpha-r_N^2\frac{\gamma}{2})(e-g+1)-Re\right)m^2 } \notag\\
	 &\hspace{1in}&	{}+ \left( r_N\beta(e-g+1)-Rd\right)m\mhat 
			{}+ \left(\delta(e-g+1) -Rn\right)m'  \notag\\
	&\hspace{1in}&	{}-\left((r_N(g-1)-\frac{r_N\gamma}{2}+\epsilon)(e-g+1)+R\right)m. \hspace{.3in}\label{formula_for_mu}
\end{eqnarray}
\end{lemma}
\begin{proof}
Suppose we have any monomials $\hat{M}_{i_{1}},\ldots ,\hat{M}_{i_{P(m,\mhat)}}$ in $\hat{B}_{m,\mhat}$ such that the set $\{ \phat^C_{m,\mhat}(\hat{M}_{i_{1}}),\ldots ,\phat^C_{m,\mhat}(\hat{M}_{i_{P(m,\mhat)}}) \}$ is a basis of $H^{0}(C,L_{W}^{m}\otimes L_{r}^{\mhat})$.  As our filtration is in order of increasing weight, a lower bound for $\sum_{j=1}^{P(m,\mhat)}w_{\lambda}(\hat{M}_{i_{j}})$ is  $\sum_{p=1}^{r_Nm}p(\bhat_{p}-\bhat_{p-1})$.  We calculate: 
\begin{eqnarray*}
\lefteqn{\displaystyle \sum_{p=1}^{r_N m}p(\bhat_{p}-\bhat_{p-1}) = r_N m\bhat_{r_N m} - \displaystyle \sum_{p=0}^{r_N m-1} \bhat_{p} } \\
	&\geq&  r_Nm(em+d\mhat-g+1)-\sum_{p=0}^{r_Nm-1}((e-\alpha)m +(d-\beta)\mhat + \gamma p + \epsilon_p) \\
	&=&  (r_N\alpha-\frac{r_N^2\gamma}{2}) m^2 + r_N\beta m\mhat - (r_N(g-1)-\frac{r_N\gamma}{2}+\epsilon)m,
\end{eqnarray*}
where $\epsilon:=\frac{1}{m}\sum_{p=0}^{r_Nm-1}\epsilon_p$.
Let $\lambda'$ be the associated 1-PS of $SL(W)$.  Thus, using Lemma \ref{numcritlemma}, we calculate:
\begin{eqnarray*}
\lefteqn{\mu^{L_{\m}}((C,\pts),\lambda') }\\
	&\geq & \left((r_N\alpha-r_N^2\frac{\gamma}{2}) m^2 + r_N\beta m\mhat - (r_N(g-1)-\frac{r_N\gamma}{2}+\epsilon)m + \delta m'\right)(e-g+1) \\
	&&	{}-(m(em+d\mhat-g+1)+nm')\sum_{i=0}^N r_i \\
	&=& \left( (r_N\alpha-r_N^2\frac{\gamma}{2})(e-g+1)-Re\right)m^2 \\
	&&	{}+ \left( r_N\beta(e-g+1)-Rd\right)m\mhat + \left(\delta(e-g+1) -Rn\right)m'\\
	&&	{}-\left((r_N(g-1)-\frac{r_N\gamma}{2}+\epsilon)(e-g+1)+R\right)m,
\end{eqnarray*}
where we have used the bounds $0\leq \sum_{i=0}^N r_i \leq R$ to estimate appropriately, according to the sign of each term.
\end{proof}
{\it Remark.}  In general, we shall assume that $m$ is very large, that $\mhat$ is proportional to $m$ and that $m'$ is proportional to $m^2$.  

Next we derive the inequality (\ref{FBI}), which is similar to Gieseker's so-called Basic Inequality.  This turns out to be an extremely useful tool.  Later we will show that this inequality holds in more generality than is stated here (see Amplification \ref{amp}). 

{\bf Notation.}  Suppose $C$ is a curve which has at least two irreducible components, and suppose it is generically reduced on any components which do not collapse under $p_{W}$.  Let $Y$ be a union of components which do not collapse under $p_W$ and let $C'$ be the closure of $C-Y$ in $C$ as constructed above.
Let \mbox{$C' \stackrel{\iota_{C'}}{\rightarrow} C \stackrel{\iota_{C}}{\rightarrow} \pwpr$ and $Y \stackrel{\iota_{Y}}{\rightarrow} C$} be the inclusion morphisms.  Let 
\begin{equation*} 
\begin{array}{rclcrcl}
L_{W \,Y } &:=& \iota_{Y}^{*}\iota_{C}^{*}p_{W}^{*}\mathcal{O}_{\Pro(W)}(1) 
		&&L_{W \, C'} &:=& \iota_{C'}^{*}\iota_{C}^{*}p_{W}^{*}\mathcal{O}_{\Pro(W)}(1) \\ \index{$L_W$!$L_{W\,C'}$}
L_{r \, Y} &:=& \iota_{Y}^{*}\iota_{C}^{*}p_{r}^{*}\mathcal{O}_{\Pro^{r}}(1) 
		&&L_{r \, C'} &:=& \iota_{C'}^{*}\iota_{C}^{*}p_{r}^{*}\mathcal{O}_{\Pro^{r}}(1) .\index{$L_r$!$L_{r\,C'}$}
\end{array}
\end{equation*}  
Let $\pi: \bar{C} \rightarrow C$ be the normalisation morphism. Let $\bar{L}_{W\,\bar{Y}} := \pi^{*}L_{W\,Y}$ and similarly define $\bar{L}_{W\,\bar{C}'},\bar{L}_{r\,\bar{Y}}$ and $\bar{L}_{r\,\bar{C}'}$.  Normalisation induces a homomorphism
\[
\pi_{m,\mhat *}:H^0(C,L_W^m\otimes L_r^{\mhat})\rightarrow H^{0}(\bar{C}, \bar{L}_{W}^{m} \otimes \bar{L}_{r}^{\mhat}).
\]
Define 
\begin{eqnarray*}
e' &:=& \deg_{\bar{C}'} \bar{L}_{W \, \bar{C}'} = \deg_{C'} L_{W \, C'}\\ \index{$e'$}
d' &:=& \deg_{\bar{C}'}\bar{L}_{r \, \bar{C'}} = \deg_{C'}L_{r \, C'}, \index{$d'$}
\end{eqnarray*}
and let $n'$ be the number of markings on $C'$.  Write $ h^{0}:=h^{0}(p_{W}(C'), \mathcal{O}_{p_{W}(C')}(1))$.    Recall that we defined $\bar{g}$  to be 
\begin{eqnarray*} 
\bar{g}:= \mbox{min} \{ 0, g_{\bar{Y}} 
	&| & \mbox{$\bar{Y}$ is the normalisation of a complete subcurve $Y$ contained} \\ 
	&& \mbox{in a connected fibre $\mathcal{C}_{h}$ for some $h \in \Hilb(\pwpr)$} \}.
\end{eqnarray*}
\begin{proposition}[\cite{G} 1.0.7] \label{lemma_prop}\label{fund_eqn_prop} Let $M\subset\tilde{M}$, where $\tilde{M}$ consists of those $(m,\mhat,m')$ such that $m,\mhat > m_{3}$ and
\[
m >(g-\frac{1}{2}+e(q_{1}+1)+q_{3}+\mu_{1}m_2)(e-g+1)
\]
with
\[ 
d\frac{\mhat}{m} + n\frac{m'}{m^2} < \frac{1}{4} e - \frac{5}{4}g +\frac{3}{4}.
\]  
Let $l\in\mathbf{H}_M(I$).  Let $(C,x_1,\ldots ,x_n)$ be a connected marked curve which is semistable with respect to $l$.  Suppose $C$ has at least two irreducible components.  Let $C'$ and $Y$ be as above; in particular, no component of $Y$ collapses under $p_W$.  The subcurve $C'$ need not be connected. Suppose $C'$ has $b$ connected components.  Suppose there exist points $P_{1},\ldots ,P_{k}$ on $\bar{Y}$ satisfying
\begin{enumerate}
\item[(i)] $\pi ( P_{i}) \in Y \cap C'$ for all $1 \leq i \leq k$
\item[(ii)] for each irreducible component $\bar{Y}_{j}$ of $\bar{Y}$,
\[\deg_{\bar{Y}_{j}}(\bar{L}_{W \, \bar{Y}}(-D)) \geq 0, 
\]
where $D = P_{1} + \cdots + P_{k}$.
\end{enumerate}
Then there exist $(m,\mhat,m')\in M$ such that
\begin{equation} 
\label{FBI}    
e' + \frac{k}{2} < \frac{h^{0}e+ (dh^{0}-d'(e-g+1))\frac{\mhat}{m} + (nh^0-n'(e-g+1))\frac{m'}{m^2}}{e-g+1} \index{fundamental inequality}
	+ \frac{S}{m}, 
\end{equation} 
where $S =g+k(2g-\frac{3}{2})+ q_{2}- \bar{g}+1$. \index{$S$}
\end{proposition}

\begin{proof}  We start by defining the key 1-PS for this case.  Let 
\[ W_{0} := \ker \{ H^{0}(\Pro(W),\mathcal{O}_{\Pro(W)}(1)) \rightarrow H^{0}(p_{W}(C'), \mathcal{O}_{p_{W}(C')}(1)) \}.
\]
Write $N_0:=\dim W_0$.  Choose a basis $w_{0},\ldots ,w_{N_0-1},w_{N_0},\ldots ,w_N$ of \\$W_1:=H^{0}(\Pro(W), \mathcal{O}_{\Pro(W)}(1))$ relative to the filtration $0 \subset W_{0} \subset W_{1}$.  Let $\lambda_{C'}$ be the 1-PS of $GL(W)$ whose action is given by
\begin{displaymath}
\begin{array}{lcrcccccc}
\lambda_{C'}(t)w_{i} & = & w_{i}, & t \in \C^{*}, & 0 & \leq & i & \leq & N_{0}-1 \\
\lambda_{C'}(t)w_{i} & = & tw_{i}, & t \in \C^{*}, & N_{0} & \leq & i & \leq & N
\end{array} 
\end{displaymath}
and let $\lambda'_{C'}$ \label{lambda'_C'}\index{$\lambda'_{C'}$} be the associated 1-PS of $SL(W)$.

\begin{claim}\label{fund_eqn_claim}If $(C,\pts)$ satisfies $\mu^{L_{\m}}((C,\pts),\lambda'_{C'})\leq 0$ for some $(m,\mhat,m')\in M$, then (\ref{FBI}) is satisfied for this choice of $(m,\mhat,m')$.  
\end{claim}
Suppose the claim is true.  Fix $l\in\mathbf{H}_M(I)$ and suppose that $(C,\pts)$ is semistable with respect to $l$.  If there do no exist $(\m)\in M$ satisfying (\ref{FBI}) then it follows from Claim \ref{fund_eqn_claim} that $\mu^{L_{\m}}((C,\pts),\lambda'_{C'})>0$ for all $(m,\mhat,m')\in M$.  But then Lemma \ref{ignore_hull} tells us that $(C,\pts)$ is not semistable with respect to $l$.  This contradiction implies the existence of such $(\m)\in M$.

It remains to prove Claim \ref{fund_eqn_claim}, so assume that $\mu^{L_{\m}}((C,\pts),\lambda'_{C'})\leq 0$.
We shall derive the fundamental inequality from this, using Lemma \ref{numcritlemma}.

Estimate the weights for $\lambda_{C'}$ coming from the marked points.  There are $n'$ of these on $C'$, so $\sum_{l=1}^n w_{\lambda_{C'}}(w_{k_l})m'\geq n'm'$.  Also,  estimate the sum of the weights.  It is clear from the definition of $\lambda_{C'}$ that $\sum_{i=0}^N w_{\lambda_{C'}}(w_i)\leq h^0$. 

Now we look at the weight coming from the curve.  We wish to estimate the sum $\sum_{j=1}^{P(m,\mhat)}w_{\lambda_{C'}}(\hat{M}_{i_{j}})$, where $\Mhat_{i_j}$ is a monomial in $H^0(\pwpr,\opwopr{m}{\mhat})$ of bidegree $(m,\mhat)$ and $\{ \phat^C_{m,\mhat}(\hat{M}_{i_{1}}),\ldots ,\phat^C_{m,\mhat}(\hat{M}_{i_{P(m,\mhat)}}) \}$ is a basis of $H^{0}(C,L_{W}^{m}\otimes L_{r}^{\mhat})$. 

Construct a filtration of $H^{0}(C, L_W^{m}\otimes L_{r}^{\mhat})$ as at line (\ref{template_filtration}).  Namely,
for $0 \leq p \leq m$ let $\Omm_p$ be the subspace of $H^{0}(\pwpr, \opwopr{m}{\mhat})$ spanned by monomials of weight less than or equal to $p$, and let $\barO_p:=\phat^C_{m,\mhat}(\Omm_p) \subset H^{0}(C,L_{W}^{m}\otimes L_{r}^{\mhat})$.
We have a filtration of $H^{0}(C, L_W^{m}\otimes L_{r}^{\mhat})$ in order of increasing weight:
\[ 
0 \subseteq \barO_0 \subseteq  \barO_1 \subseteq \cdots \subseteq \barO_m = H^{0}(C, L_W^{m}\otimes L_{r}^{\mhat}).
\]
Write $\bhat_{p} = \dim\barO_p$.

For $p=m$, it is clear that $\bhat_m=h^0(C, L_W^{m}\otimes L_{r}^{\mhat})=em+d\mhat-g+1$.  We estimate $\bhat_p$ in the case $p\neq m$.  Restriction to $Y$ induces a homomorphism 
\[
\rho^{Y,C}_{m,\mhat}:H^0(C,L_{W}^{m}\otimes L_{r}^{\mhat})\rightarrow H^0(Y,L_{WY}^{m}\otimes L_{rY}^{\mhat}).
\]
We restrict this to $\barO_p$, where $0\leq p < m$.  Note that if $\hat{M}$ is a monomial in $\barO_p$ and $p<m$ then $\hat{M}$ has at least one factor from $W_0$ and hence by definition $\hat{M}$ vanishes on $C'$.  If such $\hat{M}$ also vanishes on $Y$ then $\hat{M}$ is zero on $C$.  Hence the restriction $\rho^{Y,C}_{m,\mhat}|_{\barO_p}$ has zero kernel, so is an isomorphism of vector spaces, and thus:
\[
\dim \phat^{Y,C}_{m,\mhat}(\barO_p) = 
	\dim \barO_p = \bhat_p.
\]
We denote $\phat^{Y,C}_{m,\mhat}(\barO_p)$ by $\barO_p|_Y$.

The normalisation morphism $\pi_Y: \bar{Y} \rightarrow Y$ induces a homomorphism 
\[
\pi_{Y\,m,\mhat*}: H^{0}(Y,L_{W}^{m}\otimes L_{r}^{\mhat}) \rightarrow  H^{0}(\bar{Y},\bar{L}_{W}^{m}\otimes \bar{L}_{r}^{\mhat}). 
\]
By definition the sections in $\pi_{Y\,m,\mhat*}(\barO_p|_Y)$ vanish to order at least $m-p$ at the points $P_{1},\ldots ,P_{k}$.  Thus 
\[
\pi_{Y\,m,\mhat*}(\barO_p|_Y) 
	\subseteq H^{0}(\bar{Y},\bar{L}_{W \bar{Y}}^{m}\otimes \bar{L}_{r\bar{Y}}^{\mhat}(-(m-p)D)).
\]
Then 
\begin{eqnarray}
\bhat_{p} &=& \dim (\barO_p)_Y
	 	\ \leq \ h^{0}(\bar{Y},\bar{L}_{W \bar{Y}}^{m} \otimes \bar{L}_{r\bar{Y}}^{\mhat}(-(m-p)D)) 
		+ \dim \ker \pi_{Y\,m,\mhat*}  \notag\\
	& = & (e-e')m + (d-d')\mhat - k(m-p) - \bar{g} + 1 \notag\\
		&& {} + h^{1}(\bar{Y},\bar{L}_{W \bar{Y}}^{m}\otimes \bar{L}_{r\bar{Y}}^{\mhat} (-(m-p)D)) + \dim \ker \pi_{Y\,m,\mhat*}. \label{bp1}
\end{eqnarray}

We apply the estimates of Claim \ref{I,II,III} to our current situation.

\begin{list}
{}{\setlength{\leftmargin}{.75in} \setlength{\rightmargin}{.75in} \setlength{\itemsep}{-1.5mm} \setlength{\itemindent}{-.25in}}
\item{I.} $\dim \ker \pi_{Y\,m,\mhat *} < q_{2}$.  
\end{list}
No component of $Y$ collapses under $p_W$, so by Proposition \ref{genred} the curve $Y$ is generically reduced.  Claim \ref{I,II,III}.I may be applied to $Y\subset\pwpr$. Let $\mathcal{I}_{Y}$ denote the ideal sheaf of nilpotents in $\mathcal{O}_{Y}$.  Then $\dim\ker\pi_{Y\,m,\mhat *}<h^0(Y,\mathcal{I}_Y)$.  In Proposition \ref{7facts}, the constant $q_{2}$ was defined to have the property $h^{0}(C,\mathcal{I}_{C})<q_{2}$; hence $h^{0}(Y,\mathcal{I}_{Y})<q_{2}$ as well.  

\begin{list}
{}{\setlength{\leftmargin}{.75in} \setlength{\rightmargin}{.75in} \setlength{\itemsep}{-1.5mm} \setlength{\itemindent}{-.25in}}
\item{II.}  $h^{1}(\bar{Y},\bar{L}_{W\,\bar{Y}}^{m}\otimes\bar{L}_{r\,\bar{Y}}^{\mhat}(-(m-p)D)) \leq k(m-p) \leq km$ if $0 \leq p \leq 2g-2$. 
\end{list}
The sheaf $\bar{L}_{W\,\bar{Y}}^{m}\otimes\bar{L}_{r\,\bar{Y}}^{\mhat}$ is locally free on $\bar{Y}$, and we have chosen $m$ and $\mhat$
so that $H^1(\bar{Y},\bar{L}_{W\,\bar{Y}}^{m}\otimes\bar{L}_{r\,\bar{Y}}^{\mhat})=0$.  The hypotheses of Claim \ref{I,II,III}.II hold, and we calculate $\deg (m-p)D=k(m-p)$.  We make a coarser estimate than we could as this will be sufficient for our purposes.

\begin{list}
{}{\setlength{\leftmargin}{.75in} \setlength{\rightmargin}{.75in} \setlength{\itemsep}{-1.5mm} \setlength{\itemindent}{-.25in}}
\item{III.}  $h^{1}(\bar{Y},\bar{L}_{W\,\bar{Y}}^{m}\otimes\bar{L}_{r\,\bar{Y}}^{\mhat}(-(m-p)D))=0$ if $2g-1 \leq p \leq m-1$.  
\end{list}
$\bar{Y}$ is reduced, and is a union of disjoint irreducible (and hence integral) components $\bar{Y}_j$ of genus $g_{Y_j}\leq g$.  We apply Claim \ref{I,II,III}.III separately to each component.  Our assumption (ii) was that $\deg_{\bar{Y}_{j}}(\bar{L}_{W\,\bar{Y}}(-D))\geq 0$, so 
\[
\deg_{\bar{Y}_{j}}(\bar{L}_{W\,\bar{Y}}) \geq \deg_{\bar{Y}_j} D.
\]
If $\deg_{\bar{Y}_j} D \geq 1$ then 
\begin{eqnarray*}
\deg_{\bar{Y}_j}(\bar{L}_{W\,\bar{Y}}^m \otimes \bar{L}_{r\,\bar{Y}}^{\mhat}(-(m-p)D)) 
	&\geq& m (\deg_{\bar{Y}_j} D) -(m-p)(\deg_{\bar{Y}_j} D) \\
	&\geq& p \, \geq \, 2g-1 \, \geq \, 2g_{Y_j}-1,
\end{eqnarray*}
as required.  On the other hand, suppose that $\deg_{\bar{Y}_j} D=0$.  We assumed that no component of $Y$ collapses under $p_W$, and hence for each $j$, the degree $\deg_{\bar{Y}_j}(\bar{L}_{W\,\bar{Y}})\geq 1$.  Thus again:
\[
\deg_{\bar{Y}_j}(\bar{L}_{W\,\bar{Y}}^m \otimes \bar{L}_{r\,\bar{Y}}^{\mhat}(-(m-p)D))\geq m \geq 2g_{Y_j}-1.
\]

Combining this data with our previous formula (\ref{bp1}) we have shown:
\begin{displaymath}
\bhat_{p} \leq \left \{ \begin{array}{ll} (e-e'-k)m +(d-d')\mhat+kp-\bar{g}+1+q_{2}+km & 0 \leq p \leq 2g-2\\
					(e-e'-k)m + (d-d')\mhat +kp-\bar{g}+1+ q_{2}, 	& 2g-1 \leq p \leq m-1 .
\end{array} \right.
\end{displaymath}
Thus, we may use Lemma \ref{main_calculation}, setting $\alpha=e'+k$, $\beta=d'$, $\gamma=k$, $\delta=n'$, $\epsilon={}-\bar{g}+1+q_{2}+(2g-1)k$, $r_N=1$ and $R=h^0$. Following Lemma \ref{main_calculation}, we see
\begin{eqnarray*}
\lefteqn{ \mu^{L_{\m}}((C,\pts),\lambda')
	\geq \left((e'+\frac{k}{2})(e-g+1) -h^0e\right)m^2 } \\
	&\hspace{.5in}& {}+\left(d'(e-g+1)-dh^0\right)m\mhat + \left(n'(e-g+1)-nh^0\right)m'\\
	&\hspace{.5in}& {}-\left(g-\frac{k}{2}-\bar{g}+q_{2}+(2g-1)k\right)(e-g+1)m - h^0m. 
\end{eqnarray*}
Thus, since we assume that $\mu^{L_{\m}}((C,\pts),\lambda)\leq 0$, we may conclude that
\begin{eqnarray*}
e' + \frac{k}{2} 
	& < & \frac{h^{0}e+ (dh^{0}-d'(e-g+1))\frac{\mhat}{m} + (nh^0-n'(e-g+1))\frac{m'}{m^2}}{e-g+1} + \frac{S}{m},
\end{eqnarray*}
where $S = g+k(2g-\frac{3}{2})+q_{2}-\bar{g}+1$.
\end{proof}
This fundamental inequality allows us finally to show that no irreducible components of $C$ collapse under projection to $\Pro(W)$.

\begin{proposition}[\cite{G} 1.0.3] \label{e'>0} Let $M$ consist of those $(\m)$ such that  $m,\mhat > m_{3}$ and
\[m > \max \left \{ \begin{array}{c} 
	( g - \frac{1}{2} + e(q_{1}+1) + q_{3} +\mu_{1}m_2)(e-g+1), \\ 
	(6g+2q_2-2\bar{g})(e-g+1) 
\end{array} \right \}
\]
with
\[ 
d\frac{\mhat}{m} + n\frac{m'}{m^2} < \frac{1}{4} e - \frac{5}{4}g +\frac{3}{4}
\]  
while
\[
\frac{\mhat}{m} > 1 + \frac{\frac{3}{2}g-1+d+n\frac{m'}{m^2}}{e-g+1-d}.
\]
Let $l\in\mathbf{H}_M(I)$.  If $C$ is connected and $(C,\pts)$ is semistable with respect to $l$, then no irreducible components of $C$ collapse under $p_{W}$. 
\end{proposition}
\emph{Remark.} As the denominator $e-g+1-d$ is equal to $(2a-1)(g-1)+an+(ca-1)d$ it is evident that this is positive.

\begin{proof}This is trivial if $d=0$;  assume that $d\geq 1$.  Suppose that at least one component of $C$ collapses under $p_{W}$.  Let $C'$ be the union of all irreducible components of $C$ which collapse under $p_{W}$ and let $Y:=\overline{C-C'}$.  Suppose that $C'$ consists of $b$ connected components, namely $C'_1,\ldots ,C'_b$.  If $d'_i=\deg_{C'_i}L_{r \, C'_i}$ then $d'_i\geq 1$ since $e'_i=0$, for $i=1,\ldots ,b$.  But then $d'=\deg_{C'}L_{r\,C'}=\sum_{i=1}^b d'_i\geq b$.  Hence 
\[
1\leq b\leq d'\leq d.
\]

The curve $C$ is connected so $C' \cap Y \neq \emptyset$.  Choose one point $P \in \bar{Y}$ such that $\pi(P)\in C'\cap Y$.   We have by definition $\deg_{\bar{Y}_{j}}(\bar{L}_{W \, \bar{Y}}) \geq 1$ so $ \deg_{\bar{Y}_{j}}(\bar{L}_{W \, \bar{Y}}(-P)) \geq 0$ for each irreducible component $Y_j$ of $Y$.  The hypotheses of Proposition \ref{fund_eqn_prop} are satisfied for $k=1$, and $M$ as in the statement of this proposition.  Let $(\m)\in M$ be the integers which that corollary provides, satisfying (\ref{FBI}).

Since $C'$ consists of $b$ connected components, it is collapsed to at most $b$ distinct points under $p_W$, so we estimate 
\[
h^{0}(p_{W}(C'), \mathcal{O}_{p_{W}(C')}(1)) \leq b.  
\]
Recall that we defined $S=g + k(2g-\frac{3}{2}) + q_{2} - \bar{g}+1$.  In the current situation, $k=1$, so $S< 3g+q_2-\bar{g}$.  The hypotheses on $m$ imply then that $\frac{S}{m}(e-g+1)<\frac{1}{2}$.  Estimate $n'\geq 0$.  Then the inequality (\ref{FBI}) reads:
\begin{eqnarray*} 
0+\frac{1}{2} 
	\leq e' + \frac{k}{2} 
	&\leq &\frac{be + (bd-b(e-g+1))\frac{\mhat}{m} + bn\frac{m'}{m^2}}{e-g+1} + \frac{1}{2(e-g+1)} \\
\Rightarrow \frac{1}{2b}(e-g+1) 
	&\leq& e  + (d-(e-g+1))\frac{\mhat}{m} + n\frac{m'}{m^2} + \frac{1}{2b} \\
\Rightarrow (e-g+1-d)\frac{\mhat}{m}
	&\leq&  \left(1-\frac{1}{2b}\right)e +\frac{1}{2b}g + n\frac{m'}{m^2} 
		\ \leq \ e +\frac{1}{2}g + n\frac{m'}{m^2} \\
\Rightarrow \frac{\mhat}{m}
	& \leq & 1 + \frac{\frac{3}{2}g-1+d+n\frac{m'}{m^2}}{e-g+1-d}.
\end{eqnarray*}
We have contradicted our hypothesis that $\frac{\mhat}{m}>1 + \frac{\frac{3}{2}g-1+d+n\frac{m'}{m^2}}{e-g+1-d}$.
\end{proof}

{\it Remark.}  We have now described a range $M$ of $(\m)$, such that if $l\in\mathbf{H}_M(I)$ and $(C,\pts)$ is semistable with respect to $l$, then the map $p_{W}|_{C}: C \rightarrow p_{W}(C)$ is surjective, finite, and generically 1-1.  Further, since no components of $C$ are collapsed under $p_W$, it follows from Proposition \ref{genred} that $C$ is generically reduced.\\

One may check that there exist integers $(\m)$ such that all stable maps have a model satisfying the inequality (\ref{FBI}) of Proposition \ref{fund_eqn_prop}.  Such a calculation is carried out in \cite{e_thesis} Proposition 5.1.8.  It turns out that one may easily show that the inequality is satisfied by any complete subcurve $C'\subset C$, if $\frac{\mhat}{m} = \frac{c \cp}{2 \cp-1}$, and $\frac{m'}{m^2}=\frac{\cp}{2 \cp-1}$ for $l = 1,\ldots, n$. We will be able to use the theory of variation of GIT to show that in fact the quotient is constant in a small range around this key linearisation.  


\subsection{GIT semistability implies that the only singularities are nodes}\label{nodes}	

The next series of results provides a range $M$ of triples $(\m)$ such that if $l\in\mathbf{H}_M(I)$ and if the connected curve  $(C,x_1,\ldots ,x_n)$ is semistable with respect to $l$, then any singularities of $C$ are nodes.  First we show that $C$ has no cusps by showing normalisation morphism $\pi:\bar{C}\rightarrow C$ is unramified.  Singular points are shown to be double points by showing that the inverse image under $\pi$ of any $P\in C$ contains at most two points. We must also rule out tacnodes; these occur at double points $P$ such that the two tangent lines to $C$ at $P$ coincide.

In the hypotheses of the following lemma, note that $6g+2q_2-2\bar{g}\leq 9g+3q_2-3\bar{g}$ and that $2e-10g+6>e-9g+7$, and so in particular the hypotheses of Proposition \ref{e'>0} hold.


\begin{proposition}[cf.\ \cite{G} 1.0.5] \label{unram_prop} Let $a$ be sufficiently large that $e-9g+7=a(2g-2+n+cd)-9g+7>0$, and let $M$ consist of those $(m,\mhat,m')$ such that $m,\mhat > m_{3}$ and
\[m > \max \left \{ \begin{array}{c} 
	( g - \frac{1}{2} + e(q_{1}+1) + q_{3} +\mu_{1}m_2)(e-g+1), \\ 
	(9g+3q_{2}-3\bar{g})(e-g+1)
\end{array} \right \} 
\]
with
\[
d\frac{\mhat}{m} + n\frac{m'}{m^2} 
	< \frac{1}{8}e-\frac{9}{8}g + \frac{7}{8}
\]
while
\[
\frac{\mhat}{m} > 1 + \frac{\frac{3}{2}g-1+d+n\frac{m'}{m^2}}{e-g+1-d}.
\]
Let $l\in\mathbf{H}_M(I)$.  If $(C,x_1,\ldots ,x_n)$ is connected and semistable with respect to $l$, then the normalisation morphism $\pi: \bar{C} \rightarrow C_{\red}$ is unramified.  In particular, $C$ possesses no cusps.
\end{proposition}


\begin{proof}  Suppose $\pi$ is ramified at $P \in \bar{C}$.  Then $p_{W} \circ \pi: \bar{C} \rightarrow p_{W}(C_{\red})$ is also ramified at $P$.  Recall that by Proposition \ref{1.0.2hat}, the curve $p_{W}(C) \subset \Pro(W)$ is nondegenerate; we can think of $H^{0}(\Pro(W), \mathcal{O}_{\Pro(W)}(1))$ as a subspace of $H^{0}(p_{W}(C), L_W)$.  Define
\begin{displaymath}
\begin{array}{lll}
W_{0} & = & \{ s \in H^{0}(\Pro(W), \mathcal{O}_{\Pro(W)}(1)) | \mbox{$\pi_{*}p_{W \, *}s$ vanishes to order $\geq 3$ at $P$} \} \\
W_{1} & = & \{ s \in H^{0}(\Pro(W), \mathcal{O}_{\Pro(W)}(1)) | \mbox{$\pi_{*}p_{W \, *}s$ vanishes to order $\geq 2$ at $P$} \}
\end{array}
\end{displaymath}
and write $N_{0} := \dim W_{0}$ and $N_{1} := \dim W_{1}$.  Choose a basis 
\[ w_{1},\ldots,w_{N_{0}},w_{N_{0}+1},\ldots,w_{N_{1}},w_{N_{1}+1},\ldots,w_{N+1}
\] 
of $W_{3} := H^{0}(\Pro(W), \mathcal{O}_{\Pro(W)}(1))$ relative to the filtration $0 \subset W_{0} \subset W_{1} \subset W_{3}$.  Let $\lambda$ be the 1-PS of $GL(W)$ whose action is given by
\begin{displaymath}
\begin{array}{lcrcrcccl}
\lambda(t)w_{i} & = & w_{i}, & t \in \C^{*}, & 1 & \leq & i & \leq & N_{0} \\
\lambda(t)w_{i} & = & tw_{i}, & t \in \C^{*}, & N_{0}+1 & \leq & i & \leq & N_{1} \\
\lambda(t)w_{i} & = & t^{3}w_{i}, & t \in \C^{*}, & N_{1}+1& \leq & i &  \leq & N+1 
\end{array} 
\end{displaymath}
and let $\lambda'$ be the associated 1-PS of $SL(W)$.  
Pick some $(m,\mhat,m')\in M$.

We make an estimate for $\mu^{L_{\m}}((C,\pts),\lambda')$ using the filtration of (\ref{template_filtration}) and Lemma \ref{main_calculation}.  For $0\leq i \leq 3m$ we let $\Omm_i$ be the subspace of $H^{0}(\pwpr, \opwopr{m}{\mhat})$ spanned by monomials of weight less than or equal to $i$, and let $\barO_i:=\phat^C_{m,\mhat}(\Omm_i) \subset H^{0}(C,L_{W}^{m}\otimes L_{r}^{\mhat})$.  We need to find an estimate for $\bhat_i:=\dim\barO_i$.

As in Proposition \ref{fund_eqn_prop} we use the homomorphism
\[
\pi_{m,\mhat *}:H^0(C,L_W^m\otimes L_r^{\mhat})\rightarrow H^{0}(\bar{C}, \bar{L}_{W}^{m} \otimes \bar{L}_{r}^{\mhat})
\]
induced by the normalisation morphism.  We show that 
\begin{equation}
\pi_{m,\mhat *}(\barO_i) \subseteq H^{0}(\bar{C}, \bar{L}_{W}^{m} \otimes \bar{L}_{r}^{\mhat} \otimes \mathcal{O}_{\bar{C}}((-3m+i)P)), \label{correct_vanishing}
\end{equation}
 for $0 \leq i \leq 3m$.  When $i=0$ this follows from the definitions.  For $1 \leq i \leq 3m$ it is enough to show that monomial $\hat{M} \in  (\barO_i-\barO_{i-1})$ vanishes at $P$ to order at least $3m-i$.  Suppose that such $\hat{M}$ has $i_{0}$ factors from $W_0$, $i_{1}$ factors from $W_1$ and $i_{3}$ factors from $W_3$.  Then $i_{0}+i_{1}+i_{3}=m$ and $i_{1}+3i_{3}=i$.  By definition, $\hat{M}$ vanishes at $P$ to order at least $3i_0+2i_1$.  But 
\[
  3i_{0}+2i_{1} = 3(i_{0}+i_{1}+i_{3}) - (i_{1}+3i_{3}) = 3m-i,
\]
so the monomial vanishes as required and hence (\ref{correct_vanishing}) is satisfied.

 By (\ref{correct_vanishing}) and Riemann-Roch,
\begin{eqnarray*} 
\bhat_{i} :=  \dim \barO_i
		&\leq&   h^{0}(\bar{C},\bar{L}_{W}^{m}\otimes \bar{L}_{r}^{\mhat} \otimes \mathcal{O}_{\bar{C}}(-(3m-i)P)) 
			+ \dim \ker \pi_{m,\mhat * } \\
	&  \leq &  em+d\mhat -3m+i-\bar{g} +1 \\
		&&{}+h^{1}(\bar{C},\bar{L}_{W}^{m} \otimes \bar{L}_{r}^{\mhat} \otimes \mathcal{O}_{\bar{C}}(-(3m-i)P))+\dim \ker \pi_{m,\mhat * }.
\end{eqnarray*}
We may use Claim \ref{I,II,III} in a straightforward way to show that $\dim \ker \pi_{m,\mhat *} < q_{2}$ and that $h^{1}\left(\bar{C},\bar{L}_{W}^{m}\otimes \bar{L}_{r}^{\mhat} \otimes \mathcal{O}_{\bar{C}}(-(3m-i)P)\right) \leq 3m-i\leq 3m$ if $0 \leq i \leq 2g-2$.  More care is needed to show that the $h^1$ term vanishes for higher values of $i$.

\label{degree_at_least_3}
Let $\bar{C}_i$ be an irreducible component of $\bar{C}$.  Suppose $\bar{C}_i$ does not contain $P\in\bar{C}$.  We have shown (Proposition \ref{e'>0}) that $\deg_{\bar{C}_i}\bar{L}_{W} =\deg_{C_i}L_W\geq 1$.  Thus
\[
\deg_{\bar{C}_i}\left(\bar{L}_W^{m} \otimes \bar{L}_{r}^{\mhat} \otimes \mathcal{O}_{\bar{C}}(-(3m-i)P)\right)
	=\deg_{\bar{C}_i}\left(\bar{L}_W^{m} \otimes \bar{L}_{r}^{\mhat}\right)\geq m \geq 2g_{C_i}-1.
\]
On the other hand, suppose that $\bar{C}_i$ is the component of $\bar{C}$ on which $P$ lies.  The morphism $\bar{C}_i\rightarrow p_W(C_{i\,\red})$ is ramified at $P$, so $p_W(C_{i\,\red})$ is singular and integral in $\Pro(W)$.  Were $p_W(C_{i\,\red})$ an integral curve of degree $1$ or $2$ in $\Pro(W)$, it would be either a line or a conic and hence non-singular.  We conclude that $\deg_{\bar{C}_i}\bar{L}_W=\deg_{C_{i\,\red}}L_W\geq 3$.  Then
\[
\deg_{\bar{C}_i} \left(\bar{L}_{W}^{m} \otimes \bar{L}_{r}^{\mhat} \otimes \mathcal{O}_{\bar{C}}(-(3m-i)P) \right)\geq 3m-3m+i=i.
\]
Thus,  Claim \ref{I,II,III}.III shows that $h^{1}\left(\bar{C},\bar{L}_{W}^{m}\otimes \bar{L}_{r}^{\mhat} \otimes \mathcal{O}_{\bar{C}}(-(3m-i)P)\right) = 0$ if $2g-1\leq i\leq 3m-1$.

Combining these inequalities we have
\[
\bhat_{i} \leq \left\{ \begin{array}{ll} (e-3)m+d\mhat+i-\bar{g}+q_2+1+3m, 	& 0 \leq i \leq 2g-2 \\
					(e-3)m+d\mhat+i-\bar{g}+q_2+1, 	& 2g-1 \leq i \leq 3m-1.
\end{array} \right.
\]
Thus, in the language of Lemma \ref{main_calculation}, we shall set $\alpha=3$, $\beta=0$, $\gamma=1$, and $\epsilon={}-3\bar{g}+3q_2+6g$.  We may estimate the minimum weight of the action of $\lambda$ on the marked points $x_i$ as zero, so we set $\delta=0$.  We know that $r_N=3$.  It remains to find an upper bound for $\sum w_{\lambda}(w_i)$.

Recall that we are regarding $W$ as a subspace of $H^0(p_W(C),L_W)$.  Note that the image of $W_{0}$ under  $\pi_{*}$ is contained in $H^{0}(\bar{C}, \bar{L}_{W}(-3P))$, and the image of $W_{1}$ under  $\pi_{*}$ is contained in $H^{0}(\bar{C}, \bar{L}_{W}(-2P))$.  We have two exact sequences 
\begin{displaymath}
\begin{array}{ccccccccc}  
0 & \rightarrow & \bar{L}_{W}(-P) & \rightarrow & \bar{L}_{W} & \rightarrow & k(P) & \rightarrow & 0 \\ 
0 & \rightarrow & \bar{L}_{W}(-3P) & \rightarrow & \bar{L}_{W}(-2P) & \rightarrow & k(P) & \rightarrow & 0, \\ 
\end{array}
\end{displaymath}
which give rise to long exact sequences in cohomology
 \begin{displaymath}
\begin{array}{cccccccc}  
0 & \rightarrow & H^{0}(\bar{C},\bar{L}_{W}(-P)) & \rightarrow & H^{0}(\bar{C}, \bar{L}_{W}) & \rightarrow & H^{0}(\bar{C},k(P)) & \rightarrow \cdots \\ 
0 & \rightarrow &  H^{0}(\bar{C},\bar{L}_{W}(-3P)) & \rightarrow &  H^{0}(\bar{C},\bar{L}_{W}(-2P)) & \rightarrow &  H^{0}(\bar{C},k(P)) & \rightarrow  \cdots. \\ 
\end{array}
\end{displaymath}
The second long exact sequence implies that $\dim W_{1}/W_{0} \leq 1$.  Now recall that \mbox{ $\bar{L}_{W} := \pi^{*}(L_{W})$} and $\pi$ is ramified at $P$.  The ramification index must be at least two, so we have \mbox{$H^{0}(\bar{C},\bar{L}_{W}(-P)) =   H^{0}(\bar{C},\bar{L}_{W}(-2P))$}.  Then the first long exact sequence implies that $\dim W_{3}/W_{1} \leq 1$.  We conclude that $\sum_{i=1}^{N+1}w_{\lambda}(w_{i}) \leq 1+3 = 4=:R$.

We may now estimate the $\lambda'$-weight for $(C,\pts)$, using Lemma \ref{main_calculation}:
\begin{eqnarray*}
\lefteqn{ \mu^{L_{\m}}((C,\pts),\lambda')
	\geq \left(\frac{9}{2}(e-g+1)-4e\right)m^2 - 4dm\mhat-4nm' }\\
		&&{}-\left((9g-3\bar{g}+3q_2-\frac{9}{2})(e-g+1)+4\right)m \\
	&\geq& \left(\frac{1}{2}e-\frac{9}{2}(g-1)-4d\frac{\mhat}{m}-4n\frac{m'}{m^2}\right)m^2 - \left(9g-3\bar{g}+3q_2\right)(e-g+1)m.
\end{eqnarray*}
We assumed that $\left(9g-3\bar{g}+3q_2\right)(e-g+1)<m$, so we have shown that
\[
\mu^{L_{\m}}((C,\pts),\lambda') \geq \left(\frac{1}{2}e-\frac{9}{2}g+\frac{7}{2}-4d\frac{\mhat}{m}-4n\frac{m'}{m^2}\right)m^2.
\]
This is clearly positive, as we assumed that
\[
\frac{1}{8}e-\frac{9}{8}g+\frac{7}{8}>d\frac{\mhat}{m}+n\frac{m'}{m^2},
\]
so $m^2$ has a positive coefficient.  
 
Thus $\mu^{L_{m,\mhat,m'}}((C,\pts),\lambda')>0$.  This is true for all $(m,\mhat,m')\in M$, and therefore by Lemma \ref{ignore_hull}, the $n$-pointed curve $(C,x_1,\ldots ,x_n)$ is not semistable with respect to $l$ for any $l\in\mathbf{H}_M(I)$. 
\end{proof}
{\it Remark.}  Note that the value $e-9g+7$ is positive for any $(g,n,d)$ as long as $a\geq 10$, but smaller values of $a$ will suffice in many cases; for example, if $g\geq 3$ then $a\geq 5$ is sufficient.

\begin{proposition}[cf.\ \cite{G} 1.0.4] \label{dblpts}\label{double_points} Let $a$ be sufficiently large that $e-9g+7>0$, and let $M$ consist of those $(m,\mhat,m')$ such that $m,\mhat > m_{3}$ and
\[m > \max \left \{ \begin{array}{c} 
	( g - \frac{1}{2} + e(q_{1}+1) + q_{3} +\mu_{1}m_2)(e-g+1), \\ 
	(9g+3q_{2}-3\bar{g})(e-g+1)
\end{array} \right \} 
\]
with
\[
d\frac{\mhat}{m} + n\frac{m'}{m^2} 
	< \frac{1}{8}e-\frac{9}{8}g + \frac{7}{8}
\]
while
\[
\frac{\mhat}{m} > 1 + \frac{\frac{3}{2}g-1+d+n\frac{m'}{m^2}}{e-g+1-d}.
\]
Let $l\in\mathbf{H}_M(I)$.  If $(C,x_1,\ldots ,x_n)$ is connected, and semistable with respect to $l$, then all singular points of $C_{\red}$ are double points. 
\end{proposition}


\begin{proof}  Suppose there exists a point $P \in C$ with multiplicity $\geq 3$ on $C_{\red}.$  Let \\ $ev: H^{0}(\Pro(W), \mathcal{O}_{\Pro(W)}(1)) \rightarrow k(P)$ be the evaluation map.  Let $W_{0} = \ker ev$.  We have $N_{0} := \dim W_{0} = N$.  Choose a basis of $W_{0}$ and extend it to a basis $w_0,\ldots ,w_N$ of $W_1:=H^{0}(\Pro(W), \mathcal{O}_{\Pro(W)}(1))$.  Let $\lambda$ be the 1-PS of $GL(W)$ whose action is given by
\begin{displaymath}
\begin{array}{lcrcl}
\lambda(t)w_{i} & = & w_{i}, & t \in \C^{*}, & 0  \leq  i  \leq  N-1 \\
\lambda(t)w_{i} & = & tw_{i}, & t \in \C^{*}, & i=N
\end{array} 
\end{displaymath}
and let $\lambda'$ be the associated 1-PS of $SL(W)$.  Pick $(m,\mhat,m')\in M$.

We shall estimate $\mu^{L_{\m}}((C,\pts),\lambda')$ using Lemma \ref{main_calculation}.  Construct a filtration of $H^{0}(C, L_W^{m}\otimes L_{r}^{\mhat})$ of increasing weight as in (\ref{template_filtration}): we let $\Omm_p$ be the subspace of $H^0(\pwpr,\opwopr{m}{\mhat})$ spanned by monomials of weight less than or equal to $p$, and let $\barO_p:=\phat^C_{m,\mhat}(\Omm_p)$.  We need to find an upper bound for $\bhat_p:=\dim\barO_p$.

Define a divisor $D$ on $\bar{C}$ as follows:  Let $\pi: \bar{C} \rightarrow C$ be the normalisation morphism.
The hypotheses of Proposition \ref{unram_prop} are satisfied, \label{where to use unram} so $\pi$ is unramified; as $P$ has multiplicity at least $3$ there must be at least three distinct points in the preimage $\pi^{-1}(P)$.  Let $D = Q_{1}+Q_{2}+Q_{3}$ be three such points.  Note that if any two of these points lie on the same component $\bar{C}_1\subset\bar{C}$, then the corresponding component $C_1\subset C$ must have $\deg_{C_1}L_W\geq 3$, by the same argument as in the proof of Proposition \ref{unram_prop}.
  
The normalisation morphism induces a homomorphism 
\[ \pi_{m,\mhat*}: H^{0}(C,L_{W}^{m}\otimes L_{r}^{\mhat}) \rightarrow  H^{0}(\bar{C},\bar{L}_{W}^{m}\otimes \bar{L}_{r}^{\mhat}).
\]  
Note that $\pi_{m,\mhat*}(\barO_p) \subseteq  H^{0}(\bar{C},\bar{L}_{W}^{m} \otimes \bar{L}_{r}^{\mhat} \otimes \mathcal{O}_{\bar{C}}(-(m-p)D)).$  We have
\begin{eqnarray*} 
\bhat_{p} := \dim\barO_p
	&\leq& h^{0}(\bar{C},\bar{L}_{W}^{m} \otimes \bar{L}_{r}^{\mhat} \otimes \mathcal{O}_{\bar{C}}(-(m-p)D)) + \dim \ker \pi_{m,\mhat\,*} \\
	&  \leq &  em + d\mhat -3(m-p)-\bar{g}+1\\
	&&{}+ h^{1}(\bar{C},\bar{L}_{W}^{m}\otimes \bar{L}_{r}^{\mhat} \otimes \mathcal{O}_{\bar{C}}(-(m-p)D)) 
			+ \dim \ker \pi_{m,\mhat\, * } .
\end{eqnarray*}
We may use Claim \ref{I,II,III} parts I and II to estimate that $h^{0}(C,\mathcal{I}_{C})<q_{2}$ and that $h^{1}(\bar{C},\bar{L}_{W}^{m}(-(m-p)D) \otimes \bar{L}_{r}^{\mhat}) \leq 3(m-p)\leq 3m$ if $0\leq p \leq 2g-2$.  To show that $h^{1}(\bar{C},\bar{L}_{W}^{m}(-(m-p)D) \otimes \bar{L}_{r}^{\mhat}) =0$ if $p \geq 2g-1$, we may verify that the degree of $L_W$ on any component $C_1$ meeting $P$ implies that the hypothesis of Claim \ref{I,II,III}.III is satisfied.  

Thus:
\begin{displaymath}
\bhat_{p} \leq \left\{ \begin{array}{ll} 
	(e-3)m+d\mhat+3p -\bar{g}+q_2+1+3m, & 0 \leq p \leq 2g-2 \\
	(e-3)m+d\mhat+3p -\bar{g}+q_2+1, & 2g-1 \leq p \leq m.
\end{array} \right.
\end{displaymath}
We may apply Lemma \ref{main_calculation}, setting $\alpha=3$, $\beta=0$, $\gamma=3$ and $\epsilon={}-\bar{g}+q_2+6g-2$.  We know $r_N=1$ and may estimate the weight of the action of $\lambda$ on the marked points $\pts$ as greater than or equal to zero, so we set $\delta=0$.  Finally, note that $\sum_{i=0}^{N}w_{\lambda}(w_{i}) = 1=:R$.  Now, substituting in these values:
\begin{eqnarray*}
\lefteqn{ \mu^{L_{\m}}((C,\pts),\lambda')
	\geq  \left(\frac{3}{2}(e-g+1)-e\right)m^2 -dm\mhat - nm'}\\
		&\hspace{0in}& {}-\left((7g-\bar{g}+q_2-\frac{9}{2})(e-g+1)+1\right)m \\
	&\hspace{0in}\geq& \left(\frac{1}{2}e-\frac{3}{2}(g+1)-d\frac{\mhat}{m}-n\frac{m'}{m^2}\right)m^2
				{}-\left(7g-\bar{g}+q_2\right)(e-g+1)m.
\end{eqnarray*}
Our assumptions imply that $(7g-\bar{g}+q_2)(e-g+1)<m$.  We have shown that
\[
\mu^{L_{\m}}((C,\pts),\lambda')
	\geq \left(\frac{1}{2}e-\frac{3}{2}g+\frac{1}{2}-d\frac{\mhat}{m}-n\frac{m'}{m^2}\right)m^2.
\]
However, we assumed that
\[
d\frac{\mhat}{m}+n\frac{m'}{m^2} < \frac{1}{8}e-\frac{9}{8}g+\frac{7}{8},
\]
and $\frac{1}{8}e-\frac{9}{8}g+\frac{7}{8}<\frac{1}{2}e-\frac{3}{2}g+\frac{1}{2}$ since $e-g\geq 4$, so the coefficient of $m^2$ is positive.

Thus $\mu^{L_{m,\mhat,m'}}((C,\pts),\lambda')>0$.  This holds for all $(m,\mhat,m')\in M$, so by Lemma \ref{ignore_hull} we see $(C,x_1,\ldots ,x_n)$ is not semistable with respect to $l$ for any $l\in\mathbf{H}_M(I)$. 
\end{proof}


The remaining case we must rule out is that $C$ possesses a tacnode.  The analogous proposition in \cite{G} is 1.0.6, but the proof of that contains at least two errors (one should use $\Omega_{i}^{m}$ accounting rather than the Tata notes' $W_{i}^{m-r}W_{j}^{r}$ when the filtration has more than two stages, and tacnodes need not be separating).  These may be avoided if we simply follow \cite{HM} 4.53 instead.
\begin{proposition}[cf.\ \cite{HM} 4.53] \label{tac} Let $a$ be sufficiently large that $e-9g+7>0$, and let $M$ consist of those $(m,\mhat,m')$ such that $m,\mhat > m_{3}$ and
\[m > \max \left \{ \begin{array}{c} 
	( g - \frac{1}{2} + e(q_{1}+1) + q_{3} +\mu_{1}m_2)(e-g+1), \\ 
	(10g+3q_{2}-3\bar{g})(e-g+1) 
\end{array} \right \} 
\]
with
\[
d\frac{\mhat}{m} + n\frac{m'}{m^2} 
	< \frac{1}{8}e-\frac{9}{8}g + \frac{7}{8}
\]
while
\[
\frac{\mhat}{m} > 1 + \frac{\frac{3}{2}g-1+d+n\frac{m'}{m^2}}{e-g+1-d}.
\]
Let $l\in\mathbf{H}_M(I)$.  If $C$ is connected and $(C,\pts)$ is semistable with respect to $l$, then $C_{\red}$ does not have a tacnode. 
\end{proposition}
\begin{proof}
Suppose that $C_{\red}$ has a tacnode at $P$.  Let $\pi:\bar{C}\rightarrow C$ be the normalisation.  There exist two distinct points, $Q_1,Q_2\in\bar{C}$, such that $\pi(Q_1)=\pi(Q_2)=P$.  Moreover the two tangent lines to $C$ at $P$ are coincident.  Define the divisor $D:=Q_1+Q_2$ on $\bar{C}$.

We consider $H^0(\Pro(W),\Osheaf_{\Pro(W)}(1))$ as a subspace of $H^0(p_w(C),L_W)$.  Thus we may define subspaces:
\begin{eqnarray*}
W_0	&:=& \{s\in H^0(\Pro(W),\Osheaf_{\Pro(W)}(1))|\pi_*p_{W*}x \mbox{ vanishes to order $\geq 2$ at $Q_1$ and $Q_2$}\} \\
W_1	&:=& \{s\in H^0(\Pro(W),\Osheaf_{\Pro(W)}(1))|\pi_*p_{W*}x \mbox{ vanishes to order $\geq 1$ at $Q_1$ and $Q_2$}\} .
\end{eqnarray*}
Then $0\subset W_0\subset W_1\subset W_2:=H^0(\Pro(W),\Osheaf_{\Pro(W)}(1))$ is a filtration of $W_2$; choose a basis $w_1,\ldots ,w_{N+1}$ of $W_2$ relative to this filtration.   Write $N_i:=\dim W_i$ for $i=0,1$.  Let $\lambda$ be the 1-PS of $GL(W)$ whose action is given by
\[\begin{array}{rcrcrcccl}
\lambda(t)w_i	&=& w_i,	&t\in\C^*,	& 1	&\leq & i &\leq	& N_0\\
\lambda(t)w_i	&=& tw_i,	&t\in\C^*,	& N_0+1	&\leq & i &\leq & N_1 \\
\lambda(t)w_i	&=& t^2w_i,	&t\in\C^*,	& N_1+1	&\leq & i &\leq & N+1.
\end{array}\]
and let $\lambda'$ be the associated 1-PS of $SL(W)$.  
Fix $(m,\mhat,m')\in M$.

Construct a filtration of $H^0(C,L_W^m\otimes L_r^{\mhat})$ of subspaces of the form $\barO_i$, following (\ref{template_filtration}) as usual.  We wish to estimate $\bhat_i=\dim\barO_i$ in order to apply Lemma \ref{main_calculation}.

As in Proposition \ref{fund_eqn_prop} we use the homomorphism
\[
\pi_{m,\mhat *}:H^0(C,L_W^m\otimes L_r^{\mhat})\rightarrow H^{0}(\bar{C}, \bar{L}_{W}^{m} \otimes \bar{L}_{r}^{\mhat})
\]
induced by the normalisation morphism.  Similarly to as in Proposition \ref{unram_prop} we show that 
\begin{eqnarray}
\pi_{m,\mhat *}(\barO_i) 
	&\subseteq & H^{0}(\bar{C}, \bar{L}_{W}^{m} \otimes \bar{L}_{r}^{\mhat} \otimes \mathcal{O}_{\bar{C}}(-(2m-i)D)), \label{correct_vanishing_tac}
\end{eqnarray}
 for $0 \leq i \leq 2m$.  When $i=0$ this follows from the definitions.  For $1 \leq i \leq 2m$ it is enough to show that for any monomial $\hat{M} \in  \barO_i-\barO_{i-1}$ we have 
\[
 \pi_{m,\mhat *}(\hat{M}) \in H^{0}(\bar{C}, \bar{L}_{W}^{m}\otimes \bar{L}_{r}^{\mhat}\otimes\mathcal{O}_{\bar{C}}(-(2m-i)D)).
\]
Suppose that such $\hat{M}$ has $j_k$ factors from $W_k$, for $k=0,1,2$.  Then $j_0+j_1+j_2=m$ and $j_1+2j_2=i$.  By definition, $\pi_{m,\mhat*}(\hat{M})$ vanishes at both $Q_1$ and $Q_2$ to order at least $2j_0+j_1$.
  But 
\[
  2j_0+j_1 = 2(j_0+j_1+j_2) - (j_1+2j_2) = 2m-i,
\]
giving the required vanishing of $\pi_{m,\mhat*}(\hat{M})$.

We shall also prove at this stage that $\sum_{i=1}^{N+1}w_{\lambda}(w_{i}) \leq 3$.  We have two exact sequences 
\begin{displaymath}
\begin{array}{ccccccccc}  
0 & \rightarrow & \bar{L}_{W}(-Q_1) & \rightarrow & \bar{L}_{W} & \rightarrow & k(Q_1) & \rightarrow & 0 \\
0 & \rightarrow & \bar{L}_{W}(-2Q_1) & \rightarrow & \bar{L}_{W}(-Q_1) & \rightarrow & k(Q_1) & \rightarrow & 0 ,\\
\end{array}
\end{displaymath}
which give rise to long exact sequences in cohomology
 \begin{displaymath}
\begin{array}{ccccccc}  
0 	& \rightarrow & H^{0}(\bar{C},\bar{L}_{W}(-Q_1)) 
	& \rightarrow & H^{0}(\bar{C}, \bar{L}_{W}) 
	& \rightarrow & H^{0}(\bar{C},k(P)) \\
0 	& \rightarrow & H^{0}(\bar{C},\bar{L}_{W}(-2Q_1)) 
	& \rightarrow & H^{0}(\bar{C}, \bar{L}_{W}(-Q_1)) 
	& \rightarrow & H^{0}(\bar{C},k(P)) .\\
	
\end{array}
\end{displaymath}
If we let $\bar{W}_2$ be the image of $W_2=H^0(\Pro(W),\Osheaf_{\Pro(W)}(1))$ in $H^0(\bar{C},\bar{L}_W)$, we may intersect each of the spaces in the sequences with $W_2$.  Note that the image of $W_1$ in $H^0(\bar{C},\bar{L}_W)$ is precisely $\bar{W}_2\cap H^{0}(\bar{C},\bar{L}_{W}(-Q_1))$; if a section of $\bar{L}_W$ which vanishes at $Q_1$ is the pull-back a section of $C$, then it automatically vanishes at $Q_2$ as well.  Similarly the image of $W_0$ in $H^0(\bar{C},\bar{L}_W)$ is $\bar{W}_2\cap H^{0}(\bar{C},\bar{L}_{W}(-2Q_1))$.  We have obtained:
\[
\begin{array}{ccccccc}  
0 	& \rightarrow & \bar{W}_1
	& \rightarrow & \bar{W}_2 
	& \rightarrow & \bar{W}_2\cap H^{0}(\bar{C},k(P)) \\
0 	& \rightarrow & \bar{W}_0
	& \rightarrow & \bar{W}_1
	& \rightarrow & \bar{W}_2\cap H^{0}(\bar{C},k(P)) .
\end{array}
\]
It follows that $\dim W_2/W_1 \leq 1$ and $\dim W_1/W_0 \leq 1$.  Thus $\sum_{i=1}^{N+1}w_{\lambda}(w_{i}) \leq 1+2=3=:R$.

The weight coming from the marked points will always be estimated as zero.  In order to estimate $\bhat_p$, there are unfortunately various cases to consider, which shall give rise to different inequalities:
\begin{enumerate}
\item[1.] there is one irreducible component $C_1$ of $C$ passing through $P$;
\item[2.] there are two irreducible components $C_1$ and $C_2$ of $C$ passing through $P$, and $\deg_{C_{i\,\red}}L_{W}\otimes L_r\geq 2$ for $i=1,2$;
\item[3.] there are two irreducible components $C_1$ and $C_2$ of $C$ passing through $P$, and $\deg_{C_{1\,\red}}L_{W}\otimes L_r = 1$, while $\deg_{C_{2\,\red}}L_{W}\otimes L_r\geq 2$.
\end{enumerate}
There cannot be two degree $1$ curves meeting at a tacnode so these are the only cases.  Note that in Case 1, since $C_1$ is an irreducible curve with a tacnode, $\deg_{C_{1\,\red}}L_{W}\otimes L_r\geq 4$.  For Case 3, we know by Proposition \ref{e'>0} that $\deg_{C_{1\,\red}}L_W\geq 1$ and hence we see that $\deg_{C_{1\,\red}}L_r =0$.  


{\it Cases 1 and 2.}

We estimate $\bhat_p$.  By (\ref{correct_vanishing_tac}) and Riemann-Roch,
\begin{eqnarray*} 
\bhat_{i} :=  \dim \barO_i
		&\leq&   h^{0}(\bar{C},\bar{L}_{W}^{m}\otimes \bar{L}_{r}^{\mhat} \otimes \mathcal{O}_{\bar{C}}(-(2m-i)D)) 
			+ \dim \ker \pi_{m,\mhat * } \\
	&  \leq &  em+d\mhat -2(2m-i)-\bar{g} +1 \\
		&&{}+h^{1}(\bar{C},\bar{L}_{W}^{m} \otimes \bar{L}_{r}^{\mhat} \otimes \mathcal{O}_{\bar{C}}(-(2m-i)D))+\dim \ker \pi_{m,\mhat * }.
\end{eqnarray*}
Claim \ref{I,II,III} allows us to estimate, as usual, the upper bounds $\dim \ker \pi_{m,\mhat *} < q_{2}$ and  $h^{1}\left(\bar{C},\bar{L}_{W}^{m}\otimes \bar{L}_{r}^{\mhat} \otimes \mathcal{O}_{\bar{C}}(-(2m-i)D)\right) \leq 4m-2i\leq 4m$ if $0 \leq i \leq 2g-2$.  To show that the $h^1$ term vanishes for larger values of $i$, recall that our assumptions imply that $\mhat>m$.  Thus, for any component $C_i$ of the curve,
\[
\deg_{\bar{C}_i}\bar{L}_W^m\otimes\bar{L}_r^{\mhat}\geq\deg_{\bar{C}_i}\bar{L}_W^m\otimes\bar{L}_r^m = m\cdot \deg_{\bar{C}_i}\bar{L}_W\otimes\bar{L}_r.
\]
If we are in Case 1 then we conclude that $\deg_{\bar{C}_1}\bar{L}_W^m\otimes\bar{L}_r^{\mhat} \otimes \mathcal{O}_{\bar{C}}(-(2m-i)D)\geq 4m-2\cdot 2m+i\geq 2g-2$ if $2g-2\leq i\leq 2m-1$.  The other components of the curve do not meet $D$ and so one sees from Claim \ref{I,II,III}.III that $h^1$ is zero there.  Case 2 follows similarly.

Combining these inequalities we have
\[
\bhat_{i} \leq \left\{ \begin{array}{ll} (e-4)m+d\mhat+2i-\bar{g}+1+q_2+4m, 	& 0 \leq i \leq 2g-2 \\
					 (e-4)m+d\mhat+2i-\bar{g}+1+q_2, 		& 2g-1 \leq i \leq 3m-1.
\end{array} \right.
\]
In the language of Lemma \ref{main_calculation}, we set $\alpha=4$, $\beta=0$, $\gamma=2$, $\delta=0$, $\epsilon=8g-3\bar{g}+3q_2-1$, $r_N=2$ and $R=3$.  Now:
\begin{eqnarray*}
\lefteqn{\mu^{L_{m,\mhat,m'}}((C,\pts),\lambda') 
	 \geq (4(e-g+1)-3e)m^2 - 3dm\mhat-3nm' } \\
		&&{} - \left((10g-3\bar{g}+3q_2-5)(e-g+1)+3\right)m \\
	&\geq & \left(e-4(g-1) -3d\frac{\mhat}{m}-3n\frac{m'}{m^2}\right)m^2
		{}-\left(10g-3\bar{g}+3q_2\right)(e-g+1)m.
\end{eqnarray*}
This is clearly positive for large $m$.  In particular, as we set $m>(10g-3\bar{g}+3q_2)(e-g+1)$, we must infer that
\[
\mu^{L_{m,\mhat,m'}}((C,\pts),\lambda') 
	\geq  \left(e-4g+3 -3d\frac{\mhat}{m}-3n\frac{m'}{m^2}\right)m^2.
\]
Then, since 
\[
d\frac{\mhat}{m}+n\frac{m'}{m^2} < \frac{1}{8}e-\frac{9}{8}g + \frac{7}{8} < \frac{1}{3}e-\frac{4}{3}g+1
\]
(where the latter inequality may be seen to hold, since $e-g\geq 4$), we conclude that $\mu^{L_{m,\mhat,m'}}((C,\pts),\lambda')>0$.

{\it Case 3.} 

Now that $\deg_{\bar{C}_1}\bar{L}_W\otimes\bar{L}_r=1$, we need a new way to estimate $\bhat_{i} :=  \dim \barO_p$.  However we do know that the genus of $C_1$ is zero.  Let $Y:=\overline{C - C_1}$.  Noting in line (\ref{correct_vanishing_tac}) that $\bar{C_1}$  and $\bar{Y}$ are disjoint, we may write
\begin{eqnarray*}
\lefteqn {\pi_{m,\mhat*}(\barO_p)
	\subset H^{0}(\bar{C}, \bar{L}_{W}^{m}\otimes \bar{L}_{r}^{\mhat}\otimes\mathcal{O}_{\bar{C}}(-(2m-i)D)) }\\
	&&=  H^0(\bar{C}_1,\bar{L}_{W}^m\otimes\bar{L}_r^{\mhat}(-(2m-i)Q_1)) 
		\oplus H^0(\bar{Y},\bar{L}_{W}^m\otimes\bar{L}_r^{\mhat}(-(2m-i)Q_2)).
\end{eqnarray*}
Thus,
\begin{eqnarray} 
\bhat_{i} :=  \dim \barO_i
	&  \leq &  h^0_i +
		(e-1)m+d\mhat -(2m-i)-g_{\bar{C}} +1 \notag\\
		&&{}+h^{1}(\bar{Y},\bar{L}_{W}^{m} \otimes \bar{L}_{r}^{\mhat} \otimes \mathcal{O}_{\bar{Y}}(-(2m-i)Q_2))+\dim \ker \pi_{m,\mhat *},\notag
\end{eqnarray}
where we write $h^0_i$ for $h^0(\bar{C}_1,\bar{L}_{W}^m\otimes\bar{L}_r^{\mhat}\otimes\Osheaf_{\bar{C}_1}(-(2m-i)Q_1))$.

As usual, Claim \ref{I,II,III} allows us to estimate the upper bounds $\dim \ker \pi_{m,\mhat *} < q_{2}$ and  $h^{1}\left(\bar{Y},\bar{L}_{W}^{m}\otimes \bar{L}_{r}^{\mhat} \otimes \mathcal{O}_{\bar{Y}}(-(2m-i)Q_2)\right) \leq 2m-i\leq 2m$ if $0 \leq i \leq 2g-2$.  The assumptions for Case 3 tell us directly that we may apply Claim \ref{I,II,III}.III and obtain $h^{1}\left(\bar{Y},\bar{L}_{W}^{m}\otimes \bar{L}_{r}^{\mhat} \otimes \mathcal{O}_{\bar{Y}}(-(2m-i)Q_2)\right) = 0 $ if $2g-1 \leq i \leq 2m-1$.  

Combining these inequalities we have:
\begin{equation}
\bhat_{i} \leq \left\{ \begin{array}{ll} (e-3)m+d\mhat+i-\bar{g}+1+q_2+2m+h^0_i, 	& 0 \leq i \leq 2g-2 \\
					 (e-3)m+d\mhat+i-\bar{g}+1+q_2+h^0_i, 		& 2g-1 \leq i \leq 2m-1.
\end{array} \right.\label{beta_for_case_2}
\end{equation}
To calculate the $\epsilon$ term in the language of Lemma \ref{main_calculation}, we must calculate $\sum_{i=0}^{2m-1}h^0_i$.  We recall that $\bar{C}_1\cong \Pro^1$, that $\deg_{\bar{C}_{1}}\bar{L}_{W}=1$ and that $\deg_{\bar{C}_{1}}\bar{L}_r=0$.  Thus $\bar{L}_{W\,\bar{C}_1}=\Osheaf_{\Pro^1}(1)$ and $\bar{L}_{r\,\bar{C}_1}=\Osheaf_{\Pro^1}$, so 
\[
h^0_i=h^0(\Pro^1,\Osheaf_{\Pro^1}(m)\otimes\Osheaf_{\Pro^1}(-(2m-i)Q_1))
	=\left\{\begin{array}{ll} 
		0		& i\leq m-1 \\
		-m+i+1		& i\geq m.
		\end{array}\right.
\]
Thus in particular
\[
\sum_{i=0}^{2m-1}h^0_i = \sum_{j=1}^{m}j=\frac{1}{2}m^2+\frac{1}{2}m.\label{sum_case_b}
\]
Hence, we calculate $\epsilon$ to be $=4g-2\bar{g}+2q_2+\frac{1}{2}+\frac{1}{2}m$.
For the rest of the dictionary for Lemma \ref{main_calculation}, we set $\alpha=4$, $\beta=0$, $\gamma=1$, $\delta=0$, $r_N=2$ and $R=3$.  Then
\begin{eqnarray*}
\lefteqn {\mu^{L_{m,\mhat,m'}}((C,\pts),\lambda') 
	\geq  \left(4(e-g+1)-3e\right)m^2 -3dm\mhat -3nm' }\\
		&&{}-\left((6g-2\bar{g}+2q_2-\frac{5}{2}+\frac{1}{2}m)(e-g+1)+3\right)m \\
	&\geq & \left(\frac{1}{2}e-\frac{7}{2}(g-1) -3d\frac{\mhat}{m} -3n\frac{m'}{m^2}\right)m^2
		{}-\left(6g-2\bar{g}+2q_2\right)(e-g+1)m.
\end{eqnarray*}
Our estimates for $m$ imply that $m>(e-g+1)(6g-2\bar{g}+2q_2)$, so we have shown that
\[
\mu^{L_{m,\mhat,m'}}((C,\pts),\lambda') 
	\geq \left(\frac{1}{2}e-\frac{7}{2}g+\frac{5}{2} -3d\frac{\mhat}{m} -3n\frac{m'}{m^2}\right)m^2.
\]
This is clearly positive, as
\[
d\frac{\mhat}{m}+n\frac{m'}{m^2} < \frac{1}{8}e-\frac{9}{8}g + \frac{7}{8} < \frac{1}{6}e-\frac{7}{6}g + \frac{5}{6},
\]
where the second inequality holds since $e-g\geq 4$.

Thus $\mu^{L_{m,\mhat,m'}}((C,\pts),\lambda)>0$ for all triples $(m,\mhat,m')\in M$, and hence by Lemma \ref{ignore_hull} we see that $(C,\pts)$ is unstable with respect to $l$ for any $l\in \mathbf{H}_M(I)$.
\end{proof}


\subsection{Marked points are non-singular and distinct}\label{marked_points}

We now turn to the marked points, which we would like to be non-singular and distinct.  This is ensured in the following two propositions.
\begin{proposition}\label{points_nonsingular}Let $a$ be sufficiently large that $e-9g+7>0$, and let  $M$ consist of those $(m,\mhat,m')$ such that $m,\mhat > m_{3}$ and
\[m > \max \left \{ \begin{array}{c} 
	( g - \frac{1}{2} + e(q_{1}+1) + q_{3} +\mu_{1}m_2)(e-g+1), \\ 
	(10g+3q_{2}-3\bar{g})(e-g+1) \\
\end{array} \right \} 
\]
with
\[
d\frac{\mhat}{m} + n\frac{m'}{m^2} 
	< \frac{1}{8}e-\frac{9}{8}g +\frac{7}{8}
\]
while
\begin{eqnarray*}
\frac{\mhat}{m} 
	&>& 1 + \frac{\frac{3}{2}g-1+d+n\frac{m'}{m^2}}{e-g+1-d}\\
\frac{m'}{m^2}
	&>& \frac{g+d\frac{\mhat}{m}}{e-g+1-n}.
\end{eqnarray*}
Let $l\in\mathbf{H}_M(I)$.  If $(C,x_1,\ldots ,x_n)$ is connected, and semistable with respect to $l$, then all the marked points lie on the non-singular locus of $C$.
\end{proposition}
\emph{Remark.} As $e-g+1-n=(2a-1)(g-1)+(a-1)n+cad$ it is evident that this is positive.

\begin{proof}
Suppose there exists a point $P \in C$, which is singular and also the location of a marked point.  By Proposition \ref{dblpts} this point is a double point.  Let $ev: H^{0}(\Pro(W), \mathcal{O}_{\Pro(W)}(1)) \rightarrow k(P)$ be the evaluation map.  Let $W_{0} = \ker ev$.  We have $N_{0} := \dim W_{0} = N-1$.  Choose a basis $w_0,\ldots ,w_N$ of $W_1:=H^{0}(\Pro(W), \mathcal{O}_{\Pro(W)}(1))$ respecting the filtration $0\subset W_0\subset W_1$.  Let $\lambda$ be the 1-PS of $GL(W)$ whose action is given by
\begin{displaymath}
\begin{array}{lcrcl}
\lambda(t)w_{i} & = & w_{i}, & t \in \C^{*}, & 0  \leq  i  \leq  N-1 \\
\lambda(t)w_{i} & = & tw_{i}, & t \in \C^{*}, & i=N
\end{array} 
\end{displaymath}and let $\lambda'$ be the associated 1-PS of $SL(W)$.  
Fix $(m,\mhat,m')\in M$.

We have assumed that at least one marked point, $x_i$, say, lies at $P$.  If $w\in W_0$ then $w(x_i)=0$, so $w_{k_i}$ must be $w_{e-g}$, whose $\lambda$-weight is $1$.  Hence $\sum_{l=1}^n w_{\lambda}(w_{k_l})m'\geq m'$.
Note also that $\displaystyle \sum_{i=0}^{N}w_{\lambda}(w_{i}) = 1=:R$.
As usual, construct a filtration of $H^{0}(C, L_W^{m}\otimes L_{r}^{\mhat})$ of increasing weight as in (\ref{template_filtration}).  We need to find an upper bound for $\bhat_p=\dim\barO_p$.

Let $\pi: \bar{C} \rightarrow C$ be the normalisation morphism, which is unramified as the hypotheses of Proposition \ref{unram_prop} are satisfied.  There are two distinct points in $\pi^{-1}(P)$, by Proposition \ref{double_points}.  Let the divisor $D:=Q_1+Q_2$ on $\bar{C}$ consist of these points.  Should $Q_1$ and $Q_2$ lie on the same component $\bar{C}_1$ of $\bar{C}$, we see as in the proof of Proposition \ref{unram_prop} that $\deg_{\bar{C}_1}\bar{L}_W\geq 3$.
  
The normalisation morphism induces a homomorphism 
\[ \pi_{m,\mhat*}: H^{0}(C,L_{W}^{m}\otimes L_{r}^{\mhat}) \rightarrow  H^{0}(\bar{C},\bar{L}_{W}^{m}\otimes \bar{L}_{r}^{\mhat}),
\]  
with $\pi_{m,\mhat*}(\barO_p) \subseteq  H^{0}(\bar{C},\bar{L}_{W}^{m} \otimes \bar{L}_{r}^{\mhat} \otimes \mathcal{O}_{\bar{C}}(-(m-p)D)).$  We have
\begin{eqnarray*} 
\bhat_{p} &:= & \dim\barO_p
		\leq h^{0}(\bar{C},\bar{L}_W^m\otimes\bar{L}_r^{\mhat}\otimes\mathcal{O}_{\bar{C}}(-(m-p)D)) 
			+ \dim\ker\pi_{m,\mhat * } \\
	& =&  em+d\mhat -2(m-p)-g_{\bar{C}}+1 \\
		&&{}+ h^{1}(\bar{C},\bar{L}_{W}^{m}\otimes\bar{L}_{r}^{\mhat}\otimes\mathcal{O}_{\bar{C}}(-(m-p)D)) 
		+ \dim \ker \pi_{m,\mhat * } .
\end{eqnarray*}
We may use Claim \ref{I,II,III} as usual to establish that $\dim \ker \pi_{m,\mhat * }<q_{2}$, that $h^{1}(\bar{C},\bar{L}_{W}^{m} \otimes \bar{L}_{r}^{\mhat}(-(m-p)D)) \leq 2(m-p)\leq 2m$ if $0\leq p \leq 2g-2$, and that $h^{1}(\bar{C},\bar{L}_{W}^{m} \otimes \bar{L}_{r}^{\mhat}(-(m-p)D)) =0$ if $p \geq 2g-1$.   

We can now estimate $\bhat_p$:
\begin{displaymath}
\bhat_{p}\leq\left\{ \begin{array}{ll} (e-2)m+d\mhat+2p-\bar{g}+q_2+1+2m, 	& 0 \leq p \leq 2g-2 \\
				 	(e-2)m+d\mhat+2p-\bar{g}+q_2+1, 		& 2g-1 \leq p \leq m.
\end{array} \right.
\end{displaymath}
We may apply Lemma \ref{main_calculation}, setting $\alpha=2$, $\beta=0$, $\gamma=2$, $\delta=1$, $\epsilon=-\bar{g}+q_2+4g-1$, $r_N=1$ and $R=1$; thus we estimate:
\begin{eqnarray*}
\lefteqn{ \mu^{L_{m,\mhat,m'}}((C,\pts),\lambda')
	\geq \left((e-g+1)-e\right)m^2 -dm\mhat +((e-g+1)-n)m' } \\
	&& {}-\left((5g-\bar{g}+q_2-3)(e-g+1)+1\right)m \\
	&\geq& \left(-g+1 -d\frac{\mhat}{m} +(e-g+1-n)\frac{m'}{m^2}\right)m^2 
		{}-\left(5g-\bar{g}+q_2\right)(e-g+1)m.
\end{eqnarray*}
Our assumptions imply that $m>(5g-\bar{g}+q_2)(e-g+1)$ and so we have show that
\[
\mu^{L_{m,\mhat,m'}}((C,\pts),\lambda')
	\geq \left(-g -d\frac{\mhat}{m} +(e-g+1-n)\frac{m'}{m^2}\right)m^2 .
\]
This, however, is positive, as we assumed that
\[
\frac{m'}{m^2} > \frac{g+d\frac{\mhat}{m}}{e-g+1-n}.
\]
Thus $\mu^{L_{m,\mhat,m'}}((C,\pts),\lambda')>0$.  This is true for any $(m,\mhat,m')\in M$, so it follows by Lemma \ref{ignore_hull} that $(C,\pts)$ is not semistable with respect to $l$, for any $l\in\mathbf{H}_M(I)$.
\end{proof}


\begin{proposition}\label{points_distinct}Let  $a$ be sufficiently large that $e-9g+7>0$, and let $M$ consist of those $(m,\mhat,m')$ such that $m,\mhat > m_{3}$ and
\[m > \max \left \{ \begin{array}{c} 
	( g - \frac{1}{2} + e(q_{1}+1) + q_{3} +\mu_{1}m_2)(e-g+1), \\ 
	(10g+3q_{2}-3\bar{g})(e-g+1) 
\end{array} \right \} 
\]
with
\[
d\frac{\mhat}{m} + n\frac{m'}{m^2} 
	< \frac{1}{8}e-\frac{9}{8}g+\frac{7}{8} 
\]
while
\begin{eqnarray*}
\frac{\mhat}{m} 
	&>& 1 + \frac{\frac{3}{2}g-1+d+n\frac{m'}{m^2}}{e-g+1-d}\\
\frac{m'}{m^2}
	&>& \frac{1}{4} + \frac{g+\frac{n}{4}+ d\frac{\mhat}{m}}{2(e-g+1)-n}.
\end{eqnarray*}
Let $l\in\mathbf{H}_M(I)$.  If $(C,x_1,\ldots ,x_n)$ is connected, and semistable with respect to $l$, then all the marked points are distinct.
\end{proposition}

\emph{Remark.} The denominator $2(e-g+1)-n$ is easily checked to be positive, as it is equal to $(4a-1)(g-1)+(2a-1)n+cad$.

\begin{proof}
Suppose two marked points, $x_i$ and $x_j$ meet at $P \in C$.  The hypotheses of the previous proposition hold and so $P$ is a non-singular point.  Let $ev: H^{0}(\Pro(W), \mathcal{O}_{\Pro(W)}(1)) \rightarrow k(P)$ be the evaluation map.  Let $W_{0} = \ker ev$, and let $N_{0} := \dim W_{0} = N$.  Let $w_0,\ldots ,w_{N}$ be a basis of $W_1:=H^0(\Pro(W),\Osheaf_{\Pro(W)}(1))$ respecting the filtration $0\subset W_0\subset W_1$.  Let $\lambda$ be the 1-PS of $GL(W)$ whose action is given by
\begin{displaymath}
\begin{array}{lcrcl}
\lambda(t)w_{i} & = & w_{i}, & t \in \C^{*}, & 0  \leq  i  \leq  N-1 \\
\lambda(t)w_{i} & = & tw_{i}, & t \in \C^{*}, & i=N
\end{array} 
\end{displaymath}
and let $\lambda'$ be the associated 1-PS of $SL(W)$.  
Fix $(m,\mhat,m')\in M$.

As we assume that $x_i$ and $x_j$ lie at $P$, it follows that $\sum_{l=1}^n w_{\lambda}(w_{k_l})m'\geq 2m'$.
Again, note that $\sum_{i=0}^{N}w_{\lambda}(w_{i}) = 1=:R$.  Construct a filtration of $H^{0}(C, L_W^{m}\otimes L_{r}^{\mhat})$ of increasing weight as in (\ref{template_filtration}).  We need to find an upper bound for $\bhat_p:=\dim\barO_p$.

This time, we do not need to use the normalisation to estimate $\bhat_p$; by Proposition \ref{points_nonsingular}, we know that $C$ is smooth at $P$.  It is clear that $\barO_p \subseteq  H^{0}(C,L_{W}^{m} \otimes L_{r}^{\mhat} \otimes \mathcal{O}_{C}(-(m-p)P)).$  We have
\begin{eqnarray*} 
\bhat_{p} 
	& :=  & \dim\barO_p
		\leq   h^{0}(C,L_{W}^{m} \otimes L_{r}^{\mhat} \otimes \mathcal{O}_{C}(-(m-p)P)) \\
	&  = &  em + d\mhat - (m-p) -g+1+ 
		h^{1}({C},{L}_{W}^{m}\otimes {L}_{r}^{\mhat} \otimes \mathcal{O}_{{C}}(-(m-p)P)).
\end{eqnarray*}
We use Claim \ref{I,II,III} to estimate that $h^{1}({C},{L}_{W}^{m} \otimes {L}_{r}^{\mhat}(-(m-p)P)) \leq (m-p)\leq m$ and that $h^{1}({C},{L}_{W}^{m} \otimes {L}_{r}^{\mhat}(-(m-p)P)) =0$ if $p \geq 2g-1$.   

These give us upper bounds for $\bhat_p$:
\begin{displaymath}
\bhat_{p} \leq \left\{ \begin{array}{ll} 
			(e-1)m+d\mhat+p-g+1+m, 	& 0 \leq p \leq 2g-2 \\
 			(e-1)m+d\mhat+p-g+1, 	& 2g-1 \leq p \leq m-1.
\end{array} \right.
\end{displaymath}
We may apply Lemma \ref{main_calculation}, setting $\alpha=1$, $\beta=0$, $\gamma=1$, $\delta=2$, $\epsilon=g-1$, $r_N=1$ and $R=1$.  Thus we estimate:
\begin{eqnarray*}
\lefteqn{ \mu^{L_{\m}}((C,\pts),\lambda')
	\geq \left(\frac{1}{2}(e-g+1)-e\right)m^2 -dm\mhat  }\\
	&& {}+ \left(2(e-g+1)-n\right)m'
		-\left((2g-\frac{5}{2})(e-g+1)+1\right)m \\
	&\geq & \left(-\frac{1}{2}e-\frac{1}{2}g-\frac{1}{2}-d\frac{\mhat}{m}+(2(e-g+1)-n)\frac{m'}{m^2}\right)m^2,
\end{eqnarray*}
where we have used the fact that our assumptions imply $m>(2g-\frac{5}{2})(e-g+1)-(g-1)$.

However, this must be positive, as we assumed that
\[
\frac{m'}{m^2}> \frac{1}{4} + \frac{g+\frac{n}{4}+ d\frac{\mhat}{m}}{2(e-g+1)-n}
	 =\frac{\frac{1}{2}e+\frac{1}{2}g+\frac{1}{2}+d\frac{\mhat}{m}}{2(e-g+1)-n}.
\]
Thus $\mu^{L_{m,\mhat,m'}}((C,\pts),\lambda')>0$, and therefore by Lemma \ref{ignore_hull} we know that $(C,\pts)$ is not semistable with respect to $l$ for any $l\in\mathbf{H}_M(I)$.
\end{proof}


\subsection{GIT semistable curves are reduced and `potentially stable'}\label{reduced_sect}				

The next three results show that, if $(C,\pts)$ is semistable with respect to some $l$ in our range $\in\mathbf{H}_M(I)$ of virtual linearisations, then the curve $C$ is reduced.  We begin with a generalised Clifford's theorem.


 \begin{lemma}[cf.\ \cite{G} page 18] \label{cliff} Let $C$ be a reduced curve with only nodes, and let $L$ be a line bundle generated by global sections which is not trivial on any irreducible component of $C$.  If $H^{1}(C,L) \neq 0$ then there is a connected subcurve $C' \subset C$ such that 
\begin{equation} \label{cliffeqn} h^{0}(C',L) \leq \frac{\deg_{C'}(L)}{2}+1.
\end{equation}
 Furthermore $C' \not\cong \Pro^{1}$.  
\end{lemma}
\begin{proof}  Gieseker proves nearly all of this.  It remains only to show that that $C'$ may be taken to be connected and $C' \not\cong \Pro^{1}$.  Firstly, if (\ref{cliffeqn}) is satisfied by $C'\subset C$ then to is clear that (\ref{cliffeqn}) must be satisfied by some connected component of $C'$.  So assume $C'$ is connected and suppose that $C' \cong \Pro^{1}$.  Now, every line bundle on $\Pro^{1}$ is isomorphic to $\Osheaf_{\Pro^1}(m)$ for some $m \in \Z$.  By hypothesis $L$ is generated by global sections and non-trivial on $C'$; this implies that $m > 0$.  However, combining this with (\ref{cliffeqn}) implies that $m+1= h^{0}(C',L) \leq \frac{m}{2}+1$ which implies $m\leq 0$, a contradiction.
\end{proof}


\begin{lemma}[\cite{G} page 79]
\label{h1=0} Let $a$ be sufficiently large that $e-9g+7>0$, and let $M$ consist of those $(m,\mhat,m')$ such that $m,\mhat> m_{3}$ and
\[m > \max \left \{ \begin{array}{c} 
	( g - \frac{1}{2} + e(q_{1}+1) + q_{3} +\mu_{1}m_2)(e-g+1), \\ 
	(10g+3q_{2}-3\bar{g})(e-g+1)
\end{array} \right \} 
\]
with
\[
d\frac{\mhat}{m} + n\frac{m'}{m^2} 
	< \frac{1}{8}e - \frac{9}{8}g + \frac{7}{8}
\]
while
\begin{eqnarray*}
\frac{\mhat}{m} 
	&>& 1 + \frac{\frac{3}{2}g-1+d+n\frac{m'}{m^2}}{e-g+1-d}\\
\frac{m'}{m^2}
	&>& \frac{1}{4} + \frac{g+\frac{n}{4}+ d\frac{\mhat}{m}}{2(e-g+1)-n}.
\end{eqnarray*} 
Let $l\in\mathbf{H}_M(I)$.  If $(C,\pts)$ is connected and semistable with respect to $l$, then $H^{1}(C_{\red}, L_{W \red}) = 0$.
\end{lemma}


\begin{proof}Since $C_{red}$ is nodal, it has a dualising sheaf $\omega$.  Suppose $H^{1}(C_{\red}, L_{W \red}) \neq 0$.  Then by duality, 
\[
H^{0}(C_{\red},\omega \otimes L_{W \red}^{-1}) \cong H^{1}(C_{\red}, L_{W \red}) \neq 0.  
\]
By Proposition \ref{e'>0}, the line bundle $L_{W \red}$ is not trivial on any component of $C_{red}$.  Then by Lemma \ref{cliff} there is a connected subcurve $C' \not\cong \Pro^1$ of $C_{\red}$ for which $e'>1$ and $h^{0}(C',L_{W\,C'}) \leq \frac{e'}{2}+1$.  

Let $Y:=\overline{C-C'}_{\red}$ and pick a point $P$ on the normalisation $\bar{Y}$ so that $\pi(P)\in C'\cap Y$.  By Proposition \ref{e'>0}, we know that $\deg_{\bar{Y}_j}\bar{L}_W(-P)\geq 0$ for every component $Y_j$ of $Y$.  We may apply Proposition \ref{fund_eqn_prop}, setting $k=1$ and $b=1$.  There exists $(m,\mhat,m')\in M$ satisfying inequality (\ref{FBI}).  Estimate $d'\geq0$, and $n'\geq 0$.  Recall that in this case $S =3g + q_{2}-\bar{g}'+\frac{1}{2}$, and so the hypotheses on $m$ certainly imply that $\frac{S}{m}(e-g+1)\leq \frac{1}{2}$.  We obtain:
\begin{eqnarray*}
e' + \frac{1}{2}= e' + \frac{k}{2} 
	&\leq& \frac{(\frac{e'}{2}+1)e+ d(\frac{e'}{2}+1)\frac{\mhat}{m} 
		+ n(\frac{e'}{2}+1)\frac{m'}{m^2}}{e-g+1} + \frac{1}{2(e-g+1)} \\
\Rightarrow 0 	
	&<& \textstyle {}-\left(e'+\frac{1}{2}\right)(e-g+1) + \left(\frac{1}{2}e'+1\right)\left(e+ d\frac{\mhat}{m}+n\frac{m'}{m^2}\right) 
		+ \frac{1}{2}.
\end{eqnarray*}	
Use the bound $d\frac{\mhat}{m} + n\frac{m'}{m^2} < \frac{1}{8}e - \frac{9}{8}g+\frac{7}{8}$ to obtain
\begin{eqnarray*}
0 	&<& \textstyle {}-\left(e'+\frac{1}{2}\right)(e-g+1)
		+ \left(\frac{1}{2}e'+1\right)\left(\frac{9}{8}e -\frac{9}{8}g+\frac{7}{8}\right) +\frac{1}{2} \\
	&=& \textstyle {}-\left(\frac{7}{16}e'-\frac{5}{8}\right)(e-g)-\frac{9}{16}e'+\frac{7}{8}.
\end{eqnarray*}
Since $e'>1$ we may substitute in $e'\geq 2$:
\[
0	< {}-\frac{1}{4}(e-g)-\frac{1}{4},
\]
a contradiction.  Thus $H^1(C_{\red},L_{W\,\red})=0$.
\end{proof}
We may now, finally, show that our semistable curves are reduced.

\begin{proposition}[cf.\ \cite{G} 1.0.8] \label{reduced_prop} Let $a$ be sufficiently large that $e-9g+7>0$, and let  $M$ consist of those $(m,\mhat,m')$ such that $m,\mhat > m_{3}$ and
\[m > \max \left \{ \begin{array}{c} 
	( g - \frac{1}{2} + e(q_{1}+1) + q_{3} +\mu_{1}m_2)(e-g+1), \\ 
	(10g+3q_{2}-3\bar{g})(e-g+1), 
\end{array} \right \} 
\]
with
\[
d\frac{\mhat}{m} + n\frac{m'}{m^2} 
	< \frac{1}{8}e - \frac{9}{8}g+\frac{7}{8}
\]
while
\begin{eqnarray*}
\frac{\mhat}{m}
	&>& 1 + \frac{\frac{3}{2}g-1+d+n\frac{m'}{m^2}}{e-g+1-d}\\
\frac{m'}{m^2}
	&>& \frac{1}{4} + \frac{g+\frac{n}{4}+ d\frac{\mhat}{m}}{2(e-g+1)-n}.
\end{eqnarray*} 
Let $l\in\mathbf{H}_M(I)$.  If $(C,\pts)$ is connected and semistable with respect to $l$, then $C$ is reduced.  
\end{proposition}

\begin{proof}  Let $\iota: C_{\red} \rightarrow C$ be the canonical inclusion.  The exact sequence of sheaves on $C$ 
\[  0 \rightarrow \mathcal{I}_{C} \otimes L_{W} \rightarrow L_{W} \rightarrow \iota_{*}L_{W \red} \rightarrow 0
\]
gives rise to a long exact sequence in cohomology
\[  \cdots \rightarrow H^{1}(C,\mathcal{I}_{C} \otimes L_{W}) \rightarrow H^{1}(C,L_{W}) \rightarrow H^{1}(C,\iota_{*}L_{W \red}) \rightarrow 0.
\]

Since $C$ is generically reduced, $\mathcal{I}_{C}$ has finite support, hence $H^{1}(C,\mathcal{I}_{C} \otimes L_{W})=0$.  Lemma \ref{h1=0} tells us us that $H^{1}(C, \iota_{*}L_{W \red}) = H^{1}(C_{\red}, L_{W \red}) = 0$.  The exact sequence implies $H^{1}(C,L_{W})=0$ as well.  Next, the map 
\[ H^{0}(\Pro(W), \mathcal{O}_{\Pro(W)}(1)) \rightarrow H^{0}(p_{W}(C)_{\red}, \mathcal{O}_{p_{W}(C)_{\red}}(1)) \rightarrow H^{0}(C_{\red}, L_{W \red})
\]
is injective by Proposition \ref{1.0.2hat}.  Then 
\begin{eqnarray*}
e-g+1 	&=& h^{0}(\Pro(W),\mathcal{O}_{\Pro(W)}(1)) 
		\leq h^{0}(C_{\red}, L_{W \red}) \\ 
	&=& h^{0}(C,L_{W})-h^{0}(C,\mathcal{I}_{C} \otimes L_{W}) 
		= e-g+1-h^{0}(C,\mathcal{I}_{C} \otimes L_{W}).  
\end{eqnarray*}
	Therefore $h^{0}(C,\mathcal{I}_{C} \otimes L_{W}) =0$.  Since $\mathcal{I}_{C} \otimes L_{W}$ has finite support, $\mathcal{I}_{C}=0$, so $C$ is reduced.  
\end{proof}

Next, we improve on Proposition \ref{fund_eqn_prop}.  If $C$ is connected and $(C,\pts)$ is semistable with respect to some $l\in\mathbf{H}_M(I)$, and if $C'$ is a connected subcurve of $C$, then we know that our `fundamental inequality' is satisfied without needing to verify condition (ii).  This inequality is then an extra property of semistable curves, which we will use in Theorem \ref{jssclosed} to show that $\bar{J}^{ss}(l)\subset J$ for the right range of virtual linearisations.

We repeat the definition of $\lambda'_{C'}$: Let
\[ W_{0} := \ker \{ H^{0}(\Pro(W),\mathcal{O}_{\Pro(W)}(1)) \rightarrow H^{0}(p_{W}(C'), \mathcal{O}_{p_{W}(C')}(1)) \}.
\]
Choose a basis $w_{0},\ldots ,w_{N}$ of $H^{0}(\Pro(W), \mathcal{O}_{\Pro(W)}(1))$ relative to the filtration \\ $0 \subset W_{0} \subset W_{1} = H^{0}(\Pro(W), \mathcal{O}_{\Pro(W)}(1))$.  Let $\lambda_{C'}$ be the 1-PS of $GL(W)$ whose action is given by
\begin{displaymath}
\begin{array}{lcrcccccc}
\lambda_{C'}(t)w_{i} & = & w_{i}, & t \in \C^{*}, & 0 & \leq & i & \leq & N_{0}-1 \\
\lambda_{C'}(t)w_{i} & = & tw_{i}, & t \in \C^{*}, & N_{0} & \leq & i & \leq & N
\end{array} 
\end{displaymath}
and let $\lambda'_{C'}$\index{$\lambda'_{C'}$} be the associated 1-PS of $SL(W)$.  It is more convenient to prove that the inequality holds for linearisations $L_{m,\mhat,m'}$ before inferring the result in general.


\begin{lemma}[cf.\ \cite{G} page 83 and Proposition \ref{fund_eqn_prop} above] \label{amp_lemma} Let $a$ be sufficiently large that $e-9g+7>0$, and suppose that $m,\mhat > m_{3}$ and
\[m > \max \left \{ \begin{array}{c} 
	( g - \frac{1}{2} + e(q_{1}+1) + q_{3} +\mu_{1}m_2)(e-g+1), \\ 
	(10g+3q_{2}-3\bar{g})(e-g+1), 
\end{array} \right \} 
\]
with
\[
d\frac{\mhat}{m} + n\frac{m'}{m^2} 
	< \frac{1}{8}e - \frac{9}{8}g+\frac{7}{8}
\]
while
\begin{eqnarray*}
\frac{\mhat}{m} 
	&>& 1 + \frac{\frac{3}{2}g-1+d+n\frac{m'}{m^2}}{e-g+1-d}\\
\frac{m'}{m^2}
	&>& \frac{1}{4} + \frac{g+\frac{n}{4}+ d\frac{\mhat}{m}}{2(e-g+1)-n}.
\end{eqnarray*}
Let $(C,\pts)\subset\pwpr$ be a connected curve whose only singularities are nodes, and such that no irreducible component of $C$ collapses under projection to $\Pro(W)$.  Suppose $C$ has at least two irreducible components.  Let $C' \neq C$ be a reduced, connected, complete subcurve of $C$ and let $Y$ be the closure of $C-C'$ in $C$ with the reduced structure.  Suppose there exist points $P_{1},\ldots,P_{k}$ on $\bar{Y}$, the normalisation of $Y$, satisfying $\pi ( P_{i}) \in Y \cap C'$ for all $1 \leq i \leq k$.  Write $h^{0}(p_{W}(C'), \mathcal{O}_{p_{W}(C')}(1)) =: h^{0}.$  

Finally, suppose that
\[
\mu^{L_{m,\mhat,m'}}((C,\pts),\lambda'_{C'})\leq 0.  
\]
Then 
\begin{equation}
e' + \frac{k}{2} < \frac{h^{0}e+ (dh^{0}-d'(e-g+1))\frac{\mhat}{m} + (nh^0-n'(e-g+1))\frac{m'}{m^2}}{e-g+1} \index{fundamental inequality}
	+ \frac{S}{m}, \label{fundamental}
\end{equation}
where  $S =g + k(2g-\frac{3}{2}) + q_{2} - \bar{g}+1$.
\end{lemma}
\begin{proof}
Arguing similarly to \cite{G} pages 83-85, we prove the result by contradiction.  We first assume that $k=\#(Y\cap C')$ and then show that this implies the general case.

Let $C'$ be a connected subcurve of $C$, and let $P_1,\ldots ,P_k$ be all the points on $\bar{Y}$ satisfying $\pi ( P_{i}) \in Y \cap C'$. We assume that (\ref{fundamental}) is not satisfied for $C'$, and further that $C'$ is maximal with this property.  Namely, if $C''$ is complete and connected, and $C'\subsetneq C'' \subset C$, then (\ref{fundamental}) does hold for $C''$.  Since (\ref{fundamental}) does not hold for $C'$, 
\begin{eqnarray}
\lefteqn{(e' + \frac{k}{2})(e-g+1) 
	\geq (e'-g'+1)e+ \left(d(e'-g'+1)-d'(e-g+1)\right)\frac{\mhat}{m} }\notag\\
		&\hspace{.5in}&{}+ \left(n(e'-g'+1)-n'(e-g+1)\right)\frac{m'}{m^2} + \frac{S'}{m}(e-g+1). \hspace{.5in}\label{assum}
\end{eqnarray}
As all other hypotheses of Proposition \ref{fund_eqn_prop} have been met, we must conclude that condition (ii) there fails.  Thus there is some irreducible component $\bar{Y}_j$ of $\bar{Y}$, the normalisation of $Y$, such that 
\[
\deg_{\bar{Y}_{j}}(\bar{L}_{W \, \bar{Y}}(-(P_1+\cdots +P_k))) < 0. 
\]
Let $Y_j$ be the corresponding irreducible component of $Y$.  By assumption, $Y_j$ does not collapse under projection to $\Pro(W)$, and so $\deg_{\bar{Y}_{j}}(\bar{L}_{W \, \bar{Y}})=\deg_{Y_j}(L_W)>0$.  Putting this together,
\[
0<\deg_{\bar{Y}_{j}}(\bar{L}_{W \, \bar{Y}})<\#(\bar{Y}_{j}\cap \{P_1,\ldots ,P_k\}) = \#(Y_j\cap C') =:i_{Y_j,C'}.
\]
Thus, $i_{Y_j,C'}\geq 2$.  We define $D$ to be the connected subcurve $Y_j\cup C'$.  By the maximality assumption on $C'$, it follows that (\ref{fundamental}) does hold for $D$. We define constants $e_D$, $k_D$, $h^0_D$, $d_D$, $n_D$ and $S_D$ in analogy with the constants $e'$, $k$, $h^0$, $d'$, $n'$ and $S'$ for $C'$.  Similarly we define constants pertaining to $Y_j$. Then:
\begin{equation}
e_D + \frac{k_D}{2} < \frac{h_D^{0}e+ (dh_D^{0}-d_D(e-g+1))\frac{\mhat}{m} + (nh_D^0-n_D(e-g+1))\frac{m'}{m^2}}{e-g+1} 
	+ \frac{S_D}{m}. \label{first}
\end{equation}
Observe that $e_D=e'+e_{Y_j}$, $d_D=d'+d_{Y_j}$ and $n_D=n'+n_{Y_j}$.  The curve $C$ is nodal, so we conclude:
\begin{enumerate}
\item[(i)] $k_D=\#((Y_j\cup C')\cap(\overline{Y-Y_j}))=\#(C'\cap Y) + \#(Y_j\cap \overline{Y-Y_j}) - \# (C'\cap Y_j)$; if we set $i_{Y_j,Y}:=\#(Y_j\cap Y)$ then $k_D = k + i_{Y_j,Y}-i_{Y_j,C'}$;
\item[(ii)] $g_D=g' + g_{Y_j} + i_{Y_j,C'} -1$;
\item[(iii)] $h_D^0=e_D-g_D+1=e'+e_{Y_j}-g'-g_{Y_j}-i_{Y_j,C'}+2$.
\end{enumerate}
Note in particular that, since $i_{Y_j,C'}\geq 2$ and since $C'$ and $Y_j$ are connected, (ii) implies that $g_D\geq 1$.  But $g_D\leq g$ so if $g=0$ then we already have the required contradiction.  We henceforth assume that $g\geq 1$.  Equation (\ref{first}) may be rearranged to form
\begin{eqnarray*}
\lefteqn{\textstyle \left(e'+e_{Y_j}+\frac{k+i_{Y_j,Y}-i_{Y_j,C'}}{2}\right)(e-g+1) 
	< (e'+e_{Y_j}-g'-g_{Y_j}-i_{Y_j,C'}+2)e }\\
	&\hspace{.5in}&\textstyle	{}+ (d(e'+e_{Y_j}-g'-g_{Y_j}-i_{Y_j,C'}+2)-(d'+d_{Y_j})(e-g+1))\frac{\mhat}{m} \hspace{.5in}\\
	&&\textstyle 	{}+ (n(e'+e_{Y_j}-g'-g_{Y_j}-i_{Y_j,C'}+2)- (n'+n_{Y_j})(e-g+1))\frac{m'}{m^2} \\
	&&\textstyle	{}+ \frac{S_D}{m}(e-g+1).
\end{eqnarray*}
We subtract our assumption, line (\ref{assum}): 
\begin{eqnarray}
\lefteqn{ \textstyle \left(e_{Y_j}+\frac{i_{Y_j,Y}-i_{Y_j,C'}}{2}\right)(e-g+1)
	<(e_{Y_j}-g_{Y_j}-i_{Y_j,C'}+1)e }\notag\\
	&\hspace{.5in}&\textstyle{}+(d(e_{Y_j}-g_{Y_j}-i_{Y_j,C'}+1)-d_{Y_j}(e-g+1))\frac{\mhat}{m} \hspace{.5in}\notag\\
	&&\textstyle		{}+(n(e_{Y_j}-g_{Y_j}-i_{Y_j,C'}+1)-n_{Y_j}(e-g+1))\frac{m'}{m^2} \notag\\
	&&\textstyle		{}+  \frac{(i_{Y_j,Y}-i_{Y_j,C'})(2g-\frac{1}{2})}{m}(e-g+1). \label{subs_here}
\end{eqnarray}	
We rearrange, and use the inequality $e_{Y_j}\leq i_{Y_j,C'}-1$.
\begin{eqnarray}
\lefteqn{\textstyle \left(\frac{i_{Y_j,Y}+i_{Y_j,C'}}{2} + g_{Y_j} - 1 + d_{Y_j}\frac{\mhat}{m} + n_{Y_j}\frac{m'}{m^2} + 
		\frac{(i_{Y_j,C'}-i_{Y_j,Y})(2g-\frac{1}{2})}{m}\right)e } \notag\\
	&<& \textstyle \left(\frac{i_{Y_j,Y}+i_{Y_j,C'}}{2}-1+ d_{Y_j}\frac{\mhat}{m} + n_{Y_j}\frac{m'}{m^2}
				+ \frac{(i_{Y_j,C'}-i_{Y_j,Y})(2g-\frac{1}{2})}{m}\right)(g-1) \notag\\
	&&\textstyle	{}- g_{Y_j}\left(d\frac{\mhat}{m}+n\frac{m'}{m^2}\right), \hspace{.5in}
\end{eqnarray}
so that finally we may estimate
\begin{eqnarray}
\lefteqn{\textstyle\left(\frac{i_{Y_j,Y}+i_{Y_j,C'}}{2} + g_{Y_j} - 1 + d_{Y_j}\frac{\mhat}{m} + n_{Y_j}\frac{m'}{m^2} + 
		\frac{(i_{Y_j,C'}-i_{Y_j,Y})(2g-\frac{1}{2})}{m}\right)(e-g+1) }\notag\\
	&\hspace{1.5in}<&\textstyle {}- g_{Y_j}\left(d\frac{\mhat}{m}+n\frac{m'}{m^2}\right)\leq 0. \hspace{1.5in}	\label{estimate_e}
\end{eqnarray}
Recall that $e-g+1=\dim W>0$, that $i_{Y_j,C'}\geq 2$ and that $g\geq 1$.  Thus the left hand side of (\ref{estimate_e}) is strictly positive.
This is a contradiction.  No such $C'$ exists, i.e.\ all subcurves $C'$ of $C$ satisfy inequality (\ref{fundamental}), provided that $k=\#(C'\cap Y)$.

Finally, suppose that we choose any $k$ points $P_1\ldots ,P_{k}$ on $\bar{Y}$ such that $\pi(Y_i)\in (C'\cap Y)$.  Then $k\leq\#(C'\cap Y):=k'$.  We proved that (\ref{fundamental}) is true for $k'$, and though we must take a little care with the dependence of $S$ on $k$, it follows that (\ref{fundamental}) is true for $k$.
\end{proof}
Now we may extend this to general $l\in\mathbf{H}_M(I)$, to provide the promised extension of the fundamental basic inequality.


\begin{amplification}\label{amp} Let  $a$ be sufficiently large that $e-9g+7>0$, and let $M\subset \tilde{M}$, where $\tilde{M}$ consists of those $(m,\mhat,m')$ such that $m,\mhat > m_{3}$ and
\[m > \max \left \{ \begin{array}{c} 
	( g - \frac{1}{2} + e(q_{1}+1) + q_{3} +\mu_{1}m_2)(e-g+1), \\ 
	(10g+3q_{2}-3\bar{g})(e-g+1), 
\end{array} \right \} 
\]
with
\[
d\frac{\mhat}{m} + n\frac{m'}{m^2} 
	< \frac{1}{8}e - \frac{9}{8}g+\frac{7}{8}
\]
while
\begin{eqnarray*}
\frac{\mhat}{m} 
	&>& 1 + \frac{\frac{3}{2}g-1+d+n\frac{m'}{m^2}}{e-g+1-d}\\
\frac{m'}{m^2}
	&>& \frac{1}{4} + \frac{g+\frac{n}{4}+ d\frac{\mhat}{m}}{2(e-g+1)-n}.
\end{eqnarray*}
Let $l\in\mathbf{H}_{M}(I)$.  Let $(C,\pts)$ be semistable with respect to $l$, where $C$ is a connected curve.  Suppose $C$ has at least two irreducible components.  
Let $C' \neq C$ be a reduced, complete subcurve of $C$ and let $Y:=\overline{C-C'}$.  The subcurves $C'$ and $Y$ need not be connected; suppose $C$ has $b$ connected components.   Suppose there exist points $P_{1},\ldots,P_{k}$ on $\bar{Y}$, the normalisation of $Y$, satisfying $\pi ( P_{i}) \in Y \cap C'$ for all $1 \leq i \leq k$.  Write $h^{0}(p_{W}(C'), \mathcal{O}_{p_{W}(C')}(1)) =: h^{0}.$  Then there exists a triple $(m,\mhat,m')\in M$, such that
\begin{equation}
e' + \frac{k}{2} < \frac{h^{0}e+ (dh^{0}-d'(e-g+1))\frac{\mhat}{m} + (nh^0-n'(e-g+1))\frac{m'}{m^2}}{e-g+1} 
	+ \frac{bS}{m} \label{fund2}\index{fundamental inequality}
\end{equation}
where  $S =g + k(2g-\frac{3}{2}) + q_{2} - \bar{g}+1$.
\end{amplification}

\begin{proof}
First assume that $C'$ is connected, and suppose that inequality (\ref{fund2}) fails for all $(m,\mhat,m')\in M$.  It must follow that $(C,\pts)$ does not satisfy the hypotheses of Lemma \ref{amp_lemma}.  However, as $(C,\pts)$ is semistable with respect to $l$, all the other hypotheses of that lemma are verified, so we must conclude that $\mu^{L_{m,\mhat,m'}}((C,\pts),\lambda'_{C'})>0$ for all $(m,\mhat,m')\in M$.  It follows by Lemma \ref{ignore_hull} that $(C,\pts)$ is unstable with respect to $l$.  The contradiction implies that there do indeed exist some $(m,\mhat,m')\in M$ such that (\ref{fund2}) is satisfied.

Now, let $C'_1,\ldots,C'_b$ be the connected components of $C'$.  We may prove a version of (\ref{fund2}) for each $C_i$, for $i=1,\ldots,b$.  When we sum these inequalities over $i$, it follows that
\[
e' + \frac{k}{2} < \frac{h^{0}e+ (dh^{0}-d'(e-g+1))\frac{\mhat}{m} + (nh^0-n'(e-g+1))\frac{m'}{m^2}}{e-g+1} 
	+ \frac{bS}{m}.
\]
\end{proof}

We summarise the results of Sections \ref{first_props} to \ref{reduced_sect}.  Recall again that since $e$ is defined to be $a(2g-2+n+cd)$, and since $2g-2+n+cd$ is always at least $1$, the denominators $e-g+1-d$ and $2(e-g+1)-n$ are both positive.
\begin{theorem}\label{pot_stab_without_lines}
Let $a$ be sufficiently large that $e-9g+7>0$, and let  $M\subset \tilde{M}$, where $\tilde{M}$ consists of those $(m,\mhat,m')$ such that $m,\mhat > m_{3}$ and
\[m > \max \left \{ \begin{array}{c} 
	( g - \frac{1}{2} + e(q_{1}+1) + q_{3} +\mu_{1}m_2)(e-g+1), \\ 
	(10g+3q_{2}-3\bar{g})(e-g+1), 
\end{array} \right \} 
\]
with
\[
d\frac{\mhat}{m} + n\frac{m'}{m^2} 
	< \frac{1}{8}e - \frac{9}{8}g+\frac{7}{8}
\]
while
\begin{eqnarray*}
\frac{\mhat}{m} 
	&>& 1 + \frac{\frac{3}{2}g-1+d+n\frac{m'}{m^2}}{e-g+1-d}\\
\frac{m'}{m^2}
	&>& \frac{1}{4} + \frac{g+\frac{n}{4}+ d\frac{\mhat}{m}}{2(e-g+1)-n}.
\end{eqnarray*}
Let $l\in\mathbf{H}_M(I)$ and let $(C,\pts)$ be a connected curve, semistable with respect to $l$.  Then $(C,\pts)$ satisfies:
\begin{enumerate}
\item[(i)]  $(C,x_1,\ldots ,x_n)$ is a reduced, connected, nodal curve, and the marked points are distinct and nonsingular;
\item[(ii)] the map $C \rightarrow \Pro(W)$ collapses no component of $C$, and induces an injective map 
\[ H^{0}(\Pro(W), \mathcal{O}_{\Pro(W)}(1)) \rightarrow H^{0}(C, L_{W}) ;
\]
\item[(iii)]  $h^{1}(C, L_{W})=0$;
\item[(iv)] any complete subcurve $C' \subset C$ with $C' \neq C$ satisfies the inequality 
\[ e' + \frac{k}{2} < \frac{h^{0}e+ (dh^{0}-d'(e-g+1))\frac{\mhat}{m} + (nh^0-n'(e-g+1))\frac{m'}{m^2}}{e-g+1} + \frac{bS}{m}
\]
of Amplification \ref{amp}, where $C'$ consists of $b$ connected components and $S =g + k(2g-\frac{3}{2}) + q_{2} - \bar{g}+1$. \hfill$\Box$
\end{enumerate}
\end{theorem}
\begin{definition}\label{potentially_stable_map}
If $(C,\pts)\subset\pwpr$ satisfies conditions (i)-(iv) of Theorem \ref{pot_stab_without_lines}, then the corresponding map $(C,\pts)\stackrel{p_r}{\rightarrow}\pr$ is referred to as a \emph{potentially stable map}.
\end{definition}
{\it Remark.} Gieseker defines his `potentially stable curves' (which have no marked points) using the analogous statements, and the additional condition if the curve is not a moduli stable curve then destabilising components must have two nodes and be embedded as lines.  A similar condition can be given here, and shown to be a corollary of the fundamental inequality (for a restricted range $M$).  

Namely, for a certain $M$, we can show that if $l\in\mathbf{H}_M(I)$ and $(C,\pts)$ is semistable with respect to $l$, and if $C'$ is a rational component of $C$ which is collapsed under projection to $\pr$, then $C'$ has at least two special points; if it has precisely two then it is embedded in $\pw$ as a line.  The proof of this follows \cite{G} Proposition 1.0.9 and full details may be seen in \cite{e_thesis} Corollary 5.5.1, but it has been omitted here for brevity, as it is not needed to prove Theorem \ref{maps}.


\subsection{GIT Semistable maps represented in $\bar{J}$ are moduli stable}\label{prop_and_proj}

In the previous sections we have been studying $I^{ss}$.  In this section we focus on $\bar{J}^{ss}$.   Recall the definitions of $I$ and $J$, given at \ref{i} and \ref{j}: the scheme $I$ is the Hilbert scheme of $n$-pointed curves in $\pwpr$, and $J\subset I$ is the locally closed subscheme such that for each $(h,\pts) \in J$:
\begin{enumerate}
\item[(i)] $(\mathcal{C}_{h},\pts)$ is prestable\index{curve!prestable}, i.e.\ $\mathcal{C}_h$ is projective, connected, reduced and nodal, and the marked points are distinct and non-singular;
\item[(ii)] the projection map $\mathcal{C}_{h} \rightarrow \Pro(W)$ is a non-degenerate embedding;
\item[(iii)] the invertible sheaves $ (\mathcal{O}_{\Pro(W)}(1) \otimes \mathcal{O}_{\Pro^{r}}(1))|_{\mathcal{C}_{h}} $ and $(\omega_{\mathcal{C}_{h}}^{\otimes a}(ax_1+\cdots +ax_n) \otimes \mathcal{O}_{\Pro^{r}}(ca+1))|_{\mathcal{C}_{h}} $ are isomorphic, where $c$ is a positive integer; see note \ref{note_on_c}
\end{enumerate}
Moreover, recall from the discussion at the end of Section \ref{strategy} that we had set out to find a linearisation such that  $\bar{J}^{ss}(L)\subset J$. This, together with non-emptiness of $\bar{J}^{ss}(L)$, is sufficient to show that $\bar{J}\dblq_L SL(W)\cong \Maps$ (Theorem \ref{wishful_thinking}).
 
In this section we find a range $M$ of $(m,\mhat,m')$ such that $\bar{J}^{ss}(l)\subset J$ when $l\in\mathbf{H}_M(I)$.  This range is much narrower than those we have considered so far.

Here is the result we have been seeking:
\begin{theorem}[cf.\ \cite{HM} 4.55] \label{jssclosed} 
Let $M$ consist of those $(m,\mhat,m')$ such that $m,\mhat > m_{3}$ and
\[m > \max \left \{ \begin{array}{c} 
	( g - \frac{1}{2} + e(q_{1}+1) + q_{3} +\mu_{1}m_2)(e-g+1), \\ 
	(10g+3q_{2}-3\bar{g})(e-g+1) \\
	(6g+2q_2-2\bar{g}-1)(2a-1)
\end{array} \right \} 
\]
with
\begin{eqnarray}
\frac{\mhat}{m} &=&\frac{ca}{2a-1} + \delta \label{delta}\\
\frac{m'}{m^2}	&=&\frac{a}{2a-1} + \eta \label{eta}
\end{eqnarray}
where
\begin{equation}
| n\eta | + | d\delta | \leq \frac{1}{4a-2}-\frac{3g+q_2-\bar{g}-\frac{1}{2}}{m}. \label{delta and eta}
\end{equation}
In addition, ensure that $a$ is sufficiently large that
\begin{equation}
d\frac{\mhat}{m}+n\frac{m'}{m^2} < \frac{1}{8}e - \frac{9}{8}g + \frac{7}{8}. \label{old one}
\end{equation}
Let $l\in\mathbf{H}_M(\bar{J})$.  Then $\bar{J}^{ss}(l)$ is contained in $J$.
\end{theorem}
{\it Remark.}  The final assumption on the magnitude of $m$ ensures that the right hand side of (\ref{delta and eta}) is positive and so (\ref{delta and eta}) may be satisfied.  There may seem to be many competing bounds on the ratios $\frac{\mhat}{m}$ and $\frac{m'}{m^2}$, and on $a$.  However, one may show that (\ref{old one}) is implied by (\ref{delta}), (\ref{eta}) and (\ref{delta and eta}) for all $g$,$n$ and $d$ as long as $a\geq 10$ (cf.\ \cite{e_thesis}, proof of Theorem 5.6.1). Smaller values of $a$ are possible for most $g$, $n$ and $d$.  Once a large enough $a$ has been chosen, it is always possible to satisfy the rest of the inequalities; the simplest way is to set $\delta$ and $\eta$ to zero and pick large $m,\mhat,m'$ with the desired ratios.  
\begin{proof}
The range $M$ here is contained within the range for Theorem \ref{pot_stab_without_lines}.  
We may thus apply Theorem \ref{pot_stab_without_lines}:  if $(h,\pts)\in I^{ss}(l')$ for some $l'\in\mathbf{H}_M(I)$ then $(\mathcal{C}_h,\pts)$ nodal and reduced, and one can find $(m,\mhat,m')\in M$ satisfying inequality (\ref{fund2}).  However, this theorem in fact deals with $\mathbf{H}_M(\bar{J})$ and not $\mathbf{H}_M(I)$; on the other hand any $l\in\mathbf{H}_M(\bar{J})$ may be regarded as the restriction of some $l'\in\mathbf{H}_M(I)$ to $\bar{J}$, and then $\bar{J}^{ss}(l)=\bar{J}\cap I^{ss}(l')$.  Thus we may use all our previous results.

Suppose we can show that $J\cap\bar{J}^{ss}(l)$ is closed in $\bar{J}^{ss}(l)$.  Then if $x\in\bar{J}^{ss}(l)-J\cap\bar{J}^{ss}(l)$, there must be an open neighbourhood of $x$ in $\bar{J}^{ss}(l)-J\cap\bar{J}^{ss}(l)$, but this is a contradiction as $x$ is in $\bar{J}$ so $x$ is a limit point of $J$. It follows that $J\cap\bar{J}^{ss}(l)=\bar{J}^{ss}(l)$, i.e.\ that $\bar{J}^{ss}(l)\subset J$.  

We shall proceed by using the valuative criterion of properness to show that the inclusion $J\cap\bar{J}^{ss}(l)\hookrightarrow \bar{J}^{ss}(l)$ is proper, whence $J\cap\bar{J}^{ss}(l)$ is closed in $\bar{J}^{ss}(l)$, as required.  Let $R$ be a discrete valuation ring, with generic point $\xi$ and closed point $0$.  Let $\alpha: \Spec{R} \rightarrow \bar{J}^{ss}(l)$ be a morphism such that $\alpha(\xi) \subset J\cap\bar{J}^{ss}(l)$.  Then we will show that $\alpha(0) \in J\cap\bar{J}^{ss}(l)$.   

Define a family $\mathcal{D}$ of $n$-pointed curves in $\pwpr$ by the following pullback diagram:
\begin{displaymath}
\begin{array}{ccc}   \mathcal{D} 		& \rightarrow 			& \tilde{\mathcal{C}}|_{\bar{J}^{ss}(l)} \\
		\downarrow \uparrow \sigma_i	& 				& \downarrow \uparrow \sigma_i \\
			\Spec R 		& \stackrel{\alpha}{\rightarrow} & \bar{J}^{ss}(l)
\end{array}
\end{displaymath}  
where $\sigma_1,\ldots ,\sigma_n:\Spec R \rightarrow \mathcal{D}$ are the sections giving the marked points.  The images of the $\sigma_i$ in $\mathcal{D}$ are divisors, denoted $\sigma_i(\Spec R)$.  By definition of $J$ we have 
\[ 
\left(\mathcal{O}_{\Pro(W)}(1)\otimes \mathcal{O}_{\Pro^{r}}(1)\right)|_{\mathcal{D}_{\xi}} 
	\cong \acanon_{\mathcal{D}_\xi}(a\sigma_1(\xi)+\cdots +a\sigma_n(\xi))\otimes \mathcal{O}_{\Pro^{r}}(ca+1)|_{\mathcal{D}_{\xi}}.    
\]

We will write $(\mathcal{D}_{0},\sigma_1(0),\ldots ,\sigma_n(0)) =: (C,\pts)$, and show that its representative in the universal family $I$ is in fact in $J$.  The curve $C$ is connected, as a limit of connected curves.  We assumed that $\alpha (0)\in \bar{J}^{ss}(l)=\bar{J}\cap I^{ss}(l')$, where $l'\in\mathbf{H}_M(I)$, and so $(C,\pts)$ satisfies conditions (i) above, and the curve $p_W(C)\subset\pw$ is non-degenerate.  We will show that the line bundles in condition (iii) are isomorphic.  It follows from this that
\begin{equation}
\opw(1)|_C\cong \omega_C^a(ax_1+\cdots+ax_n)\otimes\opr(ca)|_C \label{line_bundles_for_ample}
\end{equation}
and so this line bundle has positive degree on every component of $C$.  However, we know that $\mathcal{L}=\omega_C(x_1+\cdots+x_n)\otimes\opr(c)$ has positive degree on each component of $C'$ if and only if $p_r:C\rightarrow \pw$ is a stable map, and that when this true then $\mathcal{L}^a$ is very ample.  Hence $\mathcal{L}^a$ embeds $C$ in $\pw$, i.e.\ $p_W:C\cong p_W(C)$, and we have verified condition (ii) in the definition of $J$.  Thus checking condition (iii) is sufficient to show that $(C,\pts)$ is represented in $J$.

Decompose $C = \bigcup C_{i}$ into its irreducible components.  Then we can write
\begin{eqnarray*} 
\lefteqn{\left(\mathcal{O}_{\Pro{W}}(1) \otimes \mathcal{O}_{\Pro^{r}}(1)\right)|_{\mathcal{D}} 
	 \cong \acanon_{\mathcal{D}/\Spec R}\left(a\sigma_1(\Spec R)+\cdots +a\sigma_n(\Spec R)\right) } \notag\\
	&\hspace{2in}&\otimes \mathcal{O}_{\Pro^{r}}(ca+1)|_{\mathcal{D}} \otimes \mathcal{O}_{\mathcal{D}}(\sum a_{i}C_{i}) \hspace{.5in}
\end{eqnarray*}
where the $a_{i}$ are integers.  As $\mathcal{O}_{\mathcal{D}}$-modules, $\mathcal{O}_{\mathcal{D}}(C) \cong \mathcal{O}_{\mathcal{D}}$ so we can normalise the integers $a_{i}$ so that they are all non-negative and at least one of them is zero. Separate $C$ into two sub-curves $Y := \displaystyle \bigcup_{a_{i}=0} C_{i}$ and $C' := \displaystyle \bigcup_{a_{i}>0} C_{i}$.  Since at least one of the $a_{i}$ is zero, we have $Y \neq \emptyset$ and $C' \neq C$.  Suppose for a contradiction that $C' \neq \emptyset$ (hence $Y \neq C$).  Let $k = \# (Y \cap C')$ and let $b$ be the number of connected components of $C'$.  Since C is connected, we must have $k\geq b$.  We will obtain our contradiction by showing that $\frac{k}{b}<1$.

Any local equation for the divisor $\mathcal{O}_{\mathcal{D}}(\sum a_{i}C_{i})$ must vanish identically on every component of $C'$ and on no component of $Y$.  Such an equation is zero therefore at each of the $k$ nodes in $Y \cap C'$.    Thus we obtain the inequality 
\begin{eqnarray} k & \leq & \deg_{Y} \left(\mathcal{O}_{\mathcal{D}}(-\sum a_{i}C_{i}) \right) \nonumber\\
	& = & \deg_{Y} \left(\mathcal{O}_{\Pro(W)}(1) \otimes \omega^{\otimes -a}_{\mathcal{D}_0} (-a\sigma_1(0)-\cdots -a\sigma_n(0)) 
			\otimes \mathcal{O}_{\Pro^{r}}(-ca)\right) \nonumber\\
	& = & e_{Y} -a(2g_{Y}-2+n_Y+k)-cad_{Y}.  \notag
\end{eqnarray}
 Substituting $e'=e-e_{Y}$, $d'=d-d_{Y}$, $g' = g-g_{Y}-k+1$ and $e=a(2g-2+n+cd)$, this is equivalent to
\begin{equation} 
e'-a(2g'-2+n'+cd') \leq (a-1)k. \label{ak-k}
\end{equation}

The hypotheses of Amplification \ref{amp} are satisfied for $C'$ and $k = \# (Y \cap C')$, with $M$ as in the statement of this theorem, and $l'\in\mathbf{H}_M(I)$: there exist $(m,\mhat,m')\in M$ satisfying (\ref{fundamental}).  Write $\frac{\mhat}{m}=\frac{ca}{2a-1}+\delta$ and $\frac{m'}{m^2}=\frac{a}{2a-1} + \eta$, assuming that $\delta$ and $\eta$ satisfy the hypotheses above.
\begin{eqnarray} 
e' + \frac{k}{2} 
	&<& \frac{(e'-g'+1)e + \frac{\mhat}{m}((e'-g'+1)d-(e-g+1)d')}{e-g+1} \notag\\
	  	&&{} + \frac{\frac{m'}{m^2}((e'-g'+1)n-(e-g+1)n')}{e-g+1} + \frac{bS}{m}  \notag\\
\Leftrightarrow \frac{k}{2}(e-g+1) 
	&<& e'(g-1) -e(g'-1) \notag\\
		&&{}+ {\textstyle\frac{a}{2a-1}}((e'-g'+1)(n+cd) - (e-g+1)(n'+cd')) \notag\\
		&&{} + \eta((e'-g'+1)n-(e-g+1)n') \notag\\
		&&{}+ \delta((e'-g'+1)d-(e-g+1)d') + \frac{bS}{m}(e-g+1). \hspace{.5in} \label{mess}
\end{eqnarray}
The final terms of the right hand side are already in the form we want them, so for brevity we shall only work on the first two lines. 
 \begin{eqnarray*}
\lefteqn{e'(g-1) -e(g'-1) + {\textstyle\frac{a}{2a-1}}((e'-g'+1)(n+cd) - (e-g+1)(n'+cd'))}\\
	&=& e'\left((g-1)+{\textstyle\frac{a}{2a-1}}(n+cd)\right) - (g'-1)(e + {\textstyle\frac{a}{2a-1}}(n+cd))  \hspace{1in}\\
		&&{}- {\textstyle\frac{a}{2a-1}}(e-g+1)(n'+cd')\\
	&=&\frac{e'}{2a-1}((2a-1)(g-1)+a(n+cd)) \\
		&&{}- (g'-1)(e-g+1 + {\textstyle\frac{1}{2a-1}}((2a-1)(g-1)+a(n+cd)) \\
		&&{}- {\textstyle\frac{a}{2a-1}}(e-g+1)(n'+cd') \\
	&=&\frac{e'}{2a-1}(e-g+1)-(g'-1)(e-g+1)(1+{\textstyle\frac{1}{2a-1}}) \\
		&&{}- {\textstyle\frac{a}{2a-1}}(e-g+1)(n'+cd'),
\end{eqnarray*}	
where for the last equality we have recalled that $(2a-1)(g-1)+a(n+cd)=e-g+1$.  We may substitute this back into line (\ref{mess}), multiplying through by $\frac{2a-1}{e-g+1}$ to obtain:
\begin{eqnarray*}
(2a-1)\frac{k}{2}
	&<& e' - a(2g'-2+n'+cd') + (2a-1)n\eta \left(\frac{e'-g'+1}{e-g+1} - \frac{n'}{n}\right) \\
		&& {}+(2a-1)d\delta\left(\frac{e'-g'+1}{e-g+1} - \frac{d'}{d}\right) + (2a-1)\frac{bS}{m}.
\end{eqnarray*}
Now use (\ref{ak-k}) to see
\begin{eqnarray}
\frac{k}{2} 
	&<& (2a-1)n\eta \left(\frac{e'-g'+1}{e-g+1} - \frac{n'}{n}\right)+(2a-1)d\delta \left(\frac{e'-g'+1}{e-g+1} - \frac{d'}{d}\right)\notag\\
		&&{} + (2a-1)\frac{bS}{m} \notag \\
\Rightarrow \frac{k}{2b}
	&<& (2a-1)n\eta \left(\frac{e'-g'+1}{e-g+1} - \frac{n'}{n}\right)+(2a-1)d\delta \left(\frac{e'-g'+1}{e-g+1} - \frac{d'}{d}\right) \notag\\
		&&{}+ (2a-1)\frac{S}{m} \label{I use this in the other paper}.
\end{eqnarray}
We must take care as $S$ varies with $k$; explicitly, $S=g+k(2g-\frac{3}{2})+q_2-\bar{g}+1$.  Thus the inequality we wish to contradict becomes
\begin{eqnarray}
\frac{k}{b}\left(\frac{1}{4a-2}-\frac{(2g-\frac{3}{2})}{m}\right)
	&<& n\eta \left(\frac{e'-g'+1}{e-g+1} - \frac{n'}{n}\right)
		+ d\delta \left(\frac{e'-g'+1}{e-g+1} - \frac{d'}{d}\right) \notag\\
	&& {}+ \frac{g+q_2-\bar{g}+1}{m} \label{now_estimate}.
\end{eqnarray}
It is time to use our bounds for $\eta$ and $\delta$.  Note that
\[\begin{array}{rcccl}
-1	&\leq&	\displaystyle\frac{e'-g'+1}{e-g+1} - \frac{n'}{n} 	&\leq&	1 \\
	&&&&\\
-1	&\leq&	\displaystyle\frac{e'-g'+1}{e-g+1} - \frac{d'}{d} 	&\leq&	1. 
\end{array}
\]
We assumed that $|n\eta|+|d\delta |\leq \frac{1}{4a-2}-\frac{3g+q_2-\bar{g}-\frac{1}{2}}{m}$.  It follows that
\[
n\eta \left(\frac{e'-g'+1}{e-g+1} - \frac{n'}{n}\right)+d\delta \left(\frac{e'-g'+1}{e-g+1} - \frac{d'}{d}\right)
	\leq \frac{1}{4a-2}-\frac{3g+q_2-\bar{g}-\frac{1}{2}}{m}.
\]
Hence line (\ref{now_estimate}) says
\[
\frac{k}{b}\left(\frac{1}{4a-2}-\frac{(2g-\frac{3}{2})}{m}\right)
	<\frac{1}{4a-2}-\frac{2g-\frac{3}{2}}{m}.
\]
By hypothesis $m>(6g+2q_2-2\bar{g}+1)(2a-1)>(2g-\frac{3}{2})(4a-2)$ and so we know that $\frac{1}{4a-2}-\frac{2g-\frac{3}{2}}{m}>0$.  Thus we have proved that $\frac{k}{b}<1$, a contradiction.

The contradiction implies that we cannot decompose $C$ into two strictly smaller sub-curves $C'$ and $Y$ as described.  Thus all the coefficients $a_{i}$ must be zero, and we have an isomorphism 
\begin{eqnarray*} 
\lefteqn{\left(\mathcal{O}_{\Pro(W)}(1)\otimes \mathcal{O}_{\Pro^{r}}(1)\right)|_{\mathcal{D}} }\\
	&\hspace{.5in}&\cong \acanon_{\mathcal{D}/\Spec R}(a\sigma_1(\Spec R)+\cdots +a\sigma_n(\Spec R))
			\otimes \mathcal{O}_{\Pro^{r}}(ca+1)|_{\mathcal{D}}.
\end{eqnarray*}
In particular, $(\mathcal{D}_{0},\sigma_1(0),\ldots ,\sigma_n(0))$ satisfies condition (iii) of Definition \ref{j}.  We conclude as described that it is represented in $J$, and so $\alpha(0)\in J\cap\bar{J}^{ss}$.  Hence $J\cap\bar{J}^{ss}(l)$ is closed in $\bar{J}^{ss}(l)$, which completes the proof.   
\end{proof}

{\it Remark:} A slightly larger range of values for $\frac{m'}{m^2}$ and $\frac{\mhat}{m}$ is possible; note that in fact $\frac{e'-g'+1}{e-g+1} - \frac{n'}{n}> -1$, enabling us to drop our lower bound to below $\frac{1}{4a-1}$.  It is not clear whether the upper bound can be improved.\\

Let us review what we know, given this result.  It is time to apply the theory of variation of GIT, to show that the semistable set $\bar{J}^{ss}(l)$ is the same for all $l\in\mathbf{H}_M(\bar{J})$, where $M$ is as in the statement of Theorem \ref{jssclosed}.  Recall the definitions from Section \ref{variation_of_git}.  In particular, we will make use of Proposition \ref{needed_variation}.

\begin{corollary} \label{nearly_done}
Let $M$ be as given in the statement of Theorem \ref{jssclosed}.  Let $l\in\mathbf{H}_M(I)$.  
\begin{enumerate}
\item[(i)] If $l\in\mathbf{H}_M(\bar{J})$ then $\bar{J}^{ss}(l)=\bar{J}^{s}(l)\subset J$.
\item[(ii)] If $(m,\mhat,m')\in M$ and if $\bar{J}^{ss}(L_{m,\mhat,m'})\neq\emptyset$ then, when we work over $\C$, 
\[
\bar{J}\dblq_{L_{m,\mhat,m'}} SL(W)\cong \Maps.
\]
\item[(iii)] The semistable set $\bar{J}^{ss}(l)$ is the same for all $l\in\mathbf{H}_M(\bar{J})$.
\end{enumerate}
\end{corollary}
\begin{proof}
Parts (i) and (ii) follow from Proposition \ref{Jss=Js}, Theorem \ref{wishful_thinking} and Theorem \ref{jssclosed}.  

Part (iii): The region $\mathbf{H}_{M}(\bar{J})$ is by definition convex, and lies in the ample cone $\mathbf{A}^G(X)$.  Part (i) has shown us that if $l\in\mathbf{H}_M(\bar{J})$ then $\bar{J}^{ss}(l)=\bar{J}^s(l)$.  Now the result follows from Proposition \ref{needed_variation}.
\end{proof}
We have completed the first part of the proof.  By Corollary \ref{nearly_done}, it only remains to show that $\bar{J}^{ss}(l)\neq\emptyset$ for at least one $l\in\mathbf{H}_M(\bar{J})$.


\section{The Construction Finished}\label{state_thm}\label{construction_finished}


\subsection{Statement of theorems}
We are now in a position to state the main theorem of this paper: for a specified range $M$ of values $(m,\mhat,m')$, the GIT quotient $\bar{J}\dblq_{L_{m,\mhat,m'}} SL(W)$ is isomorphic to $\Maps$.  
First we recall the notation from Section \ref{defns}.  The vector space $W$ is of dimension $e-g+1$, where $e=a(2g-2+n+cd)$, the integer $c$ being sufficiently large that this is positive.  We embed the domains of stable maps into $\Pro(W)$.  We denote by $I$ the Hilbert scheme of $n$-pointed curves in $\pwpr$ of bidegree $(e,d)$.  The subspace $J\subset I$ corresponds to $a$-canonically embedded curves, such that the projection to $\Pro^r$ is a moduli stable map; this is laid out precisely in Definition \ref{j}.

The constants $m_1$, $m_2$, $m_3$, $q_1$, $q_2$, $q_3$, $\mu_1$ and $\mu_2$ are all defined in Section \ref{constants}.  In particular we recall $m_3$ and $q_2$: if $m,\mhat\geq m_3$ then the morphism from $I$ to projective space, defined by $(h,x_1,\ldots ,x_n)\mapsto \hat{H}_{m,\mhat,m'}(h,x_1,\ldots ,x_n)$, is a closed immersion.  The constant $q_2$ is chosen so that $h^0(C,\mathcal{I}_C)\leq q_{2}$, for any curve $C\subset \pwpr$.  We also defined $\bar{g}$:
\begin{eqnarray*} 
\bar{g}:= \mbox{min} \{ 0, g_{\bar{Y}} 
	&| & \mbox{$\bar{Y}$ is the normalisation of a complete subcurve $Y$ contained} \\ 
	&& \mbox{in a connected fibre $\mathcal{C}_{h}$ for some $h \in \Hilb(\pwpr)$} \}.
\end{eqnarray*}
Recall that $\bar{g}$ is bounded below by $-(e+d)+1$.

In the statement of the theorem, note that the conditions on $m$ explicitly ensure that $\frac{1}{4a-2}-\frac{3g+q_2-\bar{g}-\frac{1}{2}}{m}>0$, and hence that the condition on $\eta$ and $\delta$ is satisfied when $\eta$ and $\delta$ are sufficiently small.  Further remarks on the bounds for $m,\mhat,m'$ and $a$ are given after the statement of Theorem \ref{jssclosed}, which concerns the same range of linearisations.
\begin{theorem}\label{maps}
Fix integers $g$, $n$ and $d\geq 0$ such there exist smooth stable $n$-pointed maps of genus $g$ and degree $d$.  Suppose $m,\mhat > m_{3}$ and
\[ m > \max \left \{ \begin{array}{c} 
	( g - \frac{1}{2} + e(q_{1}+1) + q_{3} +\mu_{1}m_2)(e-g+1), \\ 
	(10g+3q_{2}-3\bar{g})(e-g+1) \\
	(6g+2q_2-2\bar{g}-1)(2a-1)
\end{array} \right \} 
\]
with 
\begin{eqnarray*}
\frac{\mhat}{m} &=&\frac{ca}{2a-1} + \delta, \\
\frac{m'}{m^2}	&=&\frac{a}{2a-1} + \eta,
\end{eqnarray*}
where
\[
| n\eta |+| d\delta | \leq \frac{1}{4a-2}-\frac{3g+q_2-\bar{g}-\frac{1}{2}}{m},
\]
and in addition let $a$ be sufficiently large that
\[
d\frac{\mhat}{m}+n\frac{m'}{m^2} < \frac{1}{8}e - \frac{9}{8}g + \frac{7}{8}.
\]
Then, over $\C$,
\[
\bar{J}\dblq_{L_{m,\mhat,m'}} SL(W)=\Maps.
\]
\end{theorem}
\begin{corollary}[cf.\ \cite{FP} Lemma 8]\label{general_x_statement}
Let $m,\mhat,m'$ satisfy the conditions of Theorem \ref{maps}.  Let $X\stackrel{\iota}{\hookrightarrow}\Pro^r$ be a projective variety.  Let $\beta\in H_2(X)^+$ be the homology class of some stable map.  If $\beta=0$ suppose that $2g-2+n\geq 1$.  Write $d:=\iota_*(\beta)\in H_2(\Pro^r)^+$. Then there exists a closed subscheme $J_{X,\beta}$ of $J$, such that, over $\C$,
\[
\bar{J}_{X,\beta} \dblq_{L_{m,\mhat,m'}|_{\bar{J}_{X,\beta}}} SL(W) \cong \overline{\mathcal{M}}_{g,n}(X,\beta),
\]
where $\bar{J}_{X,\beta}$ is the closure of $J_{X,\beta}$ in $\bar{J}$.
\end{corollary}
Theorem \ref{maps} and Corollary \ref{general_x_statement} will have been proved when we know that $\bar{J}^{ss}(l)=\bar{J}^s(l)=J$ for $l\in\mathbf{H}_M(\bar{J})$\@.  By Corollary \ref{nearly_done}, it only remains to show non-emptiness of $\bar{J}^{ss}(l)$ for one such $l$\@.    This non-emptiness is proved by induction on $n$, the number of marked points.    The base case, $n=0$, closely follows Gieseker's method, \cite{G} Theorem 1.0.0, and was given in Swinarski's thesis \cite{Swin}.

The inductive step is different.  We take a stable map $f:(C_0,x_1,\ldots ,x_n)\rightarrow \Pro^r$, and remove one of the marked points, attaching a new genus $1$ component to $C_0$ in its place.  We extend $f$ over the new curve by defining it to contract the new component to a point.  The result is a Deligne-Mumford stable map of genus $g+1$, with $n-1$ marked points.  This is inductively known to have a GIT semi-stable model, and semi-stability of a model for $f:(C_0,x_1,\ldots ,x_n)\rightarrow \Pro^r$ follows.

In fact, the inductive step shows directly that $\bar{J}^{ss}(l)=J$ for all $l\in\mathbf{H}_M(\bar{J})$.  It follows from Proposition \ref{jssclosed} that $\bar{J}^{ss}(l)=\bar{J}^s(l)=J$ for such $l$.  Hence we know that $\bar{J}\dblq_{L_{m,\mhat,m'}} SL(W)\cong \Maps$ for $(m,\mhat,m')\in M$, without needing to appeal to independent constructions.  Of course such constructions are still needed for the base case.  However, when $r=d=0$, the base case is $\Mg$, constructed by Gieseker over $\Spec \Z$.  All the theory we are using is valid over any field, as we have extended the results we need from variation of GIT.  Thus $\curves$ is  constructed over $\Spec k$ for any field $k$.  As we shall show, this is sufficient to show that $\curves$ is in fact constructed over $\Spec\Z$.

As the constant $\mhat$ is irrelevant in the case $r=d=0$, we set it to zero and suppress it in the notation $L_{m,\mhat,m'}$.
\begin{theorem}\label{curves}
Let $g$ and $n\geq0$ be such that $2g-2+n>0$.  Set $e=a(2g-2+n)$. Suppose $m> m_{3}$ and
\[ m > \max \left \{ \begin{array}{c} 
	( g - \frac{1}{2} + e(q_{1}+1) + q_{3} +\mu_{1}m_2)(e-g+1), \\ 
	(10g+3q_{2}-3\bar{g})(e-g+1) \\
	(6g+2q_2-2\bar{g}-1)(2a-1)
\end{array} \right \} 
\]
with 
\begin{eqnarray*}
\frac{m'}{m^2}	&=&\frac{a}{2a-1} + \eta
\end{eqnarray*}
where
\[
| n\eta | \leq \frac{1}{4a-2}-\frac{3g+q_2-\bar{g}-1}{m},
\]
and in addition, ensure $a$ is sufficiently large that
\[
n\frac{m'}{m^2} < \frac{1}{8}e - \frac{9}{8}g + \frac{7}{8}.
\]
Then, as schemes over $\Spec \Z$,
\[
\bar{J}\dblq_{L_{m,m'}} SL(W)=\curves.
\]
\end{theorem}


\subsection{The base case: no marked points}\label{base_subsection}	
	
Before we can state the theorem that maps from smooth domain curves are semistable, there is a little more notation to mention.  The constant $\epsilon$ is found by Gieseker in the following lemma, which is based on \cite{morrison} Theorem 4.1.  Note that the hypothesis printed in \cite{G} is that $e \geq 20(g-1)$, but careful examination of the proof and \cite{morrison} shows that $e \geq 2g+1$ suffices.  
\begin{lemma}[\cite{G} Lemma 0.2.4] \label{morcomb} Fix two integers $g \geq 2$, $e \geq 2g+1$ and write $N = e-g$.  Then there exists $\epsilon > 0$ such that for all integers $r_{0}\leq\cdots\leq r_{N}$ (not all zero) with $\sum r_{i} = 0$ and for all integers $0 = e_{0}\leq\cdots\leq e_{N}=e$ satisfying:
\begin{enumerate}
\item[(i)] if $e_{j} > 2g-2$ then $e_{j} \geq j+g$;
\item[(ii)] if $e_{j} \leq 2g-2$ then $e_{j} \geq 2j$;  
\end{enumerate}
there exists a sequence of integers $0 = i_{1}\leq\cdots\leq i_{k}=N$ verifying the following inequality:
\[ \sum_{t=1}^{k-1} (r_{i_{t+1}}-r_{i_{t}}) (e_{i_{t+1}}+e_{i_{t}}) > 2r_{N}e + 2\epsilon(r_{n}-r_{0}). 
\] 
\hfill$\Box$
\end{lemma}
Now the statement of the theorem is:


\begin{theorem}[cf.\ \cite{G} 1.0.0] \label{1.0.0} For all $K>0$ there exist integers $p,b$ satisfying $m=(p+1)b>K$, such that for any $\mhat_1>2g-1$ satisfying $\mhat:=b\mhat_1>m_3$, if $C \subset \pwpr \rightarrow \Pro^r$ is a stable map, if $C$ is nonsingular, if the map 
\[ H^{0}(\Pro(W), \mathcal{O}_{\Pro(W)}(1)) \stackrel{\rho}{\rightarrow} H^{0}(p_{W}(C), \mathcal{O}_{p_{W}(C)}(1)) \]
is an isomorphism, and if $L_{W}$ is very ample (so that $C \cong p_{W}(C)$), 
then $C\in I^{ss}(L_{m,\mhat})$.
\end{theorem}
{\it Remark.} The values that $m$ and $\mhat$ must take will be made clear in the course of the proof.
\begin{proof}Let $C \subset \pwpr {\rightarrow} \Pro^{r}$ be such a map.  Let $\lambda$ be a 1-PS of $SL(W)$.  There exist a basis $\{ w_{0},\ldots,w_{N} \}$ of $H^{0}(\Pro(W), \mathcal{O}_{\Pro(W)}(1))$ and integers $r_{0} \leq \cdots \leq r_{N}$ such that $\sum r_{i} = 0$ and the action of $\lambda$ is given by $\lambda (t)w_{i} = t^{r_{i}}w_{i}$.  By our hypotheses the map $p_{W*}\rho: H^{0}(\Pro(W), \mathcal{O}_{\Pro(W)}(1)) \rightarrow H^{0}(C,L_{W})$ is injective. 
Write $w_{i}' := p_{W *} \rho (w_{i})$.  Let $E_{j}$ be the invertible subsheaf of $L_{W}$ generated by $w_{0}',\ldots,w_{j}'$ for $0 \leq j \leq N = e-g$, and write $e_{j} = \deg E_{j}$.  Note that $E_{N} = L_{W}$ since $L_{W}$ is very ample hence generated by global sections, $h^{0}(C,L_{W})=e-g+1$, and  $w_{0}',\ldots,w_{N}'$ are linearly independent.  The integers $e_{0},\ldots,e_{N}=e$ satisfy the following two properties:
\begin{enumerate}
\item[(i)] If $e_{j} > 2g-2$ then $e_{j} \geq j+g$.  
\item[(ii)] If $e_{j} \leq 2g-2$ then $e_{j} \geq 2j$.  
\end{enumerate}

To see this, note that since by definition $E_{j}$ is generated by $j+1$ linearly independent sections we have $h^{0}(C,E_{j}) \geq j+1$.  If $e_{j} = \deg E_{j} > 2g-2$ then $H^{1}(C,E_{j})=0$ so by Riemann-Roch $e_{j} = h^{0}-h^{1}+g-1 \geq j+g$.  If $e_{j} \leq 2g-2$ then $H^{0}(C,\omega_{C} \otimes E_{j}^{-1}) \neq 0$ so by Clifford's theorem $j+1 \leq h^{0} \leq \frac{e_{j}}{2}+1$.    

The hypotheses of Lemma \ref{morcomb} are satisfied with these $r_{i}$ and $e_{j}$, so there exist integers $0 = i_{1},\ldots,i_{k}=N$ such that   
\[ \sum_{t=1}^{k-1} (r_{i_{t+1}}-r_{i_{t}}) (e_{i_{t+1}}+e_{i_{t}}) > 2r_{N}e + 2\epsilon(r_{N}-r_{0}). 
\] 

Suppose $p$ and $b$ are positive integers, and set $m=(p+1)b$; assume $m>m_3$.   Recall that $H^{0}(\Pro(W), \mathcal{O}_{\Pro(W)}((p+1)b))$ has a basis consisting of monomials of degree $(p+1)b$ in $w_{0},\ldots,w_{N}$.  For all $1 \leq t \leq k$, let
\[
V_{i_{t}} \subset H^{0}(\Pro(W), \mathcal{O}_{\Pro(W)}(1)) 
\]
be the subspace spanned by $\{ w_{0},\ldots,w_{i_{t}} \}$.  Let $\mhat_1$ be another positive integer, such that $\mhat:=b\mhat_1>m_3$, so that
\[\phat^C_{(p+1)b, \mhat} : H^{0}(\pwpr,\mathcal{O}_{\Pro(W)}((p+1)b) \otimes \mathcal{O}_{\Pro^{r}}(\mhat)) \rightarrow H^{0}(C,L_{W}^{(p+1)b} \otimes L_{r}^{\mhat})
\] 
is surjective.  For all triples $(t_{1}, t_{2}, s)$ with $1 \leq t_{1} < t_{2} \leq k$ and $0 \leq s \leq p$ let 
\[
V_{i_{t_{1}}}^{p-s}V_{i_{t_{2}}}^{s}V_{N} \subset H^{0}(\pwpr, \opwopr{p+1}{\mhat_1})
\]
be the subspace spanned by monomials of bidegree $(p+1,\mhat_1)$, where the degree $p+1$ part has the following form:
\[
x_{j_{1}} \cdots x_{j_{p-s}} y_{j_{1}} \cdots y_{j_{s}} z_{j} 
	\mbox{ where $x_{j_k}\in V_{i_{t_1}}$, $y_{j_k}\in V_{i_{t_2}}$, $z_j\in V_N=W$}.
\]
We write $(V_{i_{t_{1}}}^{p-s}V_{i_{t_{2}}}^{s}V_{N})^b$ for  $\mbox{Sym}^b(V_{i_{t_{1}}}^{p-s}V_{i_{t_{2}}}^{s}V_{N})$, which consists of monomials of bidegree $((p+1)b,b\mhat_1)$ in  $H^{0}(\pwpr, \opwopr{(p+1)b}{\mhat})$.  Now set
\[
(\overline{V_{i_{t_{1}}}^{p-s}V_{i_{t_{2}}}^{s}V_{N}})^b 
	:= \phat^C_{(p+1)b,\mhat}((V_{i_{t_{1}}}^{p-s}V_{i_{t_{2}}}^{s}V_{N})^b).
\]

We have obtained a filtration of $H^{0}(C,L_W^{(p+1)b}\otimes L_r^{\mhat})$:
\begin{equation} \label{big filtration}
\begin{array}{lclclcccl}
0 & \subset & (\overline{V_{i_{1}}^{p}V_{i_{2}}^{0}V_{N}})^{b} & \subset& (\overline{V_{i_{1}}^{p-1}V_{i_{2}}^{1}V_{N}})^{b}&  \subset&  \cdots & \subset & (\overline{V_{i_{1}}^{1}V_{i_{2}}^{p-1}V_{N}})^{b} \\
& \subset & (\overline{V_{i_{2}}^{p}V_{i_{3}}^{0}V_{N}})^{b} & \subset& (\overline{V_{i_{2}}^{p-1}V_{i_{3}}^{1}V_{N}})^{b} &  \subset &  \cdots & \subset &  (\overline{V_{i_{2}}^{1}V_{i_{3}}^{p-1}V_{N}})^{b}\\
& \vdots & \vdots & \vdots & \vdots &  \vdots & \vdots & \vdots & \vdots \\
& \subset & (\overline{V_{i_{t}}^{p}V_{i_{t+1}}^{0}V_{N}})^{b} & \subset& (\overline{V_{i_{t}}^{p-1}V_{i_{t+1}}^{1}V_{N}})^{b}&  \subset&  \cdots & \subset& (\overline{V_{i_{t}}^{1}V_{i_{t+1}}^{p-1}V_{N}})^{b} \\
& \vdots & \vdots & \vdots & \vdots &  \vdots & \vdots & \vdots & \vdots \\
& \subset & (\overline{V_{i_{k-1}}^{p}V_{i_{k}}^{0}V_{N}})^{b} & \subset& (\overline{V_{i_{k-1}}^{p-1}V_{i_{k}}^{1}V_{N}})^{b} &  \subset&  \cdots & \subset & (\overline{V_{i_{k-1}}^{1}V_{i_{k}}^{p-1}V_{N}})^{b} \\
& \subset& (\overline{V_{i_{k-1}}^{0}V_{i_{k}}^{p}V_{N}})^{b}& =& \lefteqn {H^{0}(C,L_W^{(p+1)b}\otimes L_r^{\mhat}).}
\end{array}
\end{equation}


\begin{claim}[\cite{G} page 29]\label{smooth are stable claim} There exists an integer $b'$ which is independent of $C$, $t_{1}$, and $t_{2}$ such that if $b \geq b'$ then 
\begin{displaymath}
(\overline{V_{i_{t_{1}}}^{p-s}V_{i_{t_{2}}}^{s}V_{N}})^{b} = H^{0}(C, (E_{i_{t_{1}}}^{p-s} \otimes E_{i_{t_{2}}}^{s} \otimes L_{W})^{b} \otimes L_{r}^{\mhat}).
\end{displaymath}
\end{claim}

{\it Proof of Claim \ref{smooth are stable claim}.}

By hypothesis $L_{W}$ is very ample.  Note that $\deg L_{r}=d > 0$ since $C$ is nonsingular hence irreducible.  Thus if $\mhat_1 > 2g+1$ then $L_{r}^{\mhat_1}$ is very ample.  

By definition of the sheaves $E_{j}$, it follows that the linear system $(V_{i_{t_{1}}}^{p-s}V_{i_{t_{2}}}^{s}V_{N})^{b}$ restricted to $C$ generates 
\[
 (E_{i_{t_{1}}}^{p-s} \otimes E_{i_{t_{2}}}^{s} \otimes L_{W})^{b} \otimes L_{r}^{\mhat}
\]
and
\begin{equation} \label{inject} 
\phat^C_{(p+1)b, \mhat}((V_{i_{t_{1}}}^{p-s}V_{i_{t_{2}}}^{s}V_{N})^{b})\subseteq  H^{0}(C, (E_{i_{t_{1}}}^{p-s} \otimes E_{i_{t_{2}}}^{s} \otimes L_{W})^{b} \otimes L_{r}^{\mhat}).
\end{equation}

Recall that $\mhat = b\mhat_{1}$ and suppose $\mhat_{1} > 2g+1$ so that $L_{r}^{\mhat_{1}}$ is very ample hence generated by global sections.  $E_{i_{t_{1}}}$ and $ E_{i_{t_{2}}}$ are generated by global sections, so it follows that $\phat^C_{p+1,\mhat_{1}}(V_{i_{t_{1}}}^{p-s}V_{i_{t_{2}}}^{s}V_{N})$ is a very ample base point free linear system on $C$.

Let $\psi = \psi_{p+1,\mhat_{1}}$ be the projective embedding corresponding to the linear system $\phat^C_{p+1,\mhat_{1}}(V_{i_{t_{1}}}^{p-s}V_{i_{t_{2}}}^{s}V_{N})$.  Let $\mathcal{I}_{C/\Pro}$ be the ideal sheaf defining $C$ as a closed subscheme of $\Pro := \Pro(\phat^C_{p+1,\mhat_{1}}(V_{i_{t_{1}}}^{p-s}V_{i_{t_{2}}}^{s}V_{N}))$.  There is an exact sequence of sheaves on $\Pro$ as follows:
\[0 \rightarrow \mathcal{I}_{C/\Pro} \rightarrow \mathcal{O}_{\Pro} \rightarrow \psi_{*} \mathcal{O}_{C} \rightarrow 0.
\]  
Tensoring by the very ample sheaf $\mathcal{O}_{\Pro}(b)$ we obtain
\begin{equation} \label{psistar} 0 \rightarrow \mathcal{I}_{C/\Pro}(b) \rightarrow \mathcal{O}_{\Pro}(b) \rightarrow (\psi_{*} \mathcal{O}_{C})(b) \rightarrow 0.
\end{equation}

Write
\[ \mathcal{F} := E_{i_{t_{1}}}^{p-s} \otimes E_{i_{t_{2}}}^{s} \otimes L_{W} \otimes L_{r}^{\mhat_{1}} \cong \psi^{*} \mathcal{O}_{\Pro}(1).
\]

We have
\begin{eqnarray*} (\psi_{*} \mathcal{O}_{C})(b) 
	& := & \psi_{*} \mathcal{O}_{C} \otimes_{\mathcal{O}_{\Pro}} \mathcal{O}_{\Pro}(b) \\
	& \cong & \psi_{*}(\mathcal{O}_{C} \otimes_{\mathcal{O}_{C}} \psi^{*} \mathcal{O}_{\Pro}(1)^{b})\\
	& \cong & \psi_{*}(\mathcal{F}^{b}) \mbox{ since $\psi^{*} \mathcal{O}_{\Pro}(1) \cong \mathcal{F}$}.
\end{eqnarray*}
Now the exact sequence (\ref{psistar}) reads
\[ 
0 \rightarrow \mathcal{I}_{C/\Pro}(b) \rightarrow \mathcal{O}_{\Pro}(b) \rightarrow \psi_{*} (\mathcal{F}^{b}) \rightarrow 0.
\]
In the corresponding long exact sequence in cohomology we have 
\begin{equation} \label{h1exactseq}
\cdots \rightarrow H^{0}(\Pro, \mathcal{O}_{\Pro}(b)) \rightarrow  H^{0}(\Pro, \psi_{*}\mathcal{F}^{b} ) \rightarrow  H^{1}(\Pro, \mathcal{I}_{C/\Pro}(b)) \rightarrow \cdots.
\end{equation}

The so-called ``Uniform $m$ Lemma'' (cf. \cite{HM} Lemma 1.11) ensures that there is an integer $b' > 0$ depending on the Hilbert polynomial $P$ but not on the curve $C$ such that $ H^{1}(\Pro, \mathcal{I}_{C/\Pro}(b))=0$ if $b > b'$.  Then for such $b$ the exact sequence (\ref{h1exactseq}) implies that the map  
\[  H^{0}( \Pro,  \mathcal{O}_{\Pro}(b)) \rightarrow  H^{0}(C, \mathcal{F}^{b})
\]
is surjective.  Recall that $\Pro := \Pro(\phat^C_{p+1,\mhat_{1}}(V_{i_{t_{1}}}^{p-s}V_{i_{t_{2}}}^{s}V_{N}))$.  Then  
\[  
H^{0}(\Pro,\mathcal{O}_{\Pro}(b)) \cong \mbox{Sym}^{b}(\phat^C_{p+1,\mhat_{1}}(V_{i_{t_{1}}}^{p-s}V_{i_{t_{2}}}^{s}V_{N})).
\] 
Also there is a surjection 
\[ 
\mbox{Sym}^{b}(V_{i_{t_{1}}}^{p-s}V_{i_{t_{2}}}^{s}V_{N}) 
	\rightarrow \mbox{Sym}^{b}(\phat^C_{p+1,\mhat_{1}}(V_{i_{t_{1}}}^{p-s}V_{i_{t_{2}}}^{s}V_{N})),
\]
so putting this all together we have a surjection
\[ \label{12}
\mbox{Sym}^{b}(V_{i_{t_{1}}}^{p-s}V_{i_{t_{2}}}^{s}V_{N} ) \rightarrow  H^{0}(C,\mathcal{F}^{b}),
\]
i.e.,
\begin{equation}\label{surject}
\phat^C_{(p+1)b, \mhat}: (V_{i_{t_{1}}}^{p-s}V_{i_{t_{2}}}^{s}V_{N})^{b}
	\rightarrow H^{0}(C, (E_{i_{t_{1}}}^{p-s} \otimes E_{i_{t_{2}}}^{s} \otimes L_{W})^{b} \otimes L_{r}^{\mhat})
\end{equation}
is surjective.
It follows from lines (\ref{inject}) and (\ref{surject}) that 
\[
(\overline{V_{t_{1}}^{p-s}V_{t_{2}}^{s}V_{N}})^{b} = H^{0}(C, (E_{i_{t_{1}}}^{p-s} \otimes E_{i_{t_{2}}}^{s} \otimes L_{W})^{b} \otimes L_{r}^{\mhat}),
\]
completing the proof of Claim \ref{smooth are stable claim}.\hfill $\Box$


{\it Proof of Theorem \ref{1.0.0} continued.}  Take $b \geq 2g+1$ so that 
\\ $H^{1}(C, (E_{i_{t_{1}}}^{p-s} \otimes E_{i_{t_{2}}}^{s} \otimes L_{W})^{b} \otimes L_{r}^{\mhat}) =0.$  We use Riemann-Roch to calculate
\begin{eqnarray}
\dim {(\overline{V_{i_{t}}^{p-s}V_{i_{t+1}}^{s}V})^{b}} &  = & h^{0}(C, (E_{i_{t}}^{p-s} \otimes E_{i_{t+1}}^{s} \otimes L_{W})^{b} \otimes L_{r}^{\mhat}) \nonumber \\
& = & b((p-s)e_{i_{t}} + se_{i_{t+1}}+e) + d\mhat -g +1.
\end{eqnarray}

 We assume for the rest of the proof that $p$, $b$ and $\mhat_1$ are sufficiently large that \begin{eqnarray*} 
p 	& > & \max \{ e+g, \frac{\frac{3}{2}e+1}{\epsilon}\}, \\ 
b 	& > & \max \{ p, (2g+1)b' \}, \\ 
m 	& := &(p+1)b > \max\{m_3,K\} , \\
\mhat_1 & > & 2g+1, \\
\mhat	&:= & b\mhat_1 > m_3. \label{pandn}
\end{eqnarray*}

Choose a basis $\hat{B}_{(p+1)b,\mhat}$ of $H^{0}(\pwpr,\mathcal{O}_{\Pro(W)}((p+1)b) \otimes \mathcal{O}_{\Pro^{r}}(\mhat))$ of monomials $\hat{M}_{i}$ of bidegree $((p+1)b,\mhat)$.  Pick monomials $\hat{M}_1,\ldots,\hat{M}_{P((p+1)b,\mhat)}$ in $\hat{B}_{(p+1)b,\mhat}$ such that, $\phat^C_{(p+1)b, \mhat}(\hat{M}_{1}), \ldots, \phat^C_{(p+1)b, \mhat}(\hat{M}_{P((p+1)b,\mhat)})$ is a basis of $H^{0}(C,L_{W}^{(p+1)b} \otimes L_{r}^{\mhat})$ which respects the filtration (\ref{big filtration}).  Observe that if a monomial $\hat{M} \in {(\bar{V}_{i_{1}}^{p-s}\bar{V}_{i_{2}}^{s}\bar{V})^{b}} - {(\bar{V}_{i_{1}}^{p-s+1}\bar{V}_{i_{2}}^{s-1}\bar{V})^{b}}$ then $\hat{M}$ has $\lambda$-weight $w_{\lambda}(\hat{M}) \leq n((p-s)r_{i_{t}}+sr_{i_{t+1}}+ r_{N})$.  Moreover, as our basis respects the filtration (\ref{big filtration}), we may count how many such $\hat{M}$ there are.

We now estimate the total $\lambda$-weight of $\hat{M}_{1},\ldots,\hat{M}_{P((p+1)b,\mhat)}$, which gives an upper bound for $\mu^{L_{m,\mhat}}(C,\lambda)$, as follows:
\begin{eqnarray}
 \lefteqn{\mu^{L_{m,\mhat}}(C,\lambda) \leq \displaystyle \sum_{i=1}^{P(m)+d\mhat} w_{\lambda}(\hat{M}_{i}) 
 \leq   b(pr_{i_{1}}+r_{N}) \dim {(\overline{V_{i_{1}}^{p}V_{i_{2}}^{0}V})^{b}}} \nonumber \\
 	& {}+ \displaystyle \sum_{\stackrel{\mbox{\scriptsize $0 \leq s \leq p \hspace{14pt}$}}{{\mbox{\scriptsize $1 \leq t \leq k-1$}}}}&
  	b((p  - s)r_{i_{t}}+sr_{i_{t+1}}+ r_{N}) \left(\dim {(\overline{V_{i_{t}}^{p-s}V_{i_{t+1}}^{s}V})^{b}}\right. \notag\\
	&&	{}\left.- \dim{(\overline{V_{i_{t}}^{p-s+1}V_{i_{t}}^{s-1}V})^{b}} \right). \label{longsum}
\end{eqnarray}
The first term on the right hand side of (\ref{longsum}) is
\begin{eqnarray*} 
\lefteqn{ b(pr_{i_{1}}+r_{N})(b(pe_{i_{1}}+e_{i_k})+d\mhat-g+1)} \\
 	&=& b(pr_0+r_{N})d\mhat + b(pr_0+r_{N})(b(pe_0+e)-g+1).
\end{eqnarray*}
The factor $\dim {(\overline{V_{i_{t}}^{p-s}V_{i_{t+1}}^{s}V})^{b}} - \dim {(\overline{V_{i_{t}}^{p-s+1}V_{i_{t}}^{s-1}V})^{b}}$ of the summand is 
\begin{displaymath}
\begin{array}{l} 
\left( b((p-s)e_{i_{t}} + se_{i_{t+1}}+e) + d\mhat -g +1 \right)  \\ 
	\hspace{.25in} \mbox{} - \left( b((p-s+1)e_{i_{t}} + (s-1)e_{i_{t+1}}+e) + d\mhat -g +1 \right) 
	= b(e_{i_{t+1}}-e_{i_{t}}).
\end{array}
\end{displaymath}
Note that nearly all of the terms having $d\mhat$ as a factor have ``telescoped.''  We have 
\begin{eqnarray}
\mu^{L_{m,\mhat}}(C,\lambda) 
	& \leq & db(pr_{0}+r_{N})\mhat + b(pr_{0}+r_{N})(b(pe_{0}+e)-g+1) \notag \\
 	& & {} + \displaystyle  \sum_{\stackrel{\mbox{\scriptsize $0 \leq s \leq p \hspace{14pt}$}}{{\mbox{\scriptsize $1 \leq t \leq k-1$}}}}  b((p  -  s)r_{i_{t}}+sr_{i_{t+1}}+ r_{N}) \left(b(e_{i_{t+1}}-e_{i_{t}}) \right).\hspace{.2in} \notag\label{expan}
\end{eqnarray}
The sum of the second two terms is {\it exactly} the expression Gieseker obtains at the bottom of page 30 in \cite{G}.  Following Gieseker's calculations up to page 34, we see:
\begin{eqnarray*}
\mu^{L_{m,\mhat}}(C,\lambda)
	&<&  db(pr_{i_{1}}+r_{N})\mhat + b^{2}p(r_{N}-r_{0}) \left( - \epsilon p + \frac{3e}{2} + \frac{e+g-1}{p} \right) \\
	&<& db(pr_{i_{1}}+r_{N})\mhat.
\end{eqnarray*}
where the  last inequality follows because $p > \max \{ e+g, \frac{\frac{3}{2}e+1}{\epsilon} \}$.  Next, we may estimate $r_N=\sum_{i=0}^{N-1}-r_i\leq -Nr_0$, and we know that $r_0<0$, so we have shown that
\[
\mu^{L_{m,\mhat}}(C,\lambda)
	\leq dbr_0(p-N)d\mhat <0,
\]
as, by hypothesis, $p>e-g=N$.

Nowhere in the proof have we placed any conditions on the 1-PS $\lambda$, so the result is true for every 1-PS of $SL(W)$.  Thus $C$ is $SL(W)$-stable with respect to $L_{m,\mhat}$.
\end{proof}
Note in particular that the hypotheses are satisfied by all smooth maps represented in $J$, by definition of $\mathcal{L}$ and $a$ (cf.\ Section \ref{defns}).  We may state the base cases of our induction:
\begin{proposition}\label{base_case_maps}
Fix $n=0$. Let $M$ consist of those $(m,\mhat,m')$ such that $m'\geq 1$, and $m,\mhat > m_{3}$ with
\[m > \max \left \{ \begin{array}{c} 
	( g - \frac{1}{2} + e(q_{1}+1) + q_{3} +\mu_{1}m_2)(e-g+1), \\ 
	(10g+3q_{2}-3\bar{g})(e-g+1) \\
	(6g+2q_2-2\bar{g}-1)(2a-1)
\end{array} \right \} 
\]
while
\begin{eqnarray*}
\frac{\mhat}{m} &=&\frac{ca}{2a-1} + \delta \\
\end{eqnarray*}
where
\[
| d\delta | \leq \frac{1}{4a-2}-\frac{3g+q_2-\bar{g}-\frac{1}{2}}{m},
\]
and in addition $a$ is sufficiently large that
\[
d\frac{\mhat}{m} < \frac{1}{8}e - \frac{9}{8}g + \frac{7}{8}.
\]
Let $l\in\mathbf{H}_M(\bar{J})$.  Then, as schemes over $\C$,
\[
\bar{J}^{ss}(l)=\bar{J}^{s}(l)=J.
\]
\end{proposition}
\begin{proof}
By Theorem \ref{1.0.0} there exists $(m,\mhat,m')\in M$ such that $I^{ss}(L_{m,\mhat,m'})$ is non-empty.  In particular, any smooth curve in $J$ satisfies the hypotheses of Theorem \ref{1.0.0} and so $\bar{J}^{ss}(L_{m,\mhat,m'})$ is non-empty.  By Theorem \ref{jssclosed} we know that $\bar{J}^{ss}(L_{m,\mhat,m'})\subset J$.  Then by Theorem \ref{wishful_thinking}(ii) we see that $\bar{J}^{ss}(L_{m,\mhat,m'})=\bar{J}^s(L_{m,\mhat,m'})=J$.  The result follows for all $l\in\mathbf{H}_M(\bar{J})$ by Corollary \ref{nearly_done}.
\end{proof}
The base case for stable curves is proved over $\Spec \Z$:
\begin{proposition}\label{base_case_curves}
Fix $n=0$, and $r=d=0$.  Let $M$ consist of those $(m,m')\in\N^2$ such that $m'\geq 1$, and $m > m_{3}$ with
\[m > \max \left \{ \begin{array}{c} 
	( g - \frac{1}{2} + e(q_{1}+1) + q_{3} +\mu_{1}m_2)(e-g+1), \\ 
	(10g+3q_{2}-3\bar{g})(e-g+1) 
\end{array} \right \}. 
\]
Let $l\in\mathbf{H}_M(\bar{J})$.  Then, as schemes over $\Spec \Z$,
\[
\bar{J}^{ss}(l)=\bar{J}^{s}(l)=J.
\]
\end{proposition}
\begin{proof}
By Theorem \ref{1.0.0} there exists $(m,m')\in M$ be such that $I^{ss}(L_{m,\mhat})$ is non-empty.  In particular, any smooth curve in $J$ satisfies the hypotheses of Theorem \ref{1.0.0} and so $\bar{J}^{ss}(L_{m,\mhat})$ is non-empty.
However, by \cite{G} Theorem 2.0.2, the GIT quotient $\bar{J}\dblq_{L_{m,m'}}SL(W)$ is isomorphic to the moduli space of stable curves, $\overline{\mathcal{M}}_g$.  In particular Gieseker proves here that every stable curve is represented in $\bar{J}\dblq_{L_{m,m'}}SL(W)$.  Then, by Corollary \ref{nearly_done}(iii), we see $\bar{J}^{ss}(l)=J$ for every $l\in\mathbf{H}_M(\bar{J})$.
\end{proof}


\subsection{The general case}\label{ind_step}\label{end}	
We suppose $n>0$ and fix $M_{g,n,d}$ as in the statement of Theorem \ref{maps}, then use induction to show that $\bar{J}_{g,n,d}^{ss}(l)=\bar{J}_{g,n,d}^s(l)=J$ for all $l\in\mathbf{H}_M(\bar{J})$.  By Corollary \ref{catquodone} it follows that $\bar{J}\dblq_{L_{m,\mhat,m'}} SL(W)\cong\Maps$ for any $(m,\mhat,m')\in M$.  As our inductive work uses various spaces of maps with differing genera and numbers of marked points, we shall be more precise about using the subscripts $g,n,d$ in this section.  Note that indeed $m_1$, $m_2$, $m_3$, $q_1$, $q_2$, $q_3$, $\mu_1$ and $\mu_2$ all depend on the genus.  We may specify in addition that these are increasing as functions of the genus.

The inductive hypothesis is given in the following proposition.  We fix $g$, $n$ and $d$ such that $n\geq 1$ and smooth stable $n$-pointed maps of genus $g$ and degree $d$ do exist.  Fix an integer $c\geq 2$ (the case $(g,n,d)=(0,0,1)$ does not arise here so this will suffice according to the remark in Section \ref{bound_for_c}).  We suppose that the theorem has been proved stable maps of degree $d$, from curves of genus $g+1$ with $n-1$ marked points.  As the base case, where the map has no marked points, works for any genus this is a valid inductive hypothesis to make.  Since our assumptions imply $2g-2+n+cd>0$ it follows that $2(g+1)-2+(n-1)+cd>0$, so smooth stable maps of this type do exist.

\begin{proposition}\label{inductive_step}
Let $g$, $n$ and $d$ be such that $2g-2+n+cd>0$ and assume in addition that $n\geq 1$.  Let $M_{g,n,d}$ consist of those $(m,\mhat,m')$ such that $m,\mhat > m_{3}$ and
\begin{equation} m > \max \left \{ \begin{array}{c} 
	( g - \frac{1}{2} + e_{g,n,d}(q_{1}+1) + q_{3} +\mu_{1}m_2)(e_{g,n,d}-g+1), \\ 
	(10g+3q_{2}-3\bar{g})(e_{g,n,d}-g+1) \\
	(6g+2q_2-2\bar{g}-1)(2a-1)
\end{array} \right \} 
\label{bounds_for_m}
\end{equation}
where $m_1$, $m_2$, $m_3$, $q_1$, $q_2$, $q_3$, $\mu_1$ and $\mu_2$, are those defined in Section \ref{constants}, taken to be functions of the genus $g$; in addition
\begin{eqnarray*}
\frac{\mhat}{m} &=&\frac{ca}{2a-1} + \delta \\
\frac{m'}{m^2}	&=&\frac{a}{2a-1} + \eta
\end{eqnarray*}
where
\[
| n\eta |+ | d\delta | \leq \frac{1}{4a-2}-\frac{3g+q_2-\bar{g}-\frac{1}{2}}{m},
\]
and ensure in addition that $a$ is sufficiently large that
\[
d\frac{\mhat}{m}+n\frac{m'}{m^2} < \frac{1}{8}e - \frac{9}{8}g + \frac{7}{8}.
\]
Work over $\Spec k$.  Assume that
\[
\bar{J}_{g+1,n-1,d}^{ss}(l')=\bar{J}_{g+1,n-1,d}^{s}(l')=J_{g+1,n-1,d}
\]
for all $l'\in\mathbf{H}_{M_{g+1,n-1,d}}(\bar{J}_{g+1,n-1,d})$.  Then, for all $l\in \mathbf{H}_{M_{g,n,d}}(\bar{J}_{g,n,d})$,
\[
\bar{J}_{g,n,d}^{ss}(l)=\bar{J}_{g,n,d}^{s}(l)=J_{g,n,d}
\]
as schemes over an arbitrary field $k$.
\end{proposition}
\begin{proof}
Note that if, $n\geq 1$ then 
\begin{eqnarray*}
e_{g+1,n-1,d}-(g+1)+1
	&=&(2a-1)g+a(n-1+cd)\\
	&>&(2a-1)(g-1)+a(n+cd)=e-g+1.  
\end{eqnarray*}
We conclude, by the definition and specifications given above, that $M_{g+1,n-1,d}\subseteq M_{g,n,d}$. 
Fix specific integers $(m,\mhat,m')\in M_{g+1,n-1,d}$, satisfying
\begin{eqnarray*}
\frac{\mhat}{m} &=&\frac{ca}{2a-1} \\
\frac{m'}{m^2} &=& \frac{a}{2a-1}
\end{eqnarray*}
and also such that $\frac{m}{n}(1-S_{14})$ is an integer, where $S_{14}:=\frac{g(a-1)}{(2a-1)g+a(n-1+cd)}$.  Our inductive hypothesis implies in particular that 
\[
\bar{J}_{g+1,n-1,d}^{ss}(L_{m,\mhat,m'})=J_{g+1,n-1,d}.  
\]
We shall now find $m''$ such that $(m,\mhat,m'')\in M_{g,n,d}$, with $\bar{J}_{g,n,d}^{ss}(L_{m,\mhat,m''})=J_{g,n,d}$.

Fix some $(h,\pts)\in J_{g,n,d}$.  Write $C_0:=\mathcal{C}_h$, so that $(h,\pts)$ models a stable map $p_r:(C_0,\pts)\rightarrow \Pro^r$ in $\Maps$.  Also, fix an elliptic curve $({C}_1,y)\subset\Pro(W_{1,1,0})$ represented in $J_{1,1,0}$.

Let $ev:H^0(\Pro(W_{g,n,d}),\Osheaf_{\Pro(W_{g,n,d})}(1))\rightarrow k$ be the evaluation map at the closed point $p_{W_{g,n,d}}(x_n)\in\Pro(W_{g,n,d})$, and let $V_{g,n,d}$ be its kernel, so that $V_{g,n,d}$ is the codimension $1$ subspace of $W_{g,n,d}$ consisting of sections vanishing at $p_{W_{g,n,d}}(x_n)$.  Similarly let $V_{1,1,0}$ be the codimension one subspace of $W_{1,1,0}$ corresponding to sections vanishing at $y$.

Now note:
\begin{eqnarray*}
\dim V_{g,n,d} + \dim V_{1,1,0} + 1 
	&=& a(2g-2+n+cd)-g + a\cdot1 -1 + 1 \\
	&=& \dim W_{g+1,n-1,d}.
\end{eqnarray*}
Hence, if we let $U$ be a dimension $1$ vector space over $k$, we may pick an isomorphism
\[
W_{g+1,n-1,d}\cong V_{g,n,d}\oplus U \oplus V_{1,1,0}.
\]
We further choose isomorphisms $W_{g,n,d}\cong V_{g,n,d}\oplus U$\label{decomp_Wgn} and $W_{1,1,0}\cong U\oplus V_{1,1,0}$ which fix $V_{g,n,d}$ and $V_{1,1,0}$ respectively.  Thus we regard $W_{g,n,d}$ and $W_{1,1,0}$ as subspaces of $W_{g+1,n-1,d}$.  The most important of these identifications we shall write: 
\[
W_{g+1,n-1,d}=W_{g,n,d}\oplus V_{1,1,0}.
\]
We project $W_{g+1,n-1,d}\rightarrow W_{g,n,d}$ along $V_{1,1,0}$ and induce an embedding
\[
\Pro(W_{g,n,d})\hookrightarrow \Pro(W_{g+1,n-1,d}); 
\]
similarly $\Pro(W_{1,1,0})\hookrightarrow \Pro(W_{g+1,n-1,d})$.  Then we induce closed immersions $p_{W_{g,n,d}}(C_0)\hookrightarrow \Pro(W_{g+1,n-1,d})$ and $C_1\hookrightarrow \Pro(W_{g+1,n-1,d})$; we shall consider the curves as embedded in this space.  

If $s\in V_{1,1,0}\subset W_{g+1,n-1,d}$ is regarded as a section of $\Osheaf_{\Pro(W_{g+1,n-1,d})}(1)$, then $s(x)=0$ for any $x$ in $\Pro(W_{g,n,d})$, and in particular $s(x)=0$ for any $x$ in $p_W(C_0)$. \label{vanish}In other words, $\rho^{p_W(C_0)}(V_{1,1,0})=\{0\}$, where we write $\rho^{p_W(C_0)}$ for restriction of sections to $p_W(C_0)$.  Similarly $\rho^{C_1}(V_{g,n,d})=\{0\}$.\label{stuck_together}

The images of $\Pro(W_{g,n,d})$ and $\Pro(W_{1,1,0})$ meet only at one point, $\Pro(U)\in\Pro(W_{g+1,n-1,d})$.  We shall denote this point by $P$.  If the curves $C_0$ and $C_1$ meet, it could only be at this point.  Consider $p_W(x_n)$; by the definitions we know $s(p_W(x_n))=0$ for all $s\in V_{1,1,0}$ and for all $s\in V_{g,n,d}$.  We conclude that $p_W(x_n)\in \Pro(U)$.  Similarly $y\in\Pro(U)$.  Thus, after the curves have been embedded in $\Pro(W_{g+1,n-1,d})$, the points $p_W(x_n)$ and $y$ coincide at $P$.  We define:
\[
(C,x_1,\ldots ,x_{n-1}):=(C_0\cup C_1,x_1,\ldots ,x_{n-1}).  
\]
As the curves $C_0$ and $C_1$ are smooth at $x_n$ and $y$ respectively, and they lie in two linear subspaces meeting transversally at $P$, the singular point of $C$ at $P$ is a node. 

The map $p_r:(C_0,\pts )\rightarrow \Pro^r$ may be extended over $C$ if we define it to contract $C_1$ to the point $p_r(x_n)$.  Thus we have the graph of a pre-stable map, 
\[
(C,x_1,\ldots ,x_{n-1})\subset \Pro(W_{g+1,n-1,d})\times\Pro^r.
\]
We wish to show that $(C,x_1,\ldots ,x_{n-1})$ is represented by a point in $J_{g+1,n-1,d}$, so we check conditions (i)-(iii) of Definition \ref{j}.  Clearly (i) is satisfied: $C$ is projective, connected, reduced and nodal, and the $n-1$ marked points are distinct and non-singular.  By construction, $C\rightarrow \Pro(W_{g+1,n-1,d})$ is a non-degenerate embedding, so (ii) is satisfied.  We must check (iii), the isomorphism of line bundles.

In general, if $C$ is a nodal curve and $C'\subset C$ is a complete subcurve, meeting the rest of $C$ in only one node at $Q$, then
\[
\omega_C|_{C'}= \omega_{C'}(Q).
\]
Thus, in our situation, as $(C_0,\pts)$ is represented in $J_{g,n,d}$, 
\begin{eqnarray}
\lefteqn{ \left(\Osheaf_{\Pro(W_{g+1,n-1,d})}(1)\otimes\Osheaf_{\Pro^r}(1)\right)|_{C_0}
	= \left(\Osheaf_{\Pro(W_{g,n,d})}(1)\otimes\Osheaf_{\Pro^r}(1)\right)|_{C_0} }\notag\\
	&\hspace{1.5in}&\cong \acanon_{C_0}(ax_1+\cdots+ax_{n})\otimes\Osheaf_{\Pro^r}(ca+1)|_{C_0} \label{c_0_cong}\\
	&\hspace{1.5in}& = \left(\acanon_C(ax_1+\cdots+ax_{n-1})\otimes \Osheaf_{\Pro^r}(ca+1)|_{C}\right)|_{C_0}. \notag
\end{eqnarray}
To analyse $C_1$ we also observe that $\Osheaf_{\Pro^r}(1)|_{C_1}$ is trivial, as $f$ contracts $C_1$ to a point.
\begin{eqnarray}
\lefteqn{ \left(\Osheaf_{\Pro(W_{g+1,n-1,d})}(1)\otimes\Osheaf_{\Pro^r}(1)\right)|_{C_1}
	= \Osheaf_{\Pro(W_{1,1,0})}(1)|_{C_1} }\notag\\
	&\hspace{1.5in}&\cong \acanon_{C_1}(ay) \label{c_1_cong}\\
	&&		= \left(\acanon_C(ax_1+\cdots+ax_{n-1})\otimes\Osheaf_{\Pro^r}(ca+1)|_C\right)|_{C_1} .\notag
\end{eqnarray}
  The curve $C$ was defined as $C_0\cup C_1$.    We have an induced isomorphism of line bundles 
\begin{eqnarray*}
\lefteqn{ \left(\Osheaf_{\Pro(W_{g+1,n-1,d})}(1)\otimes\Osheaf_{\Pro^r}(1) \right)|_{C\setminus \{P\}} }\\
	&\hspace{1in}&\cong	\left(\acanon_C(ax_1+\cdots+ax_{n-1})\otimes \Osheaf_{\Pro^r}(ca+1)|_{C})\right)|_{C\setminus \{P\}}
\end{eqnarray*}	
\label{match_sheaves}found by excluding $P$ from the two isomorphisms above.  To extend this over $P$, we simply need to insist that the isomorphisms over $C_0$ and $C_1$ are consistent at $P$.
When we restrict to the fibre over $P$, the two isomorphisms (\ref{c_0_cong}) and (\ref{c_1_cong}) are scalar multiples of one another, so we obtain consistency at $P$ by multiplying the isomorphism (\ref{c_1_cong}) by a suitable non-zero scalar, once (\ref{c_0_cong}) is given.

Thus $(C,x_1,\ldots ,x_{n-1})$ is indeed represented in $J_{g+1,n-1,d}$, as required.  Now we use our inductive hypothesis: we know $\bar{J}_{g+1,n-1,d}^{ss}(L_{m,\mhat,m'})=J_{g+1,n-1,d}$, so in particular $(C,x_1,\ldots,x_{n-1})$ is semistable with respect to $L_{m,\mhat,m'}$.

For the following analysis we must clarify the notation for our standard line bundles. For $i=1,2$ let
\[
C_i \stackrel{\iota_{C_i,C}}{\rightarrow} C \stackrel{\iota_{C}}{\rightarrow} \Pro(W_{g+1,n-1,d})\times\Pro^r
\]
be the inclusion morphisms.  Note that the composition 
\[
C_0 \stackrel{\iota_{C_0,\Pro(W_{g,n,d})}}{\rightarrow} \Pro(W_{g,n,d})\times\Pro^r \stackrel{\iota_{\Pro(W_{g,n,d})}}{\rightarrow}\Pro(W_{g+1,n-1,d})\times\Pro^r
\]
is equal to $\iota_{C_0,C}\circ \iota_{C}$.  Thus,
\[
\iota_{C_0,\Pro(W_{g,n,d})}^{*}p_{W_{g,n,d}}^{*}\mathcal{O}_{\Pro(W_{g,n,d})}(1) 
	= \iota_{C_0,C}^{*}\iota_{C}^{*}p_{W_{g+1,n-1,d}}^{*}\mathcal{O}_{\Pro(W_{g+1,n-1,d})}(1).
\]
Therefore we may denote this line bundle by $L_{W\, C_0}$, as we need not specify which `$W$' space we have used.  $L_{W\, C_1}$ is defined similarly.  We let
\[\begin{array}{rclcrcl}
L_{r\, C_0} 	&:=& \iota_{C_0,C}^{*}\iota_{C}^{*}p_{r}^{*}\mathcal{O}_{\Pro^{r}}(1)
	&& L_{r\, C_1} 	&:=& \iota_{C_1,C}^{*}\iota_{C}^{*}p_{r}^{*}\mathcal{O}_{\Pro^{r}}(1),
\end{array}\]
though $L_{r\,C_1}$ is in fact trivial, since $C_1$ is collapsed by projection to $\Pro(r)$.

For $i=0,1$, define restriction maps to our subcurves:
\begin{eqnarray*}
\lefteqn{\phat^{C_i}_{m,\mhat}:H^0(\Pro(W_{g+1,n-1,d})\times\Pro^r,\Osheaf_{\Pro(W_{g+1,n-1,d})}(m)\otimes\Osheaf_{\Pro^r}(\mhat)) }\\
	&\hspace{3in}& \rightarrow H^0(C_i,L_{W\, C_i}^m\otimes L_{r\, C_i}^{\mhat}).
\end{eqnarray*}
We show these are surjective for $m$ and $\mhat$ sufficiently large.  The conditions on $m$ and $\mhat$ imply $m,\mhat > m_{3}$, so the restriction map
\begin{eqnarray*}
\lefteqn{ \phat^{C_0,\Pro(W_{g,n,d})}_{m,\mhat}:H^0(\Pro(W_{g,n,d})\times\Pro^r,\Osheaf_{\Pro(W_{g,n,d})}(m)\otimes\Osheaf_{\Pro^r}(\mhat)) }\\
	&\hspace{3in}& \rightarrow H^0(C_0,L_{W\, C_0}^{m} \otimes L_{r\, C_0}^{\mhat})
\end{eqnarray*}
is surjective by Grothendieck's uniform $m$ lemma (cf. Proposition \ref{7facts}(i)).  The restriction to a linear subspace
\begin{eqnarray*}
\lefteqn{\phat^{\Pro(W_{g,n,d})}_{m,\mhat}:H^0(\Pro(W_{g+1,n-1,d})\times\Pro^r,\Osheaf_{\Pro(W_{g+1,n-1,d})}(m)\otimes\Osheaf_{\Pro^r}(\mhat))}\\
	&\hspace{1.5in}&\rightarrow H^0(\Pro(W_{g,n,d})\times\Pro^r,\Osheaf_{\Pro(W_{g,n,d})}(m)\otimes\Osheaf_{\Pro^r}(\mhat))
\end{eqnarray*}
is surjective, hence by composition $\phat^{C_0}_{m,\mhat}= \phat^{C_0,\Pro(W_{g,n,d})}_{m,\mhat}\circ\phat^{\Pro(W_{g,n,d})}_{m,\mhat} $ is surjective.  Surjectivity of $\phat^{C_1}_{m,\mhat}$ is shown similarly.
Moreover, for $i=0,1$, the restriction $\phat^{C_i}_{m,\mhat}$ factors through maps which restrict sections of $C$ to those on one of the subcurves:
\[
\phat^{C_i,C}_{m,\mhat}:H^0(C,L_{W\, C}^m\otimes L_{r\, C}^{\mhat})
	\rightarrow H^0(C_i,L_{W\, C_i}^m\otimes L_{r\, C_i}^{\mhat}),
\]
and thus the maps $\phat^{C_0,C}_{m,\mhat}$ and $\phat^{C_1,C}_{m,\mhat}$ are also surjective.

To relate semistability of $(C,x_1,\ldots,x_{n-1})$ to that of $(C_0,\pts )$, we shall display the vector space $H^0(C,L_{W\, C}^m\otimes L_{r\, C}^{\mhat})$ as the direct sum of two subspaces.  Recall that we wrote
\[
W_{g+1,n-1,d}=W_{g,n,d}\oplus V_{1,1,0}.
\]
 Fix a basis $w_0,\ldots w_{N_{g+1,n-1,d}}$ for $W_{g+1,n-1,d}$ respecting this decomposition.  Let $\hat{B}_{m,\mhat}$ be a basis of $H^0(\Pro(W_{g+1,n-1,d})\otimes\Pro^r,\Osheaf_{\Pro(W_{g+1,n-1,d})}(m)\otimes\Osheaf_{\Pro^r}(\mhat))$ of monomials of bidegree $(m,\mhat)$, where the degree $m$ part is a monomial in $w_0,\ldots w_{N_{g+1,n-1,d}}$.

Let $\Omm_{+}$ be the subspace of $H^0(\Pro(W_{g+1,n-1,d})\times\pr,\Osheaf_{\Pro(W_{g+1,n-1,d})}(m)\otimes\opr(\mhat))$ spanned by all monomials in $\hat{B}_{m,\mhat}$ which have at least one factor from $V_{1,0,0}$.  Namely, $\Omm_+$ is spanned by monomials of bidegree $(m,\mhat)$, where the degree $m$ part has form $w_{i_1}w_{i_2}\cdots w_{i_m}$, with $w_{i_1}\in V_{1,1,0}$.  The remaining factors may come from either $W_{g,n,d}$ or $V_{1,0,0}$.

Similarly, let $\Omm_0\subset H^0(\Pro(W_{g+1,n-1,d})\times\pr,\Osheaf_{\Pro(W_{g+1,n-1,d})}(m)\otimes\opr(\mhat))$ be the subspace spanned by monomials in $\hat{B}_{m,\mhat}$ which have no factors from $V_{1,0,0}$.  In other words, all of the factors in the degree $m$ part come from $W_{g,n,d}$.

By inspecting the basis $\hat{B}_{m,\mhat}$, we see that as vector spaces
\begin{equation}
H^0(\Pro(W_{g+1,n-1,d})\times\pr,\Osheaf_{\Pro(W_{g+1,n-1,d})}(m)\otimes\opr(\mhat)) = \Omm_+ \oplus \Omm_0. \label{splits on the proj space}
\end{equation}
We wish to show that this decomposition restricts to the curve $C$.  Set
\begin{eqnarray*}
\barO_+		&:=&\phat^C_{m,\mhat}(\Omm_+)\subset H^0(C,L_{W\,C}^m\otimes L_{r\, C}^{\mhat})\\
\barO_0		&:=&\phat^C_{m,\mhat}(\Omm_0)\subset H^0(C,L_{W\,C}^m\otimes L_{r\, C}^{\mhat}).
\end{eqnarray*}

We first make some technical observations on where such sections vanish.


\begin{claim}\label{technical} Let $\barO_+$ and $\barO_0$ be defined as above.
\begin{enumerate}
\item[(i)] If $0\neq s\in \barO_+$ then $\phat^{C_0,C}_{m,\mhat}(s)=0$ and $\phat^{C_1,C}_{m,\mhat}(s)\neq 0$.
\item[(ii)] The image $\phat^{C_1,C}_{m,\mhat}(\barO_0)$ is one-dimensional, and if $s\in \barO_0$ with $s(P)=0$, then $\phat^{C_1,C}_{m,\mhat}(s)= 0$.
\end{enumerate}
\end{claim}
\begin{proof}[Proof of Claim \ref{technical}.] 
(i) Recall we observed that, for any $w\in V_{1,1,0}$, the restriction $\rho^{C_0}(w)=0$. The space $\Omm_+$ is spanned by monomials containing a factor from $V_{0,1,1}$, so it follows that if $s\in \barO_+$ then $\phat^{C_0,C}_{m,\mhat}(s)=0$.  However, if $s\in \barO_+\subset H^0(C,L_{W\, C}^m\otimes L_{r\, C}^{\mhat})$ and $s\neq 0$ then $s$ must be non-zero on one of $C_0$ and $C_1$, whence $\phat^{C_1,C}_{m,\mhat}(s)\neq 0$.

(ii) Recall that we wrote $W_{g,n,d}\cong V_{g,n,d}\oplus U$, where $V_{g,n,d}$ consists of sections vanishing at $p_W(P)$, and $U$ is spanned by a section non-vanishing at $p_W(P)$.  We saw that $\rho^{C_1}(V_{g,n,d})=\{0\}$, and so $\rho^{C_1}(W_{g,n,d})=\rho^{C_1}(U)$, which is one-dimensional.  Let $u$ span this space; $u(p_W(P))\neq 0$.

The component $C_1$ collapses to the single point $p_r(P)$ under $p_r$.
Thus we have another decomposition: 
$H^0(\Pro^r,\Osheaf_{\Pro^r}(1))=V_{\Pro^r}\oplus \hat{U}$, where sections in $V_{\Pro^r}$ vanish at $p_r(P)=p_r(C_1)$, and $\hat{U}$ is one-dimensional, spanned by a section non-vanishing at $p_r(P)$.  Again, $\rho^{C_1}(H^0(\Pro^r,\Osheaf_{\Pro^r}(1)))=\rho^{C_1}(\hat{U})$, which is one-dimensional.  Let $\hat{u}$ span this space; $\hat{u}(p_r(P))\neq 0$.

Thus, since $\Omm_+=\mbox{Sym}^m(W_{g,n,d}\otimes H^0(\pr,\opr(\mhat)))$, we have 
\[
\phat^{C_1,C}_{m,\mhat}(\barO_+)=\phat^{C_1}_{m,\mhat}(\Omm_+)=\phat^{C_1}_{m,\mhat}(\mbox{Sym}^mU\otimes \mbox{Sym}^{\mhat}\hat{U}), 
\]
which is one-dimensional since $U$ and $\hat{U}$ are. Finally, if $s\in \barO_0$ then $\phat^{C_1,C}_{m,\mhat}(s)=\alpha u^m\otimes\hat{u}^{\mhat}$ for some $\alpha\in k$.  We know that $u^m\otimes\hat{u}^{\mhat}(P)\neq 0$.  Since $s(P)=\phat^{C_1,C}_{m,\mhat}(s)(P)$, it follows that if $s(P)=0$ then $\alpha=0$ and so $\phat^{C_1,C}_{m,\mhat}(s)=0$, and the proof of Claim \ref{technical} is complete.
\end{proof}
We now give details of our decomposition of $H^0(C,L_{W\, C}^m\otimes L_{r\,C}^{\mhat})$.  
\begin{claim}\label{decomp}
Restriction to $C$ respects the decomposition (\ref{splits on the proj space}), and enables us to identify the spaces $\barO_+$ and $\barO_0$.  Precisely:
\begin{list}
{}{\setlength{\leftmargin}{.75in} \setlength{\rightmargin}{.75in} 
\setlength{\itemindent}{-.25in}}
\item[(i)] $H^0(C,L_{W\, C}^m\otimes L_{r\, C}^{\mhat})
	=\barO_0\oplus\barO_+$;
\item[(ii)] $
		\barO_0 \cong H^0(C_0,L_{W\, C_0}^m\otimes L_{r\, C_0}^{\mhat})$;
\item[(iii)] $
		\barO_+ \cong H^0(C_1,L_{W\, C_1}^m(-P))$.
\end{list}
\end{claim}
\begin{proof}[Proof of Claim \ref{decomp}.] 
(i)
By restricting (\ref{splits on the proj space}) to $C$, we see:
\[
H^0(C,L_{W\, C}^m\otimes L_{r\, C}^{\mhat})
	=\barO_0 + \barO_+.
\]
It remains to show that these spaces have zero intersection.  

Suppose $s\in \barO_0 \cap \barO_+$.  Since $s\in \barO_+$, it follows by \ref{technical}(i) that $\phat^{C_0,C}_{m,\mhat}(s)= 0$; in particular, as $P\in C_0$, we see $s(P)=0$.  We know in addition that $s\in \barO_0$, so \ref{technical}(ii) implies that $\phat^{C_1,C}_{m,\mhat}(s)=0$.  Thus $s$ restricts to zero on the whole of $C$, and we conclude that $s=0$, proving (i).

(ii) 
Recall that the morphism 
\[
\phat^{C_0,C}_{m,\mhat}:H^0(C,L_{W\, C}^m\otimes L_{r\, C}^{\mhat})\rightarrow H^0(C_0,L_{W\, C_0}^m\otimes L_{r\, C_0}^{\mhat})
\]
is surjective.  However, if $s\in\barO_+$ then by Claim \ref{technical}(i) we know that $\phat^{C_0,C}_{m,\mhat}(s)=0$.  Thus, since 
$H^0(C,L_{W\, C}^m\otimes L_{r\, C}^{\mhat})\cong \barO_0\oplus\barO_+$ we conclude that $\phat^{C,C_0}_{m,\mhat}|_{\barO_0}$ is surjective.

On the other hand, 
\[
\phat^{C_0,C}_{m,\mhat}|_{\barO_0}:\barO_0
		\rightarrow H^0(C_0,L_{W\,C_0}^m\otimes L_{r\,C_0}^{\mhat})
\]
 is injective.  For if $s \in \barO_0$ and $\phat^{C_0,C}_{m,\mhat}(s)=0$, then $s(P)=0$; then by \ref{technical}(ii) it follows that $\phat^{C_1,C}_{m,\mhat}(s)=0$.  Thus $s=0$, so $\phat^{C_0,C}_{m,\mhat}|_{\barO_0}$ has zero kernel.  We conclude that $\phat^{C_0,C}_{m,\mhat}|_{\barO_0}$ is an isomorphism of vector spaces.

(iii)
To start with, note that if $M\in\Omm_+$ then $M$ has at least one factor from $V_{0,1,1}$, so $M$ vanishes at $P$ by definition; hence  
\[
\phat^{C_1,C}_{m,\mhat}(\barO_+)
	\subset H^0(C_1,L_{W\,C_1}^m(-P)).
\]
Moreover the map $\phat^{C_1,C}_{m,\mhat}|_{\barO_+}$ is injective.  For if $s\in \barO_+$, then $\phat^{C_0,C}_{m,\mhat}(s)=0$ and so if in addition $\phat^{C_1,C}_{m,\mhat}(s)= 0$ then we must conclude that $s=0$.

On the other hand, the unrestricted morphism 
\[
\phat^{C_1,C}_{m,\mhat}:\barO_+\oplus\barO_0\rightarrow H^0(C_1,L_{W\, C_1}^m)
\]
is surjective.  In particular the subspace $H^0(C_1,L_{W\, C_1}^m(-P))\subset H^0(C_1,L_{W\, C_1}^m)$ is in its image.  So for arbitrary $s\in H^0(C_1,L_{W\, C_1}^m(-P))$ we may write $s=\phat^{C_1,C}_{m,\mhat}(s_0+s_+)$, where $s_0\in\barO_0$ and $s_+\in\barO_+$.  But then 
\[
\phat^{C_1,C}_{m,\mhat}(s_0)=s-\phat^{C_1,C}_{m,\mhat}(s_+)\in H^0(C_1,L_{W\, C_1}^m(-P)),
\]
so $s_0$ vanishes at $P$.  However, by \ref{technical}(ii), it follows that $\phat^{C_1,C}_{m,\mhat}(s_0)=0$, and so we conclude $s=\phat^{C_1,C}_{m,\mhat}(s_+)$.  In other words, the morphism $\phat^{C_1,C}_{m,\mhat}|_{\barO_+}$ is also surjective, showing:
\[
\phat^{C_1,C}_{m,\mhat}|_{\barO_+}:\barO_+
	\cong H^0(C_1,L_{W\,C_1}^m(-P)),
\]
which completes the proof of Claim \ref{decomp}.
\end{proof}

{\it Back to the proof of Proposition \ref{inductive_step}.}
We wish to show that $(C_0,\pts)$ is GIT semistable, with respect to some $l\in\mathbf{H}_{M_{g,n,d}}(\bar{J}_{g,n,d})$.  Therefore let $\lambda'$ be a 1-PS of $SL(W_{g,n,d})$ acting on $(C_0,x_1,\ldots ,x_n)\subset \Pro(W_{g,n,d})\times\Pro^r$.  Let $w_0,\ldots,w_{N_{g,n,d}}$ be a basis for $W_{g,n,d}$ diagonalising the action of $\lambda'$, so that $\lambda'$ acts with weights $r_0,\ldots ,r_{N_{g,n,d}}$. Let $r_{i_j}$ be the minimal $\lambda'$-weight at $x_j$, so
\[
r_{i_j}:=\min\{r_i|w_i(x_j)\neq 0\}.
\]
We wish to show that 
\[
\mu^{L_{m,\mhat,m''}}((C_0,\pts),\lambda')=\mu^{L_{m,\mhat}}(C,\lambda')+m''\sum_{j=1}^n r_{i_j} \leq 0
\]
for some $m''$ such that $(m,\mhat,m'')\in M_{g,n,d}$.

We have identified $W_{g+1,n-1,d}$ with $W_{g,n,d}\oplus V_{1,1,0}$.  Extend $w_0,\ldots,w_{N_{g,n,d}}$ to a basis $w_0,\ldots,w_{N_{g+1,n-1,d}}$ for $W_{g+1,n-1,d}$ respecting this decomposition. We define $\lambda$ to be the extension of $\lambda'$ over $W_{g+1,n-1,d}$, which acts with weight $r_{i_n}$ on all of $V_{1,1,0}$.  Note that $\lambda$ is a subgroup of $GL(W_{g+1,n-1,d})$, but not necessarily of $SL(W_{g+1,n-1,d})$, so we shall have to use the formula in Lemma \ref{numcritlemma}.  We already calculated $\dim V_{1,1,0}=(a-1)$ so we know the sum of the weights: 
\[
\sum_{i=0}^{N_{g+1,n-1,d}}w_{\lambda}(w_i)=(a-1)r_{i_n}.
\]

Let $\hat{B}_{m,\mhat}$ be a basis of $H^0(\Pro(W_{g+1,n-1,d})\times\Pro^r,\Osheaf_{\Pro(W_{g+1,n-1,d})}(m)\otimes\Osheaf_{\Pro^r}(\mhat))$ consisting of monomials of bidegree $(m,\mhat)$ such that the degree $m$ part is a monomial in $w_0,\ldots,w_{N_{g+1,n-1,d}}$.  Semistability of $(C,x_1,\ldots ,x_{n-1})$ implies that there exists a subset $\mathcal{B}\subset\hat{B}_{m,\mhat}$, such that
\[
\overline{\mathcal{B}}:=\{\phat_{m,\mhat}^C(M):M\in \mathcal{B} \}
\]
is a basis for $H^0(C,L_{W\,C}^m\otimes L_{r\,C}^{\mhat})$, with
\begin{eqnarray}
\lefteqn{\Big(\sum_{M\in\mathcal{B}} w_{\lambda}(M) + m'\sum_{j=1}^{n-1}r_{i_j}\Big)(e_{g+1,n-1,d}-g) }\notag\\
	&&\hspace{.5in}{}- \left(m(e_{g+1,n-1,d}m+d\mhat-g)+(n-1)m'\right)(a-1)r_{i_n} \leq 0. \hspace{.3in} \label{C_num_crit}
\end{eqnarray}

We relate this to a weight for $(C_0,\pts)\in\Pro(W_{g,n,d})\times\Pro^r$.  The basis $\hat{B}_{m,\mhat}$ respects the decomposition 
\[
H^0(\Pro(W_{g+1,n-1,d})\otimes\Pro^r,\Osheaf_{\Pro(W_{g+1,n-1,d})}(m)\otimes\Osheaf_{\Pro^r}(\mhat))
	=\Omm_0\oplus\Omm_+;
\]
hence the basis $\overline{\mathcal{B}}$ respects the decomposition 
\[
H^0(C,L_{W\, C}^m\otimes L_{r\, C}^{\mhat})
	=\barO_0\oplus\barO_+
\]
of Claim \ref{decomp}(i).  We split $\overline{\mathcal{B}}$ into two parts, $\overline{\mathcal{B}}_0$ and $\overline{\mathcal{B}}_+$, following this decomposition.  We will calculate the contribution to the weight coming from the corresponding subsets $\mathcal{B}_0$ and $\mathcal{B}_+$ of $\mathcal{B}$.

First consider $\mathcal{B}_0$.  We show the weight of this collection of monomials is a $\lambda'$-weight for $C_0$.  
Recall that the map $\phat^{\Pro(W_{g,n,d})}_{m,\mhat}$ which restricts sections to a linear subspace is surjective; this map is zero on $\Omm_+$, so 
\[
\phat^{\Pro(W_{g,n,d})}_{m,\mhat}|_{\Omm_0}:\Omm_0\rightarrow H^0(\Pro(W_{g,n,d})\times\Pro^r,\Osheaf_{\Pro(W_{g,n,d})}(m)\otimes\Osheaf_{\Pro^r}(\mhat))
\]
is an isomorphism of vector spaces (note we are speaking of the spaces without a bar, those that have not been restricted to $C$).  Moreover, if $M\in\Omm_0$ is in $\mathcal{B}_0$, then $M$ is a monomial whose degree $m$ part is in the basis $w_0,\ldots ,w_{N_{g,n,d}}$ for $W_{g,n,d}$; we may interpret $\phat^{\Pro(W_{g,n,d})}(M)$ as the same monomial in this basis.  Then, since the actions of $\lambda'$ and $\lambda$ are identical on $W_{g,n,d}$,
\begin{equation}\label{equal_weights}
w_{\lambda'}(\phat^{\Pro(W_{g,n,d})}_{m,\mhat}(M))=w_{\lambda}(M).
\end{equation}

Now $\{\phat^{\Pro(W_{g,n,d})}(M)|M\in\mathcal{B}_0\}$ is a collection of monomials of bidegree $(m,\mhat)$ in $H^0(\Pro(W_{g,n,d})\times\Pro^r,\Osheaf_{\Pro(W_{g,n,d})}(m)\otimes\Osheaf_{\Pro^r}(\mhat))$ such that when we restrict them to $C_0$, we obtain 
\[
\{\phat^{C_0,\Pro(W_{g,n,d})}_{m,\mhat}\circ \phat^{\Pro(W_{g,n,d})}_{m,\mhat}(M) | M\in\mathcal{B}_0\}
	=\{\phat^{C_0,C}_{m,\mhat}\circ\phat^C_{m,\mhat}(M)|\phat^C_{m,\mhat}(M)\in\overline{\mathcal{B}}_0\}, 
\]
which by Claim \ref{decomp}(ii) is a basis for $H^0(C_0,L_{W\,C_0}^m\otimes L_{r\,C_0}^{\mhat})$.  It follows that the restricted monomials, $\{\phat^{\Pro(W_{g,n,d})}(M)|M\in\mathcal{B}_0\}$, give a $\lambda'$-weight of $C_0$, i.e.
\[
\mu^{L_{m,\mhat}}(C,\lambda')
	\leq \sum_{M\in\mathcal{B}_0}w_{\lambda'}(\phat^{\Pro(W_{g,n,d})}(M))
	= \sum_{M\in\mathcal{B}_0}w_{\lambda}(M),
\]
where we have used (\ref{equal_weights}) to relate this back to $\lambda$.

Now consider $\mathcal{B}_+\subset\barO_+$. If $M\in\mathcal{B}_+$ then by \ref{technical}(i) the restriction  $\phat^{C_0}_{m,\mhat}(M)= 0$; basis elements are non-zero so $\phat^{C_1}_{m,\mhat}(M)\neq 0$.  The monomial $M$ has bidegree $(m,\mhat)$, where the degree $m$ part is a monomial in the basis $w_0,\ldots, w_{N_{g+1,n-1,d}}$, which respects the decomposition $W_{g+1,n-1,d}=W_{g,n,d}\oplus V_{1,1,0}$.  Any factor in $M$ from $W_{g,n,d}$ must be non-zero on $C_1$, so by the proof of \ref{technical}(ii) these are non-zero at $P$.  Hence the $\lambda$-weight of such factors is bounded below by $r_{i_n}$.  Meanwhile, the remaining factors in $M$ come from $V_{1,1,0}$, and all have $\lambda$-weight $r_{i_n}$.  We conclude that for such $M$,
\[
w_{\lambda}(M)\geq r_{i_n}m.
\]
We count the number of such $M$ using Claim \ref{decomp}(ii):
\[
\# \mathcal{B}_+=\#\overline{\mathcal{B}}_+=h^0(C_1,L_{W\, C_1}^m(-P))=am-1.
\]
Thus
\[
\sum_{M\in\mathcal{B}_+} w_{\lambda'}(M)\geq (am^2-m)r_{i_n}.
\]

We may now insert these estimates into inequality (\ref{C_num_crit}):
\begin{eqnarray*}
\lefteqn{\left(\mu^{L_{m,\mhat}}(C_0,\lambda') + (am^2-m)r_{i_n}+ m'(r_{i_1}+\cdots + r_{i_{n-1}}\right)(e_{g+1,n-1,d}-g) } \\
	&\hspace{1in}&
	\leq (m(e_{g+1,n-1,d}m+d\mhat -g)+(n-1)m')(a-1)r_{i_n} \hspace{.5in}\\
	&&
	=(e_{g+1,n+1}+d\frac{\mhat}{m}+(n-1)\frac{m'}{m^2}-\frac{g}{m})(a-1)m^2r_{i_n}.
\end{eqnarray*}
Recall that we set $\frac{\mhat}{m}=\frac{ca}{2a-1}$ and $\frac{m'}{m^2}=\frac{a}{2a-1}$.  Further, one may expand: $e_{g+1,n-1,d}-g=(2a-1)g+a(n-1+cd)$.  Thus we have shown:

\begin{eqnarray*}
\lefteqn{\mu^{L_{m,\mhat}}(C_0,\lambda') + (am^2-m)r_{i_n}+ \frac{am^2}{2a-1}(r_{i_1}+\cdots + r_{i_{n-1}})}\\
	&\hspace{1in}&
	\leq \left(1+\frac{g+d\frac{ca}{2a-1}+(n-1)\frac{a}{2a-1}-\frac{g}{m}}{(2a-1)g+a(n-1+cd)}\right)(a-1)m^2r_{i_n} \\
	&&
	= \left(1+\frac{1}{2a-1}\right)(a-1)m^2r_{i_n}-S_{14}mr_{i_n},
\end{eqnarray*}

\noindent where $S_{14}:=\frac{g(a-1)}{(2a-1)g+a(n-1+cd)}<\frac{1}{2}$.
In conclusion:
\begin{equation}
\mu^{L_{m,\mhat}}(C_0,\lambda') + \frac{am^2}{2a-1}(r_{i_1}+\cdots + r_{i_n})
	-(1-S_{14})mr_{i_n} 
		\leq 0.
\label{result!}
\end{equation}

We may repeat the argument above, attaching the elliptic curve at the location of any other of the points $x_1,\ldots ,x_{n-1}$.  This will give us similar inequalities, but with the role $x_n$ played by a different marked point.  Adding up all such inequalities, and dividing by $n$, yields:
\[
\mu^{L_{m,\mhat}}(C,\lambda')+\left(\frac{am^2}{2a-1}-\frac{1-S_{14}}{n}m\right)\sum_{j=1}^n r_{i_j}\leq 0.
\]
By the hypotheses on $m$, we know that $\frac{1-S_{14}}{n}m$ is an integer, so if we set $m''=\frac{am^2}{2a-1}-\frac{1-S_{14}}{n}m$ then this is also an integer, and we have shown that
\[
\mu^{L_{m,\mhat,m''}}((C,\pts),\lambda')\leq 0.
\]
Thus, as we made no assumptions about the 1-PS $\lambda'$, 
\[
(C_0,x_1,\ldots ,x_n)\in \bar{J}^{ss}(L_{m,\mhat,m''}).
\]  

It remains to show that $(m,\mhat,m'')\in M_{g,n,d}$, where $M_{g,n,d}$ is as in the statement of this proposition.  We observed that $M_{g+1,n-1,d}\subset M_{g,n,d}$ and so the conditions on $m$ and $\mhat$ are satisfied.   We need to check that $n|\frac{m''}{m^2}-\frac{a}{2a-1}|=n\frac{1-S_{14}}{nm}\leq \frac{1}{4a-2}-\frac{3g+q_2-\bar{g}-\frac{1}{2}}{m}$.  But $\frac{1-S_{14}}{m}<\frac{1}{m}$, and
\begin{eqnarray*}
\frac{1}{m}	&<&\frac{1}{4a-2}-\frac{3g+q_2-\bar{g}-\frac{1}{2}}{m}\\
\Longleftrightarrow m &>& (3g+q_2-\bar{g}+\frac{1}{2})(4a-2).
\end{eqnarray*}
This is implied by $(m,\mhat,m')\in M_{g+1,n-1,d}$, for the final lower bound on $m$ is $m>(6(g+1)+2q_2-2(\overline{g+1})+1)(2a-1)$.

The choice of $(h,\pts)\in J$ was arbitrary.  Hence $J_{g,n,d}\subset\bar{J}_{g,n,d}^{ss}(L_{m,\mhat,m''})$.  We proved the reverse inclusion in Theorem \ref{jssclosed}, and so  $\bar{J}_{g,n,d}^{ss}(L_{m,\mhat,m''})=J_{g,n,d}$.  Now, by Corollary \ref{nearly_done} it follows that $\bar{J}_{g,n,d}^{ss}(l)=J_{g,n,d}$ for all $l\in\mathbf{H}_{M_{g,n,d}}(\bar{J}_{g,n,d})$, and in particular for all $L_{m,\mhat,m'}$ such that $(m,\mhat,m')\in M_{g,n,d}$.  This completes the proof of Proposition \ref{inductive_step}.
\end{proof}

Now our desired results are immediate corollaries.  We no longer need to distinguish between different spaces of maps, so we stop using the subscripts for $J$ and $M$.
\begin{proof}[Proof of Theorem \ref{maps}, Corollary \ref{general_x_statement} and Theorem \ref{curves}.]
For Theorem \ref{maps} and Corollary \ref{general_x_statement} we work over $\C$.  By Proposition \ref{base_case_maps} and Proposition \ref{inductive_step} we see that $\bar{J}^{ss}(l)=\bar{J}^{s}(l)=J$,  as schemes over $\C$, for any $l\in\mathbf{H}_M(\bar{J})$.  In particular this hold for any $L_{m,\mhat,m'}$ where $(m,\mhat,m')\in M$.  Then by Corollary \ref{catquodone} it follows that 
\[
\bar{J}\dblq_{L_{m,\mhat,m'}}SL(W)\cong\Maps,
\]
and by Corollary \ref{general_x} it follows that there exists $J_{X,\beta}$ such that
\[
\bar{J}_{X,\beta}\dblq_{L_{m,\mhat,m'}|_{\bar{J}_{X,\beta}}} SL(W)\cong \overline{\mathcal{M}}_{g,n,d}(X,\beta).
\]

For Theorem \ref{curves} we note that Theorem \ref{ind_step} has for convenience been phrased in terms of vector spaces over a field $k$, and so we work at first over $k$.  Now $\bar{J}^{ss}(l)=\bar{J}^{s}(l)= J$, for schemes over $k$, follows from Proposition \ref{base_case_curves} and Proposition \ref{inductive_step}, where $k$ is any field.  

It remains to check that the open subschemes $\bar{J}^{ss}(l)\subset\bar{J}$ and $J\subset \bar{J}$ are equal over $\Z$.  As they are indeed open subschemes, it is sufficient to check that they contain the same points, since they will then automatically have the same scheme structure.  If $x$ is a point in $J$, then we may consider $x$ as a point in $J\times_{\Z}\Spec k(x)$, where $k(x)$ is the residue field of $x$.  We know that $J$ and $\bar{J}^{ss}(l)$ are equal over $\Spec k(x)$, so this provides $J\subset\bar{J}^{ss}(l)$.  The converse is shown in the same way. 

Corollary \ref{catquodone} now implies again that in this case,
\[
\bar{J}\dblq_{L_{m,m'}}SL(W)\cong \curves
\]
over $\Z$.
\end{proof}




\begin{thebibliography}{99}

\bibitem{e_thesis}
E.~Baldwin.
\newblock {\em A Geometric Invariant Theory Construction of Moduli Spaces of
  Stable Maps}.
\newblock DPhil thesis, University of Oxford, 2006.

\bibitem{e_paper_2}
E.~Baldwin.
\newblock A {GIT} construction of moduli spaces of stable maps in positive
  characteristic.
\newblock \ math.AG/0707.2050, preprint 2007.

\bibitem{second_quotient}
E.~Baldwin.
\newblock A {GIT} construction of {$\curves$} as a quotient of a subscheme of
  {$\Maps$}.
\newblock in preparation.

\bibitem{BM}
K.~Behrend and Yu. Manin.
\newblock Stacks of stable maps and {G}romov-{W}itten invariants.
\newblock {\em Duke Math. J.}, 85(1):1--60, 1996.

\bibitem{DH}
I.~V. Dolgachev and Yi~Hu.
\newblock Variation of geometric invariant theory quotients.
\newblock {\em Inst. Hautes \'Etudes Sci. Publ. Math.}, (87):5--56, 1998.
\newblock \ With an appendix by Nicolas Ressayre.

\bibitem{FP}
W.~Fulton and R.~Pandharipande.
\newblock Notes on stable maps and quantum cohomology.
\newblock In {\em Algebraic geometry---Santa Cruz 1995}, volume~62 of {\em
  Proc. Sympos. Pure Math.}, pages 45--96. Amer. Math. Soc., Providence, RI,
  1997.

\bibitem{G}
D.~Gieseker.
\newblock {\em Lectures on moduli of curves}, volume~69 of {\em Tata Institute
  of Fundamental Research Lectures on Mathematics and Physics}.
\newblock Published for the Tata Institute of Fundamental Research, Bombay,
  1982.

\bibitem{EGA}
A.~Grothendieck.
\newblock \'{E}l\'ements de g\'eom\'etrie alg\'ebrique. {III}. \'{E}tude
  cohomologique des faisceaux coh\'erents. {II}.
\newblock {\em Inst. Hautes \'Etudes Sci. Publ. Math.}, (17):91, 1963.

\bibitem{HM}
J.~Harris and I.~Morrison.
\newblock {\em Moduli of curves}, volume 187 of {\em Graduate Texts in
  Mathematics}.
\newblock Springer-Verlag, New York, 1998.

\bibitem{Hart}
R.~Hartshorne.
\newblock {\em Algebraic geometry}.
\newblock Springer-Verlag, New York, 1977.
\newblock \ Graduate Texts in Mathematics, No. 52.

\bibitem{hassett-hyeon}
B.~Hassett and D.~Hyeon.
\newblock Log canonical models for the moduli space of curves: {F}irst divisorial contraction.
\newblock \ math.AG/0607477, preprint 2006.

\bibitem{hyeon-lee}
D.~Hyeon and Y.~Lee.
\newblock Log minimal model program for the moduli space of stable curves of genus three.
\newblock \ math.AG/0703093, preprint 2007.

\bibitem{KP}
B.~Kim and R.~Pandharipande.
\newblock The connectedness of the moduli space of maps to homogeneous spaces.
\newblock In {\em Symplectic geometry and mirror symmetry (Seoul, 2000)}, pages
  187--201. World Sci. Publ., River Edge, NJ, 2001.

\bibitem{knudsen3}
F.~F.~Knudsen.
\newblock The projectivity of the moduli space of stable curves. {III}. {T}he
  line bundles on {$M\sb{g,n}$}, and a proof of the projectivity of {$\overline
  M\sb{g,n}$} in characteristic {$0$}.
\newblock {\em Math. Scand.}, 52(2):200--212, 1983.

\bibitem{morrison}
I.~Morrison.
\newblock Projective stability of ruled surfaces.
\newblock {\em Invent. Math.}, 56(3):269--304, 1980.

\bibitem{Mukai}
S.~Mukai.
\newblock {\em An introduction to invariants and moduli}, volume~81 of {\em
  Cambridge Studies in Advanced Mathematics}.
\newblock Cambridge University Press, Cambridge, 2003.
\newblock \ Translated from the 1998 and 2000 Japanese editions by W. M.
  Oxbury.

\bibitem{GIT}
D.~Mumford, J.~Fogarty, and F.~Kirwan.
\newblock {\em Geometric invariant theory}, volume~34 of {\em Ergebnisse der
  Mathematik und ihrer Grenzgebiete (2) [Results in Mathematics and Related
  Areas (2)]}.
\newblock Springer-Verlag, Berlin, third edition, 1994.

\bibitem{New}
P.~E. Newstead.
\newblock {\em Introduction to moduli problems and orbit spaces}, volume~51 of
  {\em Tata Institute of Fundamental Research Lectures on Mathematics and
  Physics}.
\newblock Tata Institute of Fundamental Research, Bombay, 1978.

\bibitem{parker}
A.~Parker.
\newblock An elementary {GIT} construction of the moduli space of stable maps.
\newblock \ math.AG/0604092, preprint 2006.

\bibitem{Schubert}
D.~Schubert.
\newblock A new compactification of the moduli space of curves.
\newblock {\em Compositio Math.}, 78(3):297--313, 1991.

\bibitem{sesh_completeness}
C.~S. Seshadri.
\newblock Quotient spaces modulo reductive algebraic groups.
\newblock {\em Ann. of Math. (2)}, 95:511--556; errata, ibid. (2) 96 (1972),
  599, 1972.

\bibitem{sesh}
C.~S. Seshadri.
\newblock Geometric reductivity over arbitrary base.
\newblock {\em Advances in Math.}, 26(3):225--274, 1977.

\bibitem{david_paper_2}
D.~Swinarski.
\newblock Geometric invariant theory stability for weighted pointed stable curves.
\newblock in preparation.

\bibitem{Swin}
D.~Swinarski.
\newblock Geometric invariant theory and moduli spaces of maps.
\newblock Master's thesis, University of Oxford, 2003.

\bibitem{Th}
M.~Thaddeus.
\newblock Geometric invariant theory and flips.
\newblock {\em J. Amer. Math. Soc.}, 9(3):691--723, 1996.

\bibitem{vakil_Mgn}
R.~Vakil.
\newblock The moduli space of curves and its tautological ring.
\newblock {\em Notices Amer. Math. Soc.}, 50(6):647--658, 2003.

\bibitem{murphy}
R.~Vakil.
\newblock Murphy's law in algebraic geometry: badly-behaved deformation spaces.
\newblock {\em Invent. Math.}, 164(3):569--590, 2006.

\end{thebibliography}
\end{document}